\documentclass{article}
\usepackage{graphicx}
\usepackage{epsfig}
\usepackage{amssymb, amsmath}
\usepackage{amsthm}
\usepackage{setspace}
\usepackage[top = 1in, bottom = 1in, left = 1.2in, right =
1.2in]{geometry}
\newtheorem{theorem}{Theorem}[section]

\newtheorem{lemma}[theorem]{Lemma}
\newtheorem{prop}[theorem]{Proposition}

\newtheorem{definition}[theorem]{Definition}
\newcommand{\li}{\int_{x_0}^1}

\begin{document}

\title{\LARGE Local Well-Posedness of Dynamics of Viscous Gaseous Stars}
\author{Juhi Jang \\ \it Brown University\\ 
 juhijang@math.brown.edu }

\singlespacing \maketitle \numberwithin{equation}{section}

\begin{abstract}
We establish the local in time well-posedness of strong solutions to
the vacuum free boundary problem of the compressible
Navier-Stokes-Poisson system in the spherically symmetric and
isentropic motion. Our result captures the physical vacuum boundary
behavior of the Lane-Emden star configurations for all adiabatic
exponents $\gamma>\frac{6}{5}$.
\end{abstract}\

\section{Formulation and Notation}
\label{section1}

The motion of self-gravitating viscous gaseous stars can be
described by the compressible Navier-Stokes-Poisson system:
\begin{equation}
\begin{split}
\frac{\partial\rho}{\partial t} + \nabla\cdot(\rho \mathbf{u})&=0,\\
\frac{\partial(\rho \mathbf{u})}{\partial t} + \nabla\cdot(\rho
\mathbf{u} \otimes \mathbf{u})
+ \nabla p &= -\rho\nabla\Phi+\mu\triangle \mathbf{u},\\
\triangle \Phi &= 4\pi \rho,\label{NSP}
\end{split}
\end{equation}
where $t\geq 0$, $\mathbf{x}\in\mathbb{R}^3$, $\rho\geq 0$ is the
density, $\mathbf{u}\in\mathbb{R}^3$ the velocity, $p$ the pressure
of the gas, $\Phi$ the potential function of the self-gravitational
force, and $\mu>0$ the constant viscosity coefficient. We consider
polytropic gases and the equation of state is given by
\[
p=A\rho^\gamma
\]
where $A$ is an entropy constant and $\gamma>1$ is an adiabatic
exponent; in this case, the motion is called barotropic, which means
the pressure does not depend on the temperature or specific entropy.
 Values of $\gamma$ have their own physical
significance \cite{ch}; for example, $\gamma=\frac{5}{3}$ stands for
monatomic gas, $\frac{7}{5}$ for diatomic gas, $\gamma \searrow 1$
for heavier molecules. These $\gamma$'s also take important part in
the existence, uniqueness, and stability of stationary solutions,
for instance, see \cite{dlyy,j2,lin} for inviscid gaseous stars
modeled by the Euler-Poisson system.

For the spherically symmetric motion, i.e.
$\rho(t,\mathbf{x})=\rho(t,r)$ and
$\mathbf{u}(t,\mathbf{x})=u(t,r)\frac{\mathbf{x}}{r}$, where $u$ is
a scalar function and $r=|\mathbf{x}|$, (\ref{NSP}) can be written
as follows:
\begin{equation}
\begin{split}
\rho_t + \frac{1}{r^2}(r^2\rho u)_r&= 0,\\
\rho u_t + \rho u u_r + p_r +\frac{4\pi\rho}{r^2}\int_{0}^{r}\rho
s^2 ds &=\mu(u_{rr}+\frac{2u_r}{r}-\frac{2u}{r^2}).\label{nspE}
\end{split}
\end{equation}
Stationary solutions $\rho=\rho_0(r)$ and $u=0$, non-moving gaseous
spheres, satisfy the following:
\begin{equation}
(p_0)_r +\frac{4\pi\rho_0}{r^2}\int_{0}^{r}\rho_0 s^2 ds=0.\label{s}
\end{equation}
Note that the viscosity has nothing to do with these static
solutions themselves, namely, they are also stationary solutions of
the Euler-Poisson system \cite{j2}. The ordinary differential
equation (\ref{s}) can be transformed into the famous Lane-Emden
equation, and solutions of (\ref{s}) can be characterized according
to $\gamma$ in the following fashion: for given finite total mass,
if $\gamma\in (\frac{6}{5},2)$, then
 there exists at least one compactly supported stationary solution
 $\rho_0$. For $\gamma\in(\frac{4}{3},2)$,
 every stationary solution is compactly supported and unique.
If $\gamma=\frac{6}{5}$, there is a unique
 solution $\rho$ with infinite support,
 and it can be written explicitly in terms of
the Lane-Emden function.  On the other hand, if
$\gamma\in(1,\frac{6}{5})$, there are no stationary solutions with
finite total mass. For $\gamma\in(\frac{6}{5},2)$, letting $r=R$ be
the finite vacuum boundary of steady stars, it is well known
\cite{lin} that
\begin{equation}
\rho_0(r)\sim (R-r)^{\frac{1}{\gamma-1}}\;\text{ if }\;r\sim
R.\label{behavior}
\end{equation}
 We
observe that for given finite total mass, the maximum value of
density $\rho_0$ is inversely proportional to $R$ the radius  of
stars. For the details of stationary solutions such as the existence
and the asymptotic behavior, according to $\gamma$, we refer to
\cite{lin}.

Our main interest is the evolution of stars with finite radii as a
free boundary problem with the vacuum boundary to (\ref{nspE}). Let
$r\in [0,R(t)]$ and $t\in[0,T]$ for $T>0$. We look for $\rho(t,r)$
and $R(t)$ so that
\begin{equation}
 \rho>0 \text{ for } r<R(t)
\text{ and }\rho(t,R(t))=0;\label{b1}
\end{equation}
in particular, we would like to capture the behavior of interesting
stationary density profiles (\ref{behavior}) near the vacuum
boundary $r=R$. On the free boundary $r=R(t)$, we impose the
kinematic boundary condition
\begin{equation}
\frac{d}{dt}R(t)=u(t, R(t)),\label{b2}
\end{equation}
and the dynamic boundary condition
\begin{equation}
(\mu u_r-p)(t,R(t))=0.\label{b3}
\end{equation}
We remark that  $\rho(t,R(t))=0$ also serves the boundary condition
in our vacuum problem; in order to see this, we formally compute the
rate of density change in $t$ along the particle path $R(t)$:
\begin{equation*}
\begin{split}
\frac{d}{dt}\rho(&t, R(t))=\rho_t(t,R(t))+(\rho_ru)(t,R(t))=
-\rho(t, R(t))\frac{(r^2 u)_r}{r^2}(t, R(t))\\
&\Rightarrow \rho(t, R(t))=\rho(0,R(0))
\exp\{-\int_0^t(u_r+\frac{2u}{r})(\tau, R(\tau)) d\tau\}\\
&\text{or } \rho(t, R(t))R(t)^2=\rho(0, R(0))R_0^2\exp \{-\int_0^t
u_r(\tau, R(\tau)) d\tau\}.
\end{split}
\end{equation*}
Thus $\rho(t, R(t))=0$ for all time if it vanishes initially. Note
that the dynamic boundary condition with the vacuum boundary
condition leads to Neumann boundary condition
$$u_r(t,R(t))=0.$$

Due to lack of current mathematical tools to deal with the vacuum
boundary in Eulerian coordinates, it is desirable to introduce
Lagrangian formulation of (\ref{nspE}). We may assume that the total
mass of stars is $4\pi$, since the total mass is preserved in time:
\[
\begin{split}
\frac{d}{dt}\int_0^{R(t)}\rho(t,r)r^2 dr
&=\rho(t,R(t))R^2(t)\frac{d}{dt}R(t)+\int_0^{R(t)}\rho_t(t,r)r^2
dr\\&=\rho(t,R(t))R^2(t) u(t,R(t))- R^2(t)\rho(t,R(t))
u(t,R(t))\\&=0
\end{split}
\]
where we have used the boundary condition (\ref{b2}) and the
continuity equation at the second equality. Now we introduce a new
independent variable $x$:
\begin{equation}
x\equiv\int_0^r \rho s^2 ds,\label{x}
\end{equation}
which is a Lagrangian (mass) variable. The domain of $x$ is $[0,1]$,
since the total mass is assumed to be $4\pi$. Denote Lagrangian
derivatives by $D_t,D_x$. By change of variables, the Lagrangian
mass coordinate system $(t,x)$ and the Eulerian coordinates $(t,r)$
obey the following relation:
\begin{equation}
\begin{split}
D_t=\partial_t+ (D_t r)\partial_r, \;D_x=
(D_xr)\partial_r.\label{cv1}
\end{split}
\end{equation}
In Lagrangian coordinates, $r$ is not an independent variable, but a
function of $t,x$. To investigate the
 dynamics of $r$, first fix $x=x_0$. Then $r=r(t,x_0)\equiv r(t)$
is a particle path, and by taking $\frac{d}{dt}$ of (\ref{x}) and by
the continuity equation, we get
\[
0= \rho (t, r(t)) r^2(t) (D_t r)+\int_0^{r(t)}\partial_t\rho s^2 ds=
\rho (t, r(t)) r^2(t) (D_t r(t))-r^2(t)\rho(t, r(t)) u(t,r(t))
\]
and hence $D_tr=u$. Since $\partial_r x=\rho r^2$ from (\ref{x}), by
using the inverse function differentiation, we obtain
$D_xr=\frac{1}{\rho r^2}$. We therefore conclude the dynamics of $r$
in Lagrangian formulation as the following:
\begin{equation}
D_tr=u,\;D_xr=\frac{1}{\rho r^2}.\label{r}
\end{equation}
The second relation formally leads to
\begin{equation}
r=\{3\int_0^x\frac{1}{\rho}dy\}^{\frac{1}{3}}.\label{rr}
\end{equation}
We notice that the dynamics of $\rho, u$ completely determines $r$.
As in (\ref{cv1}), $\partial_t,\partial_r$ can be expressed in terms
of $D_t,D_x$:
\[
\partial_t=D_t-r^2\rho uD_x,\;\partial_r=\rho r^2 D_x.
\]
By using this change of variables with (\ref{r}), one can easily
check that (\ref{nspE}) can be written in Lagrangian coordinates
$(t,x)$ as follows: for $0\leq x\leq 1$,
\begin{equation}
\begin{split}
D_t\rho + \rho^2r^2 D_xu+\frac{2\rho u}{r}=& 0,\\
D_t u + r^2D_xp +\frac{4\pi x}{r^2}+\mu\frac{2u}{\rho r^2}
=&\mu D_x(\rho r^4D_xu),\label{nspL}\\
\text{or } D_t u + r^2D_xp +\frac{4\pi x}{r^2} =\mu r^2D_x(&\rho
r^2D_xu+\frac{2u}{r}),
\end{split}
\end{equation}
with the boundary conditions
\begin{equation}
u(t,0)=0,\; (\mu\rho r^2D_xu-p)(t,1)=0,\; \text{ and }
\rho(t,1)=0.\label{bcl}
\end{equation}
We point out that $D_t\equiv\partial_t+u\partial_r$ represents the
material derivative in the spherically symmetric motion. The vacuum
boundary condition also yields  $(\rho r^2 D_xu)(t, 1)=0$, but note
that we cannot reduce to $D_x u=0$ due to the vanishing property of
$\rho$. A free boundary $R(t)$ corresponds to a fixed boundary $x=1$
in Lagrangian formulation. It is easy to check that Jacobian of its
transformation is $\rho r^2$ and hence two formulations are
equivalent as long as $\rho>0$ and $r>0$. Note that $u$ needs not to
be zero along $x=1$.

Let $\rho_0=\rho_0(r)$ be the given stationary profile in Eulerian
coordinates as a solution of (\ref{s}) with the total mass $4\pi$.
The decay rate of $\rho_0$ towards $x=1$ in Lagrangian coordinates
is asymptotically given as the following:
\begin{equation}
\rho_0\sim (1-x)^{\frac{1}{\gamma}}\;\text{ if }\; x\sim
1\label{behaviorl}
\end{equation}
from (\ref{x}) and (\ref{behavior}).

Variations of this model have been considered in the literature:
Okada and Makino \cite{om, mom} studied global weak solution,
uniqueness, and stability to the Navier-Stokes equations for gas
surrounding a solid ball (a hard core) without self-gravitation in
Lagrangian coordinates as a free boundary problem when the density
distribution contacts with the vacuum at a finite radius. More
recently, a similar study has been done by Ducomet and Zlotnik
\cite{dz} for viscous compressible flow driven by gravitation and an
outer pressure, proving interesting stabilization results. However,
all their results are restricted to cutoff domains excluding some
neighborhood of the origin, since their analysis strongly depends on
the one-dimensional structure of symmetric flows. On the other hand,
when the initial density is away from the vacuum for smooth initial
data or discontinuous data, one-dimensional or multidimensional
problem, a lot of progress has been made on the compressible
Navier-Stokes equations. However, as far as the physical vacuum is
concerned, very few rigorous results are available for compressible
flows. For one-dimensional viscous gas flow, there has been some
important progress; in particular, the vacuum interface behavior as
well as the regularity to one-dimensional Navier-Stokes free
boundary problems were investigated by Luo, Xin, and Yang
\cite{lxy}. It is important and interesting to understand the
dynamics of the Navier-Stokes-Poisson system in three space
dimension as a vacuum free boundary problem displaying the feature
of physical vacuum boundary behavior (\ref{behavior}) or
(\ref{behaviorl}).

This article concerns the local well-posedness of a free boundary
problem to the Navier-Stokes-Poisson system (\ref{nspE}) with
(\ref{b1}), (\ref{b2}), and (\ref{b3}), or (\ref{nspL}) with
(\ref{bcl}), embracing the physical vacuum boundary behavior
(\ref{behavior}) or (\ref{behaviorl}) of stationary density
profiles, especially when $\frac{6}{5}<\gamma<2$.  To establish
local in time strong solutions including the physical vacuum
boundary, we utilize both Eulerian and Lagrangian formulations. \\

Before we state the main results, we first define the suitable
energy space in which regular solutions reside.  Let the initial
density profile be given by $\rho_{in}$ satisfying the following
conditions:
\begin{equation}
\begin{split}
\text{(i)}&\;\; \rho_{in}> 0 \text{ for }r<R\;\;\text{and}\;\;\rho_{in}(R)=0, \\
\text{(ii)}&\;\int_0^R \rho_{in} s^2 ds = 1 \;(\text{total mass}=4\pi).\\
\end{split}
\end{equation}
Consider $0<r_0<r_1<r_2<R$ such that
\begin{equation}
\begin{split}
0<2d<r_0,\;0<3d<r_2-r_1,\;\frac{1}{r_0-d}\leq 1, \label{condr}
\end{split}
\end{equation}
for small fixed constant $d$. Now let $x_i$ be the initial position
in Lagrangian coordinates corresponding to $r_i$ where $i=0,1,2$:
\begin{equation}
x_i=\int_0^{r_i}\rho_{in} s^2 ds.
\end{equation}
Then by the positivity of $\rho_{in}$, we get $0<x_0<x_1<x_2<1$.
Denoting the particle path emanating from $r_i$ by $r_i(t)$,
$r_i(t)$ characterizes $x_i$:
\begin{equation}
\frac{d}{dt}\int_0^{r_i(t)}\rho(s,t)s^2 ds =0, \text{ i.e.
}\int_0^{r_i(t)}\rho(s,t)s^2 ds =x_i.
\end{equation}
This relation is the conservation of mass and can be readily
verified by using the continuity equation and
$$\frac{d}{dt}r(t)=u(r(t),t).$$
Assume $|u(r,t)|\leq K$ for all $0\leq r\leq R(t)$ and $0\leq t\leq
T$, where $T$ is sufficiently small. In particular, $d$ in
(\ref{condr}) will be chosen so that $KT\leq d$.  Note that the
smallness assumption on the time interval $T$ prevents a dramatic
change of $r$ in time.

We note that the initial data $\rho_{in}$ and $u_{in}$ given in
Eulerian coordinates can be also regarded as functions of $x$ in
Lagrangian coordinates, since (\ref{x}) is valid when $t=0$:
$x=\int_0^r\rho_{in}s^2ds$. The initial value of $r$ in Lagrangian
coordinates is given or defined as the following:
$r_{in}(x)=(3\int_0^x\frac{1}{\rho_{in}}dy)^{\frac{1}{3}}$ from
(\ref{rr}). As the first preparation of the rigorous analysis, we
introduce the following cutoff functions $\chi$ and $\zeta$. Let
$\chi\in \mathcal{C}^{\infty}[0,1]$ be a smooth function of $x$ such
that
\[
\begin{split}
\text{(i)}&\;\;0\leq \chi \leq 1\;\text{ and}
\;\;\text{supp}(\chi)\subset[x_0,1],\\
\text{(ii)}&\;\;\chi(x)=1 \;\text{ if }x_1\leq x\leq 1,\\
\text{(iii)}&\;\;|\chi'|\leq\frac{C}{x_1-x_0}\;\text{ and }\;
|\chi''|<\infty.
\end{split}
\]
$\chi$ is a function of both $r$ and $t$ in Eulerian coordinates.
Note that $$|r_i(t)-r_i|\leq d\text{ for }0\leq t\leq T,$$ since
$$r_i(t)=r_i+\int_0^t u(r(\tau),\tau)d\tau \text{ by }(\ref{r}).$$
Hence we deduce that for $0\leq t\leq T$,
\[
\chi(r,t)=0\;\;\text{if }r\leq r_0-d\;\text{ and }\;
\chi(r,t)=1\;\;\text{if }r\geq r_1+d.
\]
Similarly, construct a smooth function $\zeta$ of $r$ such that
\[
\begin{split}
\text{(i)}&\;\;0\leq \zeta \leq 1\;\text{ and}\;\;\text{supp}(\zeta)
\subset[0,r_2-d],\\
\text{(ii)}&\;\;\zeta(r)=1 \;\text{ if }0\leq r\leq r_1+d,\\
\text{(iii)}&\;\;|\zeta'|\leq\frac{C}{r_2-r_1-2d} \;\text{ and }\;
|\zeta''|<\infty.
\end{split}
\]
Then as a function of $x$ and $t$, $\zeta$ satisfies the following:
for $0\leq t\leq T$,
\[
\zeta(x,t)=1\;\;\text{if }x\leq x_1\;\text{ and }\;
\zeta(x,t)=0\;\;\text{if }x\geq x_2.
\]
We will freely view $\chi$ and $\zeta$ as functions of  $x$, $r$ and
$t$ without confusion.

Define the energy functional $\mathcal{E}(t)$ by the sum of the
Eulerian energy $\mathcal{E}_E(t)$ and the Lagrangian energy
$\mathcal{E}_L(t)$: $$\mathcal{E}(t)
=\mathcal{E}_E(t)+\mathcal{E}_L(t)$$ where
\begin{equation}
\mathcal{E}_E(t)\equiv \frac{1}{2}\sum_{|\alpha|\leq
3}\{A\gamma\int\zeta \rho^{\gamma-2}
|\partial^{\alpha}\rho|^2d\mathbf{x}+
\int\zeta\rho|\partial^{\alpha}
\mathbf{u}|^2d\mathbf{x}\}.\label{ee}
\end{equation}
Here $\partial^{\alpha}$ represents all the Eulerian mixed
derivatives and we use $\mathbf{u},\mathbf{x}$ in order to emphasize
they are vector quantities;
\begin{equation}
\begin{split}
\mathcal{E}_L(t)\equiv&\;\frac{1}{2}\int_{x_0}^1\chi
u^2dx+\frac{A}{\gamma-1}\int_{x_0}^1\chi\rho^{\gamma-1}dx+
\frac{1}{2}\sum_{j=1}^3\{\int_{x_0}^1\chi |D_t^ju|^2dx\}\\
&+\sum_{j=0}^2\{\frac{\mu}{2}\int_{x_0}^1\chi\{\rho r^4|D_t^jD_xu|^2
+\frac{2|D_t^ju|^2}{\rho r^2}\}dx+\li\chi\rho^{2\gamma-2}r^4|D_t^jD_x\rho|^2dx\}\\
&+\sum_{j=0}^1
\frac{1}{2}\li\chi\rho^{4\gamma-2}r^8|D_t^jD_{x}^2\rho|^2dx
+\frac{1}{2}\li\chi\rho^{8\gamma-2}r^{12}|D_{x}^3\rho|^2dx.\label{le}
\end{split}
\end{equation}
We first observe that some derivative terms are missing in
$\mathcal{E}_L(t)$ and will show that missing terms are controlled
by $\mathcal{E}_L(t)$ via the equations. More specifically, for the
$u$-part, the terms involving more than one spatial derivative can
be estimated by the momentum equation since they are closely related
to the viscosity term, which is represented by the sum of lower
derivative terms. For the $\rho$-part, pure time derivative terms
can be directly estimated by the continuity equation. We also
observe that the energy defined in the above involves the pressure
($p=A\rho^{\gamma}$) rather than $\rho$ itself; for instance, the
first derivative energy of the $\rho$-part
$\li\chi\rho^{2\gamma-2}r^4|D_x\rho|^2dx $ is equal to
$\frac{1}{(A\gamma)^2}\li\chi r^4|D_xp|^2dx $. Indeed, the pressure
turns out to be the right quantity to look at, and moreover, because
the behavior of $p$ $(\sim 1-x)$ around the vacuum boundary in
Lagrangian coordinates does not depend on $\gamma$, our analysis can
embrace all the physically interesting $\gamma$'s. In particular,
the different weights in front of different spatial derivatives of
$\rho$ have been carefully chosen so that not only the energy space
can capture the behavior of stationary profiles, but also the energy
estimates can be closed.

In turn we define the dissipation
$\mathcal{D}(t)\equiv\mathcal{D}_E(t)+\mathcal{D}_L(t)$ by
\begin{equation}
\begin{split}
{\mu}\sum_{|\alpha|\leq
3}\int\zeta|\nabla\partial^{\alpha}\mathbf{u}|^2
d\mathbf{x}+{\mu}\sum_{j=0}^3\li\chi\{\rho r^4
|D_t^jD_xu|^2+\frac{2|D_t^ju|^2}{\rho
r^2}\}dx+\sum_{j=1}^3\li\chi|D_t^iu|^2dx.\label{diss}
\end{split}
\end{equation}
Next, we introduce the following assumption ($K$):
\begin{equation}
\begin{split}
\sup_{x_0\leq x\leq 1}\{\rho, |\frac{u}{r}|, |\rho
r^2D_xu+\frac{2u}{r}|=|\frac{D_t\rho}{\rho}|, |\rho r^2D_xu|,|\rho
r^2 D_tD_{x}u|,|\frac{D_tu}{r}|,|\rho^{2\gamma-1}r^2D_x\rho| \}\leq K\\
\sup_{0\leq r\leq r_2-d}\{\rho, |u|, |\partial_ru|,
|\frac{\partial_t\rho}{\rho}|, |\frac{\partial_r\rho}{\rho}|,
|\partial_t u| \}\leq K\label{K}
\end{split}
\end{equation}
which indicates what regularity strong solutions should enjoy. It is
shown in Lemma \ref{weaving} that $K$ is bounded by the energy
$\mathcal{E}$. We
are now ready to state the main a priori estimates.\\

\begin{theorem}
Suppose $\rho, u$ are smooth solutions to the free boundary problem
of the Navier-Stokes-Poisson system (\ref{nspE}) with (\ref{b1}),
(\ref{b2}), and (\ref{b3}), or (\ref{nspL}) with (\ref{bcl}) for
given initial data such that $\mathcal{E}(0)\equiv
\mathcal{E}(\rho_{in},u_{in})$ is bounded. Then there exist
$\mathcal{C}_1=C_1(K),\;\mathcal{C}_2>0$ such that the following
energy inequality holds for $0\leq t<\frac{d}{K}$,
\begin{equation}
\frac{d}{dt}\mathcal{E}(t)+\frac{1}{2}\mathcal{D}(t)\leq
\mathcal{C}_1\mathcal{E}(t)+\mathcal{C}_2(\mathcal{E}(t))^2\label{ei}
\end{equation}
Moreover, there exist $T>0$ and $A=A(T,C_1,C_2,\mathcal{E}(0))>0$
such that  $$\sup_{0\leq t\leq T}\mathcal{E}(t)\leq A.
$$\label{prop}
\end{theorem}


In the next theorem, we establish the local in time well-posedness
of strong solutions
to the Navier-Stokes-Poisson system. \\

\begin{theorem} Let the initial data $\rho_{in},u_{in}$ be given
such that $\mathcal{E}(0)\equiv
\mathcal{E}(\rho_{in},u_{in})<\infty$.
 There exists $T^{\ast}>0$ such that there exists a unique solution  $R(t)$,
$\rho(t,r)$, $u(t,r)$ to the Navier-Stokes-Poisson system
(\ref{nspE}) with (\ref{b1}), (\ref{b2}), (\ref{b3}) in
$[0,T^{\ast})\times[0,R(t)]$  such that $$\sup_{0\leq t\leq
T^{\ast}}\mathcal{E}(t)\leq 2\mathcal{E}(0).$$ Moreover,
$\rho(t,x)$, $u(t,x)$, $r(t,x)$ where $x$ is a Lagrangian variable
defined in (\ref{x}), serve a unique solution to (\ref{nspL}) with
(\ref{bcl}) in $[0,T^{\ast})\times[0,1]$.\label{thm}
\end{theorem}\

We remark that the energy $\mathcal{E}$ of the stationary solutions
$\rho_0$ is bounded for all $\gamma$, and therefore strong solutions
constructed in Theorem \ref{thm} can capture the physical vacuum
boundary behavior (\ref{behavior}) or (\ref{behaviorl}) locally in
time. We believe that this local well-posedness result provides the
foundation towards further interesting study such as global
well-posedness under the same energy space and nonlinear stability
questions of Lane-Emden steady stars.

The central difficulty in this article is to deal with the vacuum
free boundary where the density vanishes at certain rate, which
makes the system degenerate along the boundary. Since the free
boundary gets fixed in Lagrangian coordinates, it is desirable to
work in Lagrangian framework and the degeneracy from the vacuum is
overcome with the density-weighted energy estimates. We note that
these energy estimates can be closed due to the presence of the
viscosity; the smoothing effect of the velocity is  transferred to
the density. We should also point out that the vacuum boundary
condition $\rho(t,R(t))=0$ solely cannot define the free boundary
properly but rather the dynamic boundary condition (\ref{b3})
characterizes the free boundary problem: indeed, it enables to
integrate by parts up to the boundary. On the other hand, the
three-dimensional structure of symmetric flows is still prevalent
around the origin, the coordinate singularity from symmetry is
worrisome, and thus the cooperation of Eulerian formulation seems
necessary. We perform Lagrangian and Eulerian energy estimates
separately by using the cutoff functions. Overlapping terms
involving $\chi'$ and $\zeta'(\nabla\zeta)$ do not cause any other
essential difficulty and are controlled by the opposite energies.

Another key idea is to implement an appropriate iteration scheme
whose approximate solutions converge to the desired strong solution
to the system under the energy norm $\mathcal{E}$. This can be done
by the separation of the density and the velocity from the system in
Lagrangian coordinates. With given fixed $\rho^n$, $r^n$ related
coefficients, we obtain $u^{n+1}$ by solving linear parabolic PDEs.
$\rho^{n+1}$ is defined from the continuity equation along the
particle path by using $\rho^n$, $r^n$, and $u^{n+1}$. In turn,
$r^{n+1}$ is
 determined by $\rho^{n+1}$ from the dynamics of $r$.
 The $(n+1)$-th total energy $\mathcal{E}^{n+1}$ is
accordingly further separated into $\rho$-part and $u$-part and we
prove that they are uniformly bounded. We remark that this whole
energy separation works due to the viscosity term, in order words,
the viscosity dominates the given vacuum problem.

The method developed in this paper is lucid and concrete, and
moreover, it provides the critical quantities, such as $\rho r^2
D_xu+\frac{2u}{r}$, which govern the whole dynamics. We believe that
this method itself will have rich applications to other interesting
problems. We note that our results strongly depend on the given
initial data in that the initial density is explicitly embedded in
the construction of the cutoff functions, which are important
components of the energy $\mathcal{E}$.

It would be very interesting, challenging both physically and
mathematically to study the full system without the symmetry
assumption as a free boundary problem; in the general case, no
result is known for compressible gas flows with the free boundary.
We note that the existence problem for the pure compressible
Navier-Stokes system in three space dimension is still open.
Currently, the above argument does not seem to extend directly to
the general case. However, we believe that the methodology can make
a contribution to the study of the full system with the vacuum
boundary, along with the recent progress on the free boundary
problems from other contexts. We will leave them for future study.

The article proceeds as follows. In the first half of the paper, we
establish the a priori estimates, which should shed some light on
the construction of strong solutions. In Section \ref{section3},
boundary estimates are performed in Lagrangian coordinates and in
Section \ref{section4}, interior estimates are presented in Eulerian
coordinates. By weaving those two estimates and verifying that the
smoothness assumption can be closed under the energy $\mathcal{E}$,
we conclude the a priori estimates in Section \ref{section5}. In the
rest of the article, strong solutions are constructed by
implementing an iteration scheme, based on the separation of the
density and the velocity from the system. In Section \ref{section6},
the approximate scheme is displayed and approximate solutions at
each step are shown to have the same regularity of the previous
approximations. Section \ref{section8} is devoted to
obtaining uniform energy estimates and Theorem \ref{thm} is proven. \\

\section{Boundary Estimates in Lagrangian Coordinates}
\label{section3}

First, we observe  that the behavior of the density $\rho(t,x)$ in
Lagrangian coordinates should be inherited from the initial profile
$\rho_{in}$. This can be readily seen from the continuity equation
in (\ref{nspL}):
\begin{equation}
\rho(t,x)=\rho_{in}(x)\exp\{-\int_0^t(\rho r^2
D_xu+\frac{2u}{r})(\tau,x) d\tau\}\label{rho1}
\end{equation}
As long as $|\rho r^2 D_xu+\frac{2u}{r}|$ is bounded, so is $\rho$;
in Section \ref{section5}, we will show that $$\sup_{0< x< 1}|\rho
r^2 D_xu+\frac{2u}{r}|\leq
C_{in}(\mathcal{E}(t)^{\frac{1}{2}}+\mathcal{E}(t))$$ where $C_{in}$
depends only on the initial density $\rho_{in}$. Indeed, the
quantity $|\rho r^2 D_xu+\frac{2u}{r}|$ as well as the formula
$\rho$ (\ref{rho1}) are essential to establish the local
well-posedness. Let
\[
M\equiv \sup_{0< x< 1}|\rho r^2 D_xu+\frac{2u}{r}|.
\]
Note that $M\leq K$, where $K$ appears in (\ref{K}). We also observe
that $\rho(t,x)$ can be controlled by $\rho_{in}$ and $M$ from
(\ref{rho1}): for $0\leq t\leq T$
\begin{equation}
\rho_{in}(x)e^{-MT}\leq\rho(t,x)\leq \rho_{in}(x)e^{MT}\label{rhob}
\end{equation}
Thus for sufficiently small $T$, $\rho$ stays close to initial
density $\rho_{in}$. Integrating the momentum equation in $x$, one
might hope to get $ \mu(\rho r^2
D_xu+\frac{2u}{r})=p+\int_1^x(\frac{D_tu}{r^2}+\frac{4\pi
y}{r^4})dy. $ However, we remark that this computation is absurd
because $u$ does not have to be zero on $x=1$, while $(\mu\rho
r^2D_xu-p)(t,1)=0.$

We establish the estimates of $\mathcal{E}_L(t)$ by a series of
lemmas. The following lemma treats $t$-derivatives of
$u$.\\

\begin{lemma}
There exists $C_K>0$ such that
\begin{equation}
\frac{1}{2}\frac{d}{dt}\sum_{i=0}^3\li\chi|D_t^iu|^2dx+\frac{3\mu}{4}\sum_{i=0}
^3\li\chi \{\rho r^4|D_t^iD_xu|^2+\frac{2|D_t^iu|^2}{\rho r^2}\}dx
\leq C_K \mathcal{E}_L+ \mathcal{OL}_1,\label{lemma31}
\end{equation}
where $\mathcal{OL}_1\leq \widetilde{C}_K\mathcal{E}_E$ for some
constant $\widetilde{C}_K$ \label{tl}.
\end{lemma}

\begin{proof} It is instructive to see how the estimates go in detail at the
energy level . Let $i=0$. Multiply the momentum equation in
(\ref{nspL}) by $\chi u$ and integrate to get
\begin{equation}
\frac{1}{2}\frac{d}{dt}\int_{x_0}^1 \chi u^2 dx + \int_{x_0}^1\chi
ur^2 D_xp dx-\mu \int_{x_0}^1 \chi u D_x(\rho r^4 D_xu)
dx+\int_{x_0}^1 \chi u \frac{4\pi
x}{r^2}dx+\mu\li\chi\frac{2u^2}{\rho r^2}dx =0.\label{en0}
\end{equation}
The second and third terms in (\ref{en0}): by integrating by parts
and using the boundary condition (\ref{bcl})
\[
\begin{split}
\int_{x_0}^1\chi ur^2 D_xp dx-\mu \int_{x_0}^1 \chi u D_x(\rho r^4
D_xu) dx&=-\int_{x_0}^1D_x(\chi u r^2)p dx+
\mu\li D_x(\chi u)\rho r^4 D_xu dx\\
=-\int_{x_0}^1\chi' u r^2 p dx+\mu\li\chi'u\rho r^4 D_xu dx&
\underline{-\int_{x_0}^1\chi D_x(u r^2)p dx}_{\star}+\mu\li\chi\rho
r^4|D_xu|^2dx
\end{split}
\]
Use the continuity equation to reduce $\star$ to
\[
\star=\frac{A}{\gamma-1}\frac{d}{dt}
\int_{x_0}^1\chi\rho^{\gamma-1}dx.
\]
Overlapping terms can be controlled as follows:
\[
\begin{split}
&-\int_{x_0}^{x_1}\chi' u r^2 p dx+\mu\int_{x_0}^{x_1}\chi'u\rho r^4
D_xu dx\leq\frac{C}{x_1-x_0}\{\int_{r_0(t)}^{r_1(t)}
\rho^{\gamma+1}|u| r^4dr+\mu\int_{r_0(t)}^{r_1(t)}\rho|u\partial_ru|
r^4dr\}\\
&\leq \frac{C(r_1+d)^2}{x_1-x_0}\{\sup_{r\leq
r_1+d}\rho^{\gamma+1}(\int_{r_0(t)}^{r_1(t)}\rho^{\gamma}
r^2dr+\int_{r_0(t)}^{r_1(t)}\rho
u^2r^2dr)+\int_{r_0(t)}^{r_1(t)}\rho
u^2r^2dr+\int_{r_0(t)}^{r_1(t)} \rho |\partial_ru|^2r^2dr\}\\
&\leq \widetilde{C}_K\mathcal{E}_E
\end{split}
\]
where we have used the relations $dx=\rho r^2 dr;\;
D_x=\frac{1}{\rho r^2}\partial_r$ as well as the Cauchy-Schwarz
inequality. As for the fourth term in (\ref{en0}), we apply the
Cauchy-Schwarz inequality:
\[
\begin{split}
\int_{x_0}^1 \chi u \frac{4\pi x}{r^2}dx&\leq
\frac{\mu}{2}\li\chi\frac{u^2}{\rho
r^2}dx+\frac{8\pi^2}{\mu}\li\chi\rho\frac{x^2}{r^2}dx\\
&\leq \frac{\mu}{2}\li\chi\frac{u^2}{\rho
r^2}dx+\frac{C}{(r_0-d)^2}\sup_{x_0\leq x\leq
1}|\rho^{2-\gamma}|\li\chi\rho^{\gamma-1}dx
\end{split}
\]
where we have used $\rho=\rho^{2-\gamma}\rho^{\gamma-1}$ at the last
inequality. After absorbing the viscosity term into the LHS and
using the assumption (\ref{K}), we get the following:
\[
\begin{split}
\frac{d}{dt}&\int_{x_0}^1\chi\{\frac{1}{2} u^2+\frac{A}{\gamma-1}
\rho^{\gamma-1}\}dx+ \frac{3\mu}{4} \int_{x_0}^1 \chi\{\rho
r^4|D_xu|^2+
\frac{2u^2}{\rho r^2}\}dx\\
\leq&\; C_K\li\chi\rho^{\gamma-1}dx +\widetilde{C}_K\mathcal{E}_E
\end{split}
\]
Note that from the continuity equation, we also get
\begin{equation}
\li\chi \rho^{-3}|D_t\rho|^2dx\leq 3\int_{x_0}^1 \chi \{\rho
r^4|D_xu|^2+ \frac{2u^2}{\rho r^2}\}dx.\label{rt}
\end{equation}
Now let $i=1$. Take $D_t$ of the momentum equation to get
\[
\begin{split}
D_t^{2}u+r^2D_xD_tp&-\mu D_x(\rho r^4D_x D_tu)+\mu\frac{2D_tu}{\rho
r^2}\\&=\mu D_x(D_t(\rho r^4)D_xu)-D_tr^2 D_xp-D_t(\frac{2\mu}{\rho
r^2})u-D_t(\frac{4\pi x}{r^2})
\end{split}
\]
Multiply by $\chi D_tu$ and integrate in $x$ to get
\[
\begin{split}
&\frac{1}{2}\frac{d}{dt}\li\chi|D_tu|^2dx+\underline{\int_{x_0}^1\chi
D_tur^2 D_xD_tp dx-\mu \int_{x_0}^1 \chi D_tu D_x(\rho r^4
D_tD_xu) dx}_{(i)}+\mu\li\chi\frac{2|D_tu|^2}{\rho r^2}dx\\
&=\underline{\li\chi\{\mu D_x(D_t(\rho r^4)D_xu)-D_tr^2
D_xp-D_t(\frac{2\mu}{\rho r^2})u-D_t(\frac{4\pi
x}{r^2})\}D_tudx}_{(ii)}
\end{split}
\]
The LHS can be estimated in the same way as in the zeroth order
case. The RHS has new terms but it consists of lower order
derivative terms. We will provide the detailed computation in the
below. First, integrate by parts by using the boundary condition
(\ref{bcl}) in $(i)$ to have
\[
\begin{split}
(i)=&-A\gamma\li\chi'D_tur^2 \rho^{\gamma-1}D_t\rho
dx+\mu\underline{\li\chi' D_tu\rho r^4 D_xD_tu
dx}_{(a)}\\&-\underline{A\gamma\li\chi (r^2D_tD_xu
+\frac{2D_tu}{\rho r})\rho^{\gamma-1}D_t\rho dx}_{(b)}+\mu\li\chi
\rho r^4|D_tD_xu|^2dx
\end{split}
\]
Note that $D_xD_tu$ in $(a)$ is not worrisome:  by another
integration by parts, $(a)$ becomes
\[
(a)=-\int_{x_0}^{x_1}\chi''\rho r^4|D_tu|^2dx-\int_{x_0}^{x_1}\chi'
D_x(\rho r^4)|D_tu|^2dx,
\]
which is bounded by the Eulerian energy. For $(b)$, we apply the
Cauchy-Schwarz inequality:
\begin{equation}
\begin{split}
(b)&\leq \frac{\mu}{8}\li\chi\{\rho
r^4|D_tD_xu|^2+\frac{|D_tu|^2}{\rho
r^2}\}dx+\frac{10A^2\gamma^2}{\mu}\li\chi
\rho^{2\gamma-3}|D_t\rho|^2dx\\
&\leq\frac{\mu}{4}\li\chi\{\rho r^4|D_tD_xu|^2+\frac{|D_tu|^2}{\rho
r^2}\}dx +\frac{30A^2\gamma^2}{\mu}\sup_{x_0\leq x\leq
1}|\rho^{2\gamma}|\int_{x_0}^1 \chi \{\rho r^4|D_xu|^2+
\frac{2u^2}{\rho r^2}\}dx\label{(b)}
\end{split}
\end{equation}
where we have used (\ref{rt}). Each term in the LHS $(ii)$ is
treated as follows: For the first term we integrate by parts using
the boundary condition.
\[
\mu\li\chi D_x(D_t(\rho r^4)D_xu)D_tudx=-\mu\li\chi'D_t(\rho
r^4)D_xuD_tudx\underline{-\mu\li\chi D_t(\rho r^4)D_xu
D_xD_tudx}_{(c)}
\]
The first term in the RHS is an overlapping term, bounded by
$\widetilde{C}_K\mathcal{E}_E$, and we apply the Cauchy-Schwarz
inequality for the second term.
\[
(c)\leq\frac{\mu}{8}\li\chi\rho r^4|D_tD_xu|^2dx+ 2\mu\sup_{x_0\leq
x\leq 1}|\frac{D_t\rho}{\rho}+\frac{4u}{r}|^2\li\chi\rho
r^4|D_xu|^2dx
\]
By recalling $D_tr=u$ and applying the Cauchy-Schwarz inequality,
the second term in $(ii)$ is bounded by $CK\mathcal{E}_L$ as
follows:
\[
\begin{split}
-\li\chi D_tr^2D_xpD_tudx&=-2A\li\chi\frac{u}{r}
r^2\rho^{\gamma-1}D_x\rho D_tudx\\
&\leq A\sup_{x_0\leq x\leq
1}|\frac{u}{r}|\{\li\chi\rho^{2\gamma-2}r^4|D_x\rho|^2dx
+\li\chi|D_tu|^2 dx\}
\end{split}
\]
In the same way, one can check that the third term in $(ii)$ is
bounded by
\[
-2\mu\li\chi D_t(\frac{1}{\rho r^2})u D_tudx\leq \frac{\mu}{8}\li
\chi \frac{|D_tu|^2}{\rho r^2}dx+8\mu \sup_{x_0\leq x\leq
1}|\frac{D_t\rho}{\rho}+\frac{2u}{r}|^2\li\chi\frac{u^2}{\rho r^2}
dx
\]
The last term in $(ii)$ is estimated in the following
\[
-4\pi\li\chi D_t(\frac{x}{r^2})D_tudx=8\pi\li\chi \frac{u
x}{r^3}D_tudx\leq\frac{\mu}{8}\li \chi \frac{|D_tu|^2}{\rho r^2}dx+
\frac{32\pi^2}{\mu}\sup_{x_0\leq x\leq 1}|\frac{\rho}{r^4}|\li\chi
u^2dx
\]
Thus, after absorbing the viscosity term into the LHS and using the
assumption (\ref{K}), we get (\ref{lemma31}) for $i=1$. Now let
$i=2$ or 3. The spirit is same as in the case $i=1$. Take $D_t^i$ of
the momentum equation to get
\begin{equation}
\begin{split}
D_t^{i+1}u+r^2D_xD_t^ip-\mu D_x(\rho r^4D_x
D_t^iu)&+\mu\frac{2D_t^iu}{\rho r^2}\\=\sum_{j=0}^{i-1}\{\mu
D_x(D_t^{i-j}(\rho r^4)D_t^jD_xu)&-D_t^{i-j}r^2
D_xD_t^jp-D_t^{i-j}(\frac{2\mu}{\rho r^2})D_t^ju\}-D_t^i(\frac{4\pi
x}{r^2})\label{mom}
\end{split}
\end{equation}
Multiply by $\chi D_t^iu$ and integrate in $x$ to get
\[
\begin{split}
&\frac{1}{2}\frac{d}{dt}\li\chi|D_t^iu|^2dx+\int_{x_0}^1\chi
D_t^iur^2 D_xD_t^ip dx-\mu \int_{x_0}^1 \chi D_t^iu D_x(\rho r^4
D_xD_t^iu) dx+\mu\li\chi\frac{2|D_t^iu|^2}{\rho r^2}dx\\
&=\sum_{j=0}^{i-1}\li\chi\{\mu D_x(D_t^{i-j}(\rho
r^4)D_t^jD_xu)-D_t^{i-j}r^2 D_xD_t^jp-D_t^{i-j}(\frac{2\mu}{\rho
r^2})D_t^ju\}D_t^iudx\\
&\;\;\;-\int_{x_0}^1 \chi D_t^iu D_t^i(\frac{4\pi x}{r^2})dx
\end{split}
\]
Note each term in the RHS involves only lower order derivatives.
The second and third terms in the LHS can be written as the
following: by integrating by parts and using the boundary condition
(\ref{bcl})
\[
\begin{split}
-\int_{x_0}^1\chi' D_t^iu r^2 D_t^ip dx+\mu\li\chi'D_t^iu\rho r^4
D_xD_t^iu dx\underline{-\int_{x_0}^1\chi D_x(r^2D_t^iu)D_t^ip
dx}_{\star}+\mu\li\chi\rho r^4|D_xD_t^iu|^2dx
\end{split}
\]
The first two terms are overlapping terms and they can be bounded by
$\widetilde{C}_K\mathcal{E}_E$ by using the change of variable
$D_t=\partial_t+u\partial_r$. One might try to use the continuity
equation to reduce $\star$ to a $t$-derivative as in the zeroth
order case, but we point out that it is rather complicated. Instead,
we use the Cauchy-Schwarz inequality and take advantage of the
viscosity term to control it as we did in (\ref{(b)}). In order to
do so, it is enough to derive the following: for $0\leq i\leq 2$,
\[
D_t^{i+1}\rho=-\rho\{ \rho
r^2D_xD_t^{i}u+\frac{2D_t^{i}u}{r}\}-\sum_{0\leq j<i, 0\leq k\leq
i}D_t^{i-j-k}\rho^2D_x(D_t^kr^2D_t^{j}u)
\]
\begin{equation}
\Rightarrow \li\chi \rho^{-3}|D_t^{j}\rho|^2dx\leq
C_K\mathcal{E}_L\text{ for } 1\leq j\leq 3.
\end{equation}
Hence, by summing over $i$, we get the following:
\[
\frac{1}{2}\frac{d}{dt}\sum_{i=1}^3\li\chi|D_t^iu|^2dx+\frac{3\mu}{4}
\sum_{i=1}^3\li\chi\{\rho r^4|D_xD_t^iu|^2+\frac{2|D_t^iu|^2}{\rho
r^2}\}dx\leq  C_K \mathcal{E}_L+\widetilde{C}_K\mathcal{E}_E.
\]
This finishes the proof of Lemma \ref{tl}.
\end{proof}\

Next, we estimate mixed derivatives with only one spatial
derivative.\\

\begin{lemma}
There exists $C_K>0$ such that
\begin{equation}
\begin{split}
\frac{\mu}{2}\frac{d}{dt}\sum_{i=0}^2\li\chi\{\rho&
r^4|D_xD_t^iu|^2+\frac{2|D_t^iu|^2}{\rho
r^2}\}dx+\frac{1}{2}\sum_{i=0}^2\li\chi |D_t^{i+1}u|^2dx\\
&\leq \frac{\mu}{4}\sum_{i=1}^3\li\chi\{\rho
r^4|D_xD_t^iu|^2+\frac{2|D_t^iu|^2}{\rho
r^2}\}dx+C_K\mathcal{E}_L+\mathcal{OL}_2\label{onex}
\end{split}
\end{equation}
where $\mathcal{OL}_2\leq \widetilde{C}_K\mathcal{E}_E$.
\end{lemma}\

Note that part of the dissipation energy, which involves one more
derivative terms than the total energy, appears in the RHS of the
inequality (\ref{onex}). This is not a problem because it will
be absorbed into the dissipation obtained in Lemma \ref{tl}.\\

\begin{proof} We provide the estimate for $i=0$ in detail. Higher order
terms can be estimated in the same way, since taking $t$ derivative
does not destroy the structure of the equations. Multiply the
momentum equation by $D_t u$ and integrate to get:
\begin{equation*}
\begin{split}
\li\chi |D_t u|^2dx+\li\chi(r^2D_xp)D_tudx-\mu\li \chi D_x(\rho
r^4D_xu)D_tu dx +\mu\li\chi\frac{2uD_t u}{\rho r^2}dx\\
+\li\chi\frac{4\pi x}{r^2}D_tudx =0\label{ux}
\end{split}
\end{equation*}
The second term: by the Cauchy-Schwarz inequality
\[
\begin{split}
|\li\chi(r^2D_xp)D_tudx|&\leq \li\chi
r^4|D_xp|^2dx+\frac{1}{4}\li\chi|D_tu|^2 dx\\&=(A\gamma)^2\li\chi
\rho^{2\gamma-2}|D_x\rho|^2dx+\frac{1}{4}\li\chi|D_tu|^2 dx
\end{split}
\]
The third term: by integrating by parts and using the boundary
condition (\ref{bcl})
\[
\begin{split}
-&\mu\li\chi D_x(\rho r^4D_xu)D_tu dx=\mu\int_{x_0}^{x_1} \chi'\rho
r^4D_xu D_tudx+\underline{\mu\li\chi\rho r^4D_xuD_{x}D_tudx}_{\star}\\
\star &=\frac{\mu}{2}\frac{d}{dt}\li\chi\rho r^4|D_xu|^2dx-
\frac{\mu}{2}\li\chi D_t\rho r^4|D_xu|^2dx- 2\mu\li\chi \rho r^3
u|D_xu|^2dx
\end{split}
\]
Recall $\sup_{x_0\leq x\leq 1}|\frac{D_t\rho}{\rho}|\leq K$  and
$\sup_{x_0\leq x\leq 1}|\frac{u}{r}|\leq K$ from the assumption
(\ref{K}). Thus we get
\[
\star\leq \frac{\mu}{2}\frac{d}{dt}\li\chi\rho r^4|D_xu|^2dx
+\frac{5\mu K}{2}\li\chi \rho r^4|D_xu|^2dx.
\]
The fourth term: first apply the product rule in $t$-variable and
use the assumption (\ref{K})
\[
\begin{split}
\mu\li\chi\frac{2uD_t u}{\rho
r^2}dx&=\frac{\mu}{2}\frac{d}{dt}\li\chi \frac{2u^2}{\rho
r^2}dx+\mu\li\chi\frac{D_t\rho u^2}{\rho^2 r^2}dx
+2\mu\li\chi\frac{u^3}{\rho r^3}dx \\
&\leq\frac{\mu}{2}\frac{d}{dt}\li\chi \frac{2u^2}{\rho
r^2}dx+\frac{3\mu K}{2}\li\chi \frac{2u^2}{\rho r^2}dx
\end{split}
\]
The fifth term: apply the Cauchy-Schwarz inequality
\[
\begin{split}
\li\chi\frac{4\pi x}{r^2}D_tudx&\leq \frac{\mu}{4}\li\chi
\frac{|D_tu|^2}{\rho r^2}dx +\frac{1}{\mu}\li\chi\frac{16\pi^2
x^2\rho}{r^2}dx\\
&\leq\frac{\mu}{4}\li\chi \frac{|D_tu|^2}{\rho r^2}dx +\frac{C}{
(r_0-d)^2}\sup_{x_0\leq x\leq 1}|\rho^{2-\gamma}|\li\chi
\rho^{\gamma-1}dx
\end{split}
\]
This completes the inequality (\ref{onex}) for $i=0$. For $i=1,2$,
multiply (\ref{mom}) by $\chi D_t^{i+1}u$ and integrate to get
\[
\begin{split}
\li\chi |D_t^{i+1}u|^2dx +\li\chi r^2D_xD_t^ip
D_tu^{i+1}dx-\mu\li\chi D_x(\rho r^4D_x
D_t^iu)D_tu^{i+1}dx\\+\mu\li\chi\frac{2D_t^iuD_tu^{i+1}}{\rho
r^2}dx=-\li\chi D_t^i(\frac{4\pi x}{r^2})D_t^{i+1}udx\\
+\sum_{j=0}^{i-1}\li\chi\{\mu D_x(D_t^{i-j}(\rho
r^4)D_t^jD_xu)-D_t^{i-j}r^2 D_xD_t^jp-D_t^{i-j}(\frac{2\mu}{\rho
r^2})D_t^ju\}D_t^{i+1}udx
\end{split}
\]
The LHS has the same structure as the case $i=0$ and the RHS
consists of lower order terms. The computation is similar as in the
previous lemma. Hence, by applying the Cauchy-Schwarz inequality and
using the assumption $(K)$, we get the desired results.
\end{proof}\

Before we move onto higher spatial derivatives of $u$, we present
how to estimate $D_x\rho$, a crucial step to close the energy
estimates. Note that for stationary solutions, due to
(\ref{behaviorl}), $$|D_x\rho_0^{\frac{\gamma}{2}}|^2\sim
\frac{1}{1-x}$$ which is not integrable, and hence weighted
estimates come into play. In order to estimate spatial derivatives
of $\rho$, we make use of both the continuity equation and an
integral form of the momentum equation. We point out that the
estimate
does not yield overlapping terms. \\

\begin{lemma}
There exists $C_K>0$ such that
\begin{equation}
\frac{1}{2}\frac{d}{dt}\sum_{i=0}^2\li\chi
\rho^{2\gamma-2}r^4|D_t^iD_x\rho|^2dx\leq
C_K\mathcal{E}_L.\label{rhox}
\end{equation}\label{rhoxl}
\end{lemma}

\begin{proof}
 Integrate the momentum equation in $x$ from 1 to $x$ for $x\geq x_0$: due to the
 boundary condition (\ref{bcl}),
\begin{equation}
-\rho r^2D_xu=-\frac{1}{\mu}p-
\frac{1}{\mu}\int_1^{x}\{\frac{D_tu}{r^2}+\frac{4\pi y}{r^4}+\mu
\frac{2u}{\rho r^4}-\mu\frac{2D_xu}{r}\} dy.\label{m2}
\end{equation}
We will estimate $\rho^{\gamma}$ rather than $\rho$. Multiply by
$\gamma \rho^{\gamma}$ and use the continuity equation:
\[
\gamma\rho^{\gamma}\frac{D_t\rho}{\rho}=-\frac{\gamma}{\mu}\rho^{\gamma}
p-\gamma\rho^{\gamma}\frac{2u}{r}-
\frac{\gamma}{\mu}\rho^{\gamma}\int_1^{x}\{\frac{D_tu}{r^2}+\frac{4\pi
y}{r^4}+\mu \frac{2u}{\rho r^4}-\mu\frac{2D_xu}{r}\} dy
\]
Apply $D_x$:
\[
\begin{split}
D_tD_x(\rho^{\gamma})=-\frac{2A\gamma}{\mu}\rho^{\gamma}D_x(\rho^{\gamma})
-\gamma D_x(\rho^{\gamma})\frac{2u}{r}-\gamma\rho^{\gamma}
D_x(\frac{2u}{r})-\frac{\gamma}{\mu}\rho^{\gamma}(\frac{D_tu}{r^2}+\frac{4\pi
x}{r^4}+\mu \frac{2u}{\rho
r^4}-\mu\frac{2D_xu}{r})\\-\frac{\gamma}{\mu}D_x(\rho^{\gamma})
\int_1^{x}\{\frac{D_tu}{r^2}+\frac{4\pi y}{r^4}+\mu \frac{2u}{\rho
r^4}-\mu\frac{2D_xu}{r}\} dy
\end{split}
\]
Use the identity $D_x(\frac{u}{r})=\frac{D_xu}{r}-\frac{u}{\rho
r^4}$ and (\ref{m2}) in order to simplify the above as the
following:
\[
D_tD_x(\rho^{\gamma})=-\frac{A\gamma}{\mu}\rho^{\gamma}D_x(\rho^{\gamma})-
\gamma D_x(\rho^{\gamma})(\rho r^2
D_xu+\frac{2u}{r})-\frac{\gamma}{\mu}\rho^{\gamma}(\frac{D_tu}{r^2}+\frac{4\pi
x}{r^4})
\]
 Multiply by $r^4D_x(\rho^{\gamma})$, and integrate to get:
\begin{equation}
\begin{split}
&\frac{1}{2}\frac{d}{dt}\li \chi r^4|D_x(\rho^{\gamma})|^2
dx-2\li\chi r^3u|D_x(\rho^{\gamma})|^2dx+\frac{A\gamma}{\mu}\li\chi
\rho^{\gamma}r^4|D_x(\rho^{\gamma})|^2 dx\\
&=-\gamma\li\chi(\rho r^2D_xu+
\frac{2u}{r})r^4|D_x(\rho^{\gamma})|^2dx
-\frac{\gamma}{\mu}\li\chi\rho^{\gamma}(D_tu+\frac{4\pi x}{r^2})
r^2D_x(\rho^{\gamma})dx\label{rrr}
\end{split}
\end{equation}
The last term can be bounded by
\[
\frac{\gamma}{\mu}\li\chi\rho^{\gamma}(D_tu+\frac{4\pi x}{r^2})
r^2D_x(\rho^{\gamma})dx\leq\frac{A\gamma}{2\mu}\li\chi \rho^{\gamma}
r^4 |D_x(\rho^{\gamma})|^2dx+\frac{1}{A\mu }\li\chi\rho^{\gamma}
(|D_tu|^2+\frac{16\pi^2x^2}{r^4})dx
\]
by the Cauchy-Schwarz inequality. Hence by the assumption (\ref{K}),
(\ref{rrr}) becomes
\[
\begin{split}
 &\frac{1}{2}\frac{d}{dt}\li \chi r^4|D_x(\rho^{\gamma})|^2
dx+\frac{A\gamma}{2\mu}\li\chi
\rho^{\gamma}r^4|D_x(\rho^{\gamma})|^2 dx\\
&\leq (2+\gamma)K\li\chi r^4|D_x(\rho^{\gamma})|^2dx+C\sup_{x_0\leq
x\leq 1}|\rho^{\gamma}|\li\chi |D_tu|^2dx+C\sup_{x_0\leq x\leq
1}|\rho|\li\chi \rho^{\gamma-1}dx
\end{split}
\]
This proves (\ref{rhox}) when $i=0$. Note that taking $t$
derivatives does not destroy the structure of equations and thus, by
using the previous lemmas for $t$ derivatives of $u$, one can derive
the similar estimates for $i=1,2.$
\end{proof}\

As a corollary of Lemma \ref{rhoxl}, we derive the estimates of
mixed derivatives of $u$ with two spatial derivatives, which will be
needed to estimate higher order derivatives of $\rho$ and to close
the energy estimates later in Section \ref{section5}, directly from
the equation. The $L^2$ estimate of $D_x(\rho r^2D_xu+\frac{2u}{r})$
easily follows from the momentum equation: square it and integrate
to get
\begin{equation}
\begin{split}
\li\chi \rho|r^2D_x(\rho r^2D_xu+\frac{2u}{r})|^2dx&\leq
\frac{2}{\mu^2}\li\chi\{
\rho|D_tu|^2+\rho r^4|D_xp|^2+\frac{16\pi^2x^2\rho}{r^4}\}dx\\
&\leq C_K \mathcal{E}_L\label{2s}
\end{split}
\end{equation}
Each term in the RHS has been already estimated and it is bounded by
$\mathcal{E}_L$. Next, in order to estimate $D_x(\rho r^2 D_xu)$, we
rewrite the momentum equation:
\[
\mu D_x(\rho r^2D_xu)=\frac{D_tu}{r^2}+D_xp+\frac{4\pi x}{r^4}
+{\mu\frac{2u}{\rho r^4}-\mu\frac{2D_xu}{r}}
\]
Note that due to the singularity of the last two terms at the vacuum
boundary in RHS we need an appropriate weight for $L^2$ estimate of
$D_x(\rho r^2D_xu)$: in the same spirit as in (\ref{2s}), we obtain
\begin{equation}
\begin{split}
\frac{\mu^2}{2}\li\chi\rho|r^2D_x(\rho r^2D_xu)|^2dx&\leq
\li\chi\{\rho|D_tu|^2+\rho r^4|D_xp|^2+\frac{16\pi^2
x^2\rho}{r^4}\\&\;\;\;\;\;\;\;\;\;\;\;\;\;+\frac{4\mu^2 u^2}{\rho
r^4}+4\mu^2\rho r^2|D_xu|^2\} dx\\
&\leq C_K\mathcal{E}_L+ \frac{4\mu^2}{(r_0-d)^2}\li\chi\{\rho
r^4|D_xu|^2+\frac{2u^2}{\rho r^2}\}dx
\end{split}
\end{equation}
Similarly, one can derive the following: for $i=1,2$
\begin{equation}
\begin{split}
\li\chi \rho|r^2D_x(\rho r^2D_t^iD_xu+\frac{2D_t^iu}{r})|^2dx\leq
C_K
\mathcal{E}_L\\
\li\chi\rho|r^2D_x(\rho r^2D_t^iD_xu)|^2dx\leq
(C_K+C_{r_0})\mathcal{E}_L
\end{split}
\end{equation}
where $C_{r_0}$ is a constant depending on $r_0$. But this
dependence on $r_0$ can be ignored because we have chosen $r_0$ and
$d$ so that $\frac{1}{r_0-d}\leq 1$.

Next, we present the weighted energy estimate of $D_x^2\rho$.\\

\begin{lemma}
There exists $C_K>0$ such that
\begin{equation}
\frac{1}{2}\frac{d}{dt}\sum_{i=0}^1\li\chi\rho^{4\gamma-2}r^8|D_t^iD_x^2\rho|^2
dx\leq C_K\mathcal{E}_L.\label{rhoxdl}
\end{equation}\label{rhoxdll}
\end{lemma}

\begin{proof}
Take $D_x$ of the continuity equation and use the momentum equation
to substitute $D_x(\rho r^2D_xu+\frac{2u}{r})$ by
$\frac{1}{\mu}(\frac{D_t u}{r^2}
+A\gamma\rho^{\gamma-1}D_x\rho+\frac{4\pi x}{r^4})$:
\begin{equation}
D_{t}D_x\rho+D_x\rho(\rho
r^2D_xu+\frac{2u}{r})+\frac{\rho}{\mu}(\frac{D_t u}{r^2}
+A\gamma\rho^{\gamma-1}D_x\rho+\frac{4\pi x}{r^4})=0 \label{rhotx2}
\end{equation}
Take one more $D_x$ of (\ref{rhotx2}), use the momentum equation
again to substitute $D_x(\rho r^2D_xu+\frac{2u}{r})$ by
$\frac{1}{\mu}(\frac{D_t u}{r^2}
+A\gamma\rho^{\gamma-1}D_x\rho+\frac{4\pi x}{r^4})$, and rearrange
terms to get:
\begin{equation}
\begin{split}
D_tD_x^2\rho+\frac{A\gamma}{\mu}\rho^{\gamma}D_{x}^2\rho
=&-D_{x}^2\rho(\rho r^2D_x
u+\frac{2u}{r})-\frac{A\gamma(\gamma+1)}{\mu}\rho^{\gamma-1}|D_x\rho|^2
\\
&-\frac{2D_x\rho}{\mu}(\frac{D_tu}{r^2}+\frac{4\pi x}{r^4})
-\frac{\rho}{\mu}D_x(\frac{D_tu}{r^2}+\frac{4\pi x}{r^4})
\label{rhoxx}
\end{split}
\end{equation}
Multiply by $\chi\rho^{4\gamma-2}r^8D_{x}^2\rho$ and integrate in
$x$ to get
\[
\begin{split}
&\frac{1}{2}\frac{d}{dt}\li\chi\rho^{4\gamma-2}r^8|D_{x}^2\rho|^2dx
-\underline{\frac{1}{2}\li\chi
D_t(\rho^{4\gamma-2}r^8)|D_{x}^2\rho|^2dx}_{(i)}+
\frac{A\gamma}{\mu}\li\chi\rho^{5\gamma-2}|D_{x}^2\rho|^2dx\\
&=\underline{-\li\chi(\rho
r^2D_xu+\frac{2u}{r})\rho^{4\gamma-2}r^8|D_{x}^2\rho|^2dx}_{(ii)}
\underline{-\frac{A\gamma(\gamma+1)}{\mu}\li\chi
\rho^{5\gamma-3}r^8|D_x\rho|^2 D_{x}^2\rho dx}_{(iii)} \\&\;\;\;\;
\underline{-\frac{2}{\mu}\li\chi\rho^{4\gamma-2}r^8
D_x\rho(\frac{D_tu}{r^2}+\frac{4\pi x}{r^4})D_{x}^2\rho
dx}_{(iv)}\\
&\;\;\;\;\underline{-\frac{1}{\mu}
\li\chi\rho^{4\gamma-1}r^8(\frac{D_tD_xu}{r^2}-\frac{2D_tu}{\rho
r^5}+\frac{4\pi }{r^4}-\frac{16\pi x}{\rho r^7}) D_{x}^2\rho
dx}_{(v)}
\end{split}
\]
For $(i)$ and $(ii)$, we use the assumption (\ref{K}):
\[
\begin{split}
(i)&=(2\gamma-1)\li\chi D_t\rho
\rho^{4\gamma-3}r^8|D_{x}^2\rho|^2dx+
4\li\chi\rho^{4\gamma-2}r^7u|D_{x}^2\rho|^2dx\\
&\leq (2\gamma-1)\sup_{x_0\leq x\leq 1}|\frac{D_t\rho}{\rho}|
\li\chi\rho^{4\gamma-2}r^8|D_{x}^2\rho|^2dx+4\sup_{x_0\leq x\leq 1}
|\frac{u}{r}|\li\chi\rho^{4\gamma-2}r^8|D_{x}^2\rho|^2dx\\
&\leq CK\li\chi\rho^{4\gamma-2}r^8|D_{x}^2\rho|^2dx\\
(ii)&\leq \sup_{x_0\leq x\leq 1}|\rho
r^2D_xu+\frac{2u}{r}|\li\chi\rho^{4\gamma-2}r^8|D_{x}^2\rho|^2dx
\leq K\li\chi\rho^{4\gamma-2}r^8|D_{x}^2\rho|^2dx
\end{split}
\]
For $(iii)$ and $(iv)$, we use the assumption (\ref{K}) and apply
the Cauchy-Schwarz inequality:
\[
\begin{split}
&(iii)\leq C\sup_{x_0\leq x\leq 1}|\rho^{2\gamma-1}r^2D_x\rho|
(\li\chi\rho^{2\gamma-2}r^4|D_x\rho|^2dx
+\li\chi\rho^{4\gamma-2}r^8|D_{x}^2\rho|^2dx)\\
&(iv)\leq C\sup_{x_0\leq x\leq
1}|\rho^{2\gamma-1}r^2D_x\rho|(\li\chi(|D_tu|^2+\frac{16\pi^2
x^2}{r^4})dx +\li\chi\rho^{4\gamma-2}r^8|D_{x}^2\rho|^2dx)
\end{split}
\]
For $(v)$, we separate terms into two cases. Apply the
Cauchy-Schwarz inequality to the first two terms after taking the
sup of $\rho^{2\gamma-\frac{1}{2}}$:
\[
\begin{split}
\li\chi&\rho^{4\gamma-1}r^8(\frac{D_tD_xu}{r^2}-\frac{2D_tu}{\rho
r^5})D_{x}^2\rho dx\\
&\leq\sup_{x_0\leq x\leq 1}|\rho^{2\gamma-\frac{1}{2}}|(\li\chi
(\rho r^4|D_tD_xu|^2+\frac{2|D_tu|^2}{\rho
r^2})dx+\li\chi\rho^{4\gamma-2}r^8|D_{x}^2\rho|^2dx)
\end{split}
\]
For the last two terms, before applying the Cauchy-Schwarz
inequality, we take the sup of different powers of $\rho$ in order
to have them bounded by $\mathcal{E}_L$:
\[
\begin{split}
4\pi\li\chi\rho^{4\gamma-1}r^4 D_{x}^2\rho dx&\leq 2\pi\sup_{x_0\leq
x\leq 1}|\rho^{\frac{3\gamma+1}{2}}|
(\li\chi\rho^{\gamma-1}dx+\li\chi\rho^{4\gamma-2}r^8|D_{x}^2\rho|^2dx)\\
-16\pi\li\chi\rho^{4\gamma-2}rx D_{x}^2\rho dx&\leq 8\pi
\sup_{x_0\leq x\leq 1}|\rho^{\frac{3\gamma-1}{2}}|
(\frac{1}{(r_0-d)^6}\li\chi\rho^{\gamma-1}dx+\li\chi\rho^{4\gamma-2}r^8
|D_{x}^2\rho|^2dx)
\end{split}
\]
This completes (\ref{rhoxdl}) for $i=0$. One can derive the similar
inequality when $i=1$.
\end{proof}\

As a direct result of the previous lemmas, we can get the bounds of
mixed derivatives of $u$ with three spatial derivatives. Multiply
the momentum equation by $\rho r^2$ and take  $D_x$ to get
\[
\begin{split}
\mu D_x(\rho r^2D_x(\rho r^2D_xu))=&\rho D_tD_xu+D_x\rho D_tu+
A\gamma\rho^{\gamma}
r^2D_x^2\rho+A\gamma^2\rho^{\gamma-1}r^2|D_x\rho|^2+\frac{2A\gamma\rho
^{\gamma-1}D_x\rho}{r}\\
&+\frac{4\pi \rho}{r^2}+\frac{4\pi xD_x\rho}{r^2}-\frac{8\pi x}{r^5}
+\frac{4\mu D_xu}{r^2}-\frac{4\mu u}{\rho r^5}-\frac{2\mu D_x(\rho
r^2D_xu)}{r}.
\end{split}
\]
Each term in RHS has been estimated with suitable weights, and thus
one can check
\begin{equation}
\mu^2\li\chi |\rho r^2D_x(\rho r^2D_x(\rho r^2D_xu))|^2dx\leq C_K
\mathcal{E}_L.
\end{equation}
Similarly, one gets
\[
\mu^2\li\chi |\rho r^2D_x(\rho r^2D_x(\rho r^2D_tD_xu))|^2dx\leq C_K
\mathcal{E}_L.
\]

To complete the estimate for $\mathcal{E}_L$, it remains to estimate
pure spatial derivative terms $D_x^3\rho$ and $D_x^4u$. The spirit
is same as in Lemma \ref{rhoxdll}: take $D_x$ of (\ref{rhoxx}) get
\[
\begin{split}
D_tD_x^3\rho+\frac{A\gamma}{\mu}\rho^{\gamma}D_{x}^3\rho
=&-D_{x}^3\rho(\rho r^2D_x
u+\frac{2u}{r})-\frac{3A\gamma(\gamma+1)}{\mu}\rho^{\gamma-1}D_x\rho
D_x^2\rho-\frac{A\gamma(\gamma^2-1)}{\mu}\underline{\rho^{\gamma-2}(D_x\rho)^3}_{(i)}
\\
&-\frac{3D_x^2\rho}{\mu}(\frac{D_tu}{r^2}+\frac{4\pi
x}{r^4})-\frac{3D_x\rho}{\mu}D_x(\frac{D_tu}{r^2}+\frac{4\pi
x}{r^4})
-\underline{\frac{\rho}{\mu}D_x^2(\frac{D_tu}{r^2}+\frac{4\pi
x}{r^4})}_{(ii)}
\end{split}
\]
Each term in RHS has been already treated with appropriate weights
and thus by multiplying by $\rho^{8\gamma-2}r^{12}D_x^3\rho$ and
integrating, we conclude the following:
\begin{equation}
\frac{1}{2}\frac{d}{dt}\li\chi\rho^{8\gamma-2}r^{12}|D_x^3\rho|^2dx\leq
C_K\mathcal{E}_L
\end{equation}
Note that the weight $\rho^{8\gamma-2}$ is carefully chosen such
that the RHS can be bounded by $C_K\mathcal{E}_L$. For instance,
$(i),\;(ii)$ can be estimated as follows. Indeed, the estimate of
$(i)$ will provide a good reason for the choice of weights. First,
for $(i)$, we take the sup of $\rho^{2\gamma-1}r^2 D_x\rho$ and
apply the Cauchy-Schwarz inequality.
\[
\begin{split}
\li\chi&\rho^{8\gamma-2}r^{12}\rho^{\gamma-2}(D_x\rho)^3  D_x^3\rho
dx=\li\chi (\rho^{2\gamma-1}r^2
D_x\rho)^2(\rho^{\gamma-1}r^2D_x\rho)(\rho^{4\gamma-1}r^6D_x^3\rho)dx\\
&\leq \sup_{x_0\leq x\leq 1}|\rho^{2\gamma-1}r^2
D_x\rho|^2(\li\chi\rho^{2\gamma-2}r^4|D_x\rho|^2dx
+\li\chi\rho^{8\gamma-2}r^{12}|D_x^3\rho|^2dx)\leq K^2\mathcal{E}_L
\end{split}
\]
For $(ii)$, first compute the second derivative of
$\frac{D_tu}{r^2}+\frac{4\pi x}{r^4}$.
\[
\begin{split}
\li\chi \rho^{8\gamma-2}r^{12}\rho D_x^2(\frac{D_tu}{r^2}+\frac{4\pi
x}{r^4})D_x^3\rho dx=\li\chi\rho^{4\gamma-1}\{r^2D_x(\rho
r^2D_tD_xu)-(\rho r^2D_{t}D_xu)(\frac{D_x\rho
r^2}{\rho}+\frac{6}{\rho r})\\+2D_tu(\frac{D_x\rho
r}{\rho}+\frac{5}{\rho r^2})-\frac{32\pi}{r}+16\pi
x(\frac{D_x\rho}{\rho r}+\frac{7}{\rho
r^4})\}(\rho^{4\gamma-1}r^6D_x^3\rho) dx
\end{split}
\]
Each term in the RHS has been already estimated. We repeat the same
argument, namely, extract an appropriate factor of $\rho$ possibly
with some other terms outside the integral to make use of the
assumption $(K)$, and apply the Cauchy-Schwarz inequality to
conclude that it is bounded by $C_K\mathcal{E}_L$ for some constant
$C_K$ depending on $K$. In turn, from the momentum equation, a
weighted $L^2$ estimate of $D_x(\rho r^2D_x(\rho r^2D_x(\rho
r^2D_xu)))$ can be obtained. This
completes the estimate of $\mathcal{E}_L(t)$.\\

\section{Interior Estimates in Eulerian Coordinates}
\label{section4}

In this section, we present the interior estimates of
$\mathcal{E}_E(t)$. Away from the vacuum boundary, $\rho$ is
expected to be strictly positive and the classical results of the
Navier-Stokes theory can be applied. The potential terms need
attention. Interior domain containing the origin retains
three-dimensional structure despite the symmetry and $H^3$ interior
estimate is necessary. The Eulerian description is used to establish
the interior energy estimates. The energy estimates will be
performed in the cartesian product space in order to avoid the
coordinate singularity at the origin coming from the symmetry.
Recall the full Navier-Stokes-Poisson system (\ref{NSP}):
\begin{equation*}
\begin{split}
\rho_t+\nabla\cdot(\rho \mathbf{u})=0\\
\rho(\mathbf{u}_t+(\mathbf{u}\cdot \nabla)\mathbf{u})+\nabla
p+\rho\nabla \Phi=\mu\triangle
\mathbf{u}\\
\triangle\Phi=4\pi\rho\label{nspeg}
\end{split}
\end{equation*}
where $\mathbf{u}(\mathbf{x})=u(r)\frac{\mathbf{x}}{r},$ where
$\mathbf{x}=(x_1,x_2,x_3)$, $r=|\mathbf{x}|$.

We will not put the potential term into conservation form in the
energy estimates because it makes the energy negative definite and
in principle, it serves lower order term and it can be controlled by
main terms $\rho, u$. The symmetrizer of the system is used for the
energy estimates to get cleaner estimates,
while other weight functions may work by the aid of the dissipation.\\

\begin{lemma} There exist constants $ C_K, C>0$ such that
\begin{equation}
\frac{1}{2}\frac{d}{dt}\mathcal{E}_E+\frac{1}{2}\mathcal{D}_E\leq
C_K\mathcal{E}_E+C(\mathcal{E}_E)^2+\mathcal{OL}_3
\end{equation}
where $\mathcal{OL}_3\leq \widetilde{C}_K\mathcal{E}_L$ for some
$\widetilde{C}_K$.
\end{lemma}

\begin{proof} It is standard to perform the energy estimates to the
Navier-Stokes equations, and so we point out the gist of the
estimates rather than provide all the details.
 First, the zeroth order energy estimate is carried out.
Let $W_0\equiv\frac{A\gamma}{\gamma-1}\rho^{\gamma-1}$ be a
symmetrizing weight function.  Multiply (\ref{NSP}) by $\zeta W_0$
and $\zeta \mathbf{u}$, and integrate to get:
\[
\begin{split}
\frac{d}{dt}\{\frac{1}{2}\int\zeta\rho
|\mathbf{u}|^2d\mathbf{x}+\frac{A}{\gamma-1}\int\zeta\rho^{\gamma}
d\mathbf{x}\}+\mu\int\zeta|\nabla
\mathbf{u}|^2d\mathbf{x}\\=\frac{1}{4\pi}\int\zeta\nabla\Phi\cdot\nabla
\partial_t\Phi d\mathbf{x}
+\frac{1}{2}\int\zeta\partial_t\rho|\mathbf{u}|^2d\mathbf{x}
-\int\zeta\rho(\mathbf{u}\cdot\nabla)\mathbf{u}\cdot
\mathbf{u}d\mathbf{x}\\+\frac{A\gamma}{\gamma-1}\int\nabla
\zeta\cdot
\rho^{\gamma}\mathbf{u}d\mathbf{x}+\mu\int\nabla\zeta\cdot\nabla
\mathbf{u} \mathbf{u}d\mathbf{x} +\int\nabla\zeta\Phi\rho
\mathbf{u}d\mathbf{x}+\frac{1}{4\pi}\int\nabla\zeta\Phi\nabla\Phi_t
d\mathbf{x}
\end{split}
\]
The RHS consists of the first potential term,  the next two
nonlinear terms, and remaining overlapping terms. Here is the
estimate of the potential term: by symmetry
\[
\begin{split}
|\frac{1}{4\pi}\int\zeta\nabla\Phi\cdot\nabla\partial_t\Phi
d\mathbf{x}|&=|\int_0^{r_2-d}\zeta(\frac{1}{r^3}\int_0^r\rho s^2ds)
(\frac{1}{r^3}\int_0^r\partial_t\rho s^2ds)r^2dr|\\&=|-
\int_0^{r_2-d}\zeta(\frac{1}{r^3}\int_0^r\rho s^2ds) (\frac{\rho
u}{r})r^2dr|\\
&\leq \frac{1}{\mu}\int_0^{r_2-d}\zeta(\frac{1}{r^3}\int_0^r\rho
s^2ds)^2 \rho^2 r^2dr +\frac{\mu}{4}\int\zeta|\frac{u}{r}|^2r^2dr\\
&\leq \frac{1}{9\mu}\sup_{0\leq r\leq
r_2-d}|\rho^{4-\gamma}|\int\zeta\rho^{\gamma}r^2dr
+\frac{\mu}{4}\int\zeta|\frac{u}{r}|^2r^2dr
\end{split}
\]
We have applied the Cauchy-Schwarz inequality and have used for
$r\leq r_2-d$,
$$\frac{1}{r^3}\int_0^r\rho s^2ds\leq \frac{1}{3}\sup_{0\leq r\leq
r_2-d}|\rho|.$$ Next, we compute overlapping terms in Lagrangian
coordinates. To illustrate the idea, we present the estimate of the
first two terms: by change of variables, Eulerian integrals are
converted into Lagrangian integrals and then we apply the
Cauchy-Schwarz inequality.
\[
\begin{split}
|\int\nabla \zeta\cdot
\rho^{\gamma}\mathbf{u}d\mathbf{x}|&\leq\frac{C}{r_2-r_1-2d}
\int_{x_1}^{x_2}\rho^{\gamma-1}|u|dx\\
&\leq \frac{C} {r_2-r_1-2d}\sup_{x_1\leq x\leq
x_2}|\rho^{\frac{\gamma-1}{2}}|\{\int_{x_1}^{x_2}\rho^{\gamma-1} dx+
\int_{x_1}^{x_2}|u|^2 dx\}\\
|\mu\int\nabla\zeta\cdot\nabla \mathbf{u}
\mathbf{u}d\mathbf{x}|&=|-\frac{\mu}{2}\int\triangle\zeta
|\mathbf{u}|^2d\mathbf{x}|
\\&\leq C\mu
\int_{x_1}^{x_2} u^2\frac{dx}{\rho} \leq C\mu\sup_{x_1\leq x\leq
x_2}|\frac{1}{\rho}|\int_{x_1}^{x_2}|u|^2 dx
\end{split}
\]
Note that
\[
\sup_{x_1\leq x\leq x_2}|\frac{1}{\rho}|\leq \sup_{x_1\leq x\leq
x_2}|\frac{1}{\rho_{in}}|e^{MT}\text{ from }(\ref{rhob}).
\]
 Hence we get the following zeroth-order estimate:
\begin{equation}
\frac{1}{2}\frac{d}{dt}\{\int\zeta\rho
|\mathbf{u}|^2d\mathbf{x}+\frac{A}{\gamma-1}\int\zeta\rho^{\gamma}d\mathbf{x}\}
+\frac{3\mu}{4}\int\zeta|\nabla \mathbf{u}|^2d\mathbf{x}\leq
C_K\mathcal{E}_E+\mathcal{OL}\label{e0}
\end{equation}
where $\mathcal{OL}\leq C_{K,in}\mathcal{E}_L$ and $C_{K,in}$
depends also on the initial density $\rho_{in}$. The higher temporal
derivatives (up to third) can be estimated in the similar way. Let
$\partial$ be any Eulerian derivative.
\[
\begin{split}
\partial_t(\partial\rho)+\nabla\partial\rho\cdot
\mathbf{u}+\nabla\rho\cdot\partial
\mathbf{u}+\partial\rho\nabla\cdot \mathbf{u}
+\rho\nabla\cdot\partial \mathbf{u}&=0\\
\rho\partial_t(\partial \mathbf{u})+\partial \rho
\partial_t\mathbf{u}+\partial[\rho(\mathbf{u}\cdot\nabla)\mathbf{u}]+\nabla\partial
p+\rho\nabla
\partial\Phi+\partial\rho\nabla\Phi&=\mu\triangle\partial \mathbf{u}
\end{split}
\]
Define $W_1\equiv A\gamma\rho^{\gamma-2}\partial\rho$.
\[
\begin{split}
\frac{1}{2}\frac{d}{dt}\{A\gamma\int\zeta\rho^{\gamma-2}(\partial\rho)^2
d\mathbf{x} +\int\zeta\rho|\partial \mathbf{u}
|^2d\mathbf{x}\}+\mu\int\zeta|\nabla\partial
\mathbf{u}|^2d\mathbf{x}
+\mu\int\nabla\zeta\nabla\partial \mathbf{u}\partial \mathbf{u}d\mathbf{x}\\
=\frac{A\gamma(\gamma-2)}{2}\int\zeta\rho^{\gamma-3}\partial_t\rho(\partial\rho)^2
d\mathbf{x} +\frac{1}{2}\int\zeta\partial_t\rho|\partial
\mathbf{u}|^2d\mathbf{x}-A\gamma\int\zeta\rho^{\gamma-2}\partial\rho
\nabla\partial\rho\cdot \mathbf{u}d\mathbf{x}\\
-A\gamma\int\zeta\rho^{\gamma-2}\partial\rho\nabla\rho\cdot \partial
\mathbf{u}d\mathbf{x}-A\gamma\int
\zeta\rho^{\gamma-2}(\partial\rho)^2\nabla\cdot
\mathbf{u}d\mathbf{x}-A\gamma\int\zeta\rho^{\gamma-1}
\partial\rho\nabla\cdot\partial \mathbf{u}d\mathbf{x}\\
-\int\zeta\partial\rho \partial_t\mathbf{u}\partial
\mathbf{u}d\mathbf{x}-\int\zeta
\partial[\rho(\mathbf{u}\cdot\nabla)\mathbf{u}] \partial \mathbf{u}d\mathbf{x}
-\int\zeta \nabla\partial p\partial
\mathbf{u}d\mathbf{x}-\int\zeta\rho\nabla\partial\Phi\partial
\mathbf{u}d\mathbf{x}-\int\zeta\partial\rho \nabla\Phi\partial
\mathbf{u}d\mathbf{x}
\end{split}
\]
Note that the pressure term $-\int\zeta\nabla\partial p\partial
\mathbf{u}d\mathbf{x}$ on the RHS contains the undesirable second
derivative of $\rho$ but it is canceled with
$-A\gamma\int\zeta\rho^{\gamma-1}\partial\rho\nabla\cdot\partial
\mathbf{u}d\mathbf{x}$ due to the choice of symmetrizing weight
$W_1$:
\[
-A\gamma\int\zeta\rho^{\gamma-1}\partial\rho\nabla\cdot\partial
\mathbf{u}d\mathbf{x} -\int\zeta\nabla\partial p\partial
\mathbf{u}d\mathbf{x}=-\int\zeta\partial p\nabla\cdot\partial
\mathbf{u}d\mathbf{x}- \int\zeta\nabla\partial p\partial
\mathbf{u}d\mathbf{x}=\int\nabla\zeta\partial p\cdot\partial
\mathbf{u}d\mathbf{x}
\]
For another second derivative term
$-A\gamma\int\zeta\rho^{\gamma-2}\partial\rho
\nabla\partial\rho\cdot \mathbf{u}d\mathbf{x}$, we integrate it by
parts:
\[
\frac{A\gamma}{2}\int\nabla\zeta\rho^{\gamma-2}(\partial\rho)^2\mathbf{u}d\mathbf{x}+
\frac{A\gamma(\gamma-2)}{2}\int\zeta\rho^{\gamma-3}\nabla\rho(\partial\rho)^2\mathbf{u}
d\mathbf{x}
+\frac{A\gamma}{2}\int\zeta\rho^{\gamma-2}(\partial\rho)^2\nabla\cdot
\mathbf{u}d\mathbf{x}
\]
Potential terms are, in principle, lower order and the  $L^2$
estimate $||\partial^2\Phi||_{L^2}\leq C||\rho||_{L^2}$ is useful.
Hence we get the following first order estimate:
\begin{equation}
\frac{1}{2}\frac{d}{dt}\{A\gamma\int\zeta\rho^{\gamma-2}(\partial\rho)^2
d\mathbf{x} +\int\zeta\rho|\partial \mathbf{u}
|^2d\mathbf{x}\}+\frac{3\mu}{4}\int\zeta|\nabla\partial
\mathbf{u}|^2d\mathbf{x} \leq
C_K\mathcal{E}_E+\mathcal{OL}\label{e1}
\end{equation}
In order to apply the Sobolev imbedding theorem, we need to estimate
up to the 3rd derivatives. Take one more derivative of the
equations:
\[
\begin{split}
\partial_t(\partial^2\rho)+\nabla\partial^2\rho\cdot
\mathbf{u}+2\nabla\partial\rho\cdot\partial \mathbf{u}
+\nabla\rho\cdot\partial^2\mathbf{u}+\partial^2\rho\nabla\cdot
\mathbf{u}+2\partial\rho\nabla\cdot
\partial \mathbf{u}+\rho\nabla\cdot\partial^2\mathbf{u}=0\\
\rho\partial_t(\partial^2
\mathbf{u})+2\partial\rho\partial_t(\partial
\mathbf{u})+\partial^2\rho
\partial_t\mathbf{u}+\partial^2[\rho(\mathbf{u}\cdot
\nabla)\mathbf{u}]+\nabla\partial^2p+2\partial\rho\nabla\partial\Phi
+\rho\nabla\partial^2\Phi
+\partial^2\rho\nabla\Phi=\mu\triangle\partial^2\mathbf{u}
\end{split}
\]
Let $W_2\equiv A\gamma\rho^{\gamma-2}\partial^2\rho$.
\[
\begin{split}
\frac{1}{2}\frac{d}{dt}\{\int\zeta\rho^{\gamma-2}(\partial^2\rho)^2d\mathbf{x}
+\int\zeta\rho |\partial^2
\mathbf{u}|^2d\mathbf{x}\}+\mu\int\zeta|\nabla\partial^2\mathbf{u}|^2d\mathbf{x}
+\mu\int\nabla\zeta
\nabla\partial^2\mathbf{u}\partial^2\mathbf{u}d\mathbf{x}\\
=\frac{A\gamma(\gamma-2)}{2}\int\zeta\rho^{\gamma-3}\partial_t\rho(\partial^2\rho)^2
d\mathbf{x}+
\frac{1}{2}\int\zeta\partial_t\rho|\partial^2\mathbf{u}|^2d\mathbf{x}
-A\gamma\int\zeta\rho^{\gamma-2}\partial^2\rho
\nabla\partial^2\rho\cdot
\mathbf{u}d\mathbf{x}\\-2A\gamma\int\zeta\rho^{\gamma-2}
\partial^2\rho\nabla\partial\rho\cdot\partial \mathbf{u}d\mathbf{x}
-A\gamma\int\zeta\rho^{\gamma-2}
\partial^2\rho\nabla\rho\cdot\partial^2\mathbf{u}d\mathbf{x}
-A\gamma\int\zeta\rho^{\gamma-2}
(\partial^2\rho)^2\nabla\cdot \mathbf{u}d\mathbf{x}\\
-2A\gamma\int\zeta\rho^{\gamma-2}\partial^2\rho\partial\rho\nabla\cdot
\partial \mathbf{u}d\mathbf{x}
-A\gamma\int\zeta\rho^{\gamma-1}\partial^2\rho\nabla\cdot\partial^2\mathbf{u}d\mathbf{x}
-2\int\zeta
\partial\rho\partial_t(\partial \mathbf{u})\partial^2\mathbf{u}d\mathbf{x}\\
-\int\zeta\partial^2\rho
\partial_t\mathbf{u}\partial^2\mathbf{u}d\mathbf{x}
-\int\zeta\partial^2[\rho(\mathbf{u}\cdot\nabla)\mathbf{u}]\partial^2\mathbf{u}
d\mathbf{x}-\int\zeta\nabla\partial^2p\cdot\partial^2\mathbf{u}d\mathbf{x}
-2\int\zeta\partial\rho\nabla\partial\Phi\cdot\partial^2\mathbf{u}d\mathbf{x}\\
-\int\zeta\rho\nabla\partial^2\Phi
\cdot\partial^2\mathbf{u}d\mathbf{x}-\int\zeta\partial^2\rho\nabla\Phi\cdot\partial^2
\mathbf{u}d\mathbf{x}
\end{split}
\]
As in the first order estimates, for higher order derivative terms,
either we use the integration by parts or they cancel each other.
Eventually, we get the following:
\begin{equation}
\begin{split}
\frac{1}{2}\frac{d}{dt}\{A\gamma\int\zeta\rho^{\gamma-2}(\partial^2\rho)^2
d\mathbf{x} +\int\zeta\rho|\partial^2 \mathbf{u}
|^2d\mathbf{x}\}+\frac{3\mu}{4}\int\zeta|\nabla\partial^2
\mathbf{u}|^2 d\mathbf{x}\leq
C_K\mathcal{E}_E+\mathcal{OL}\label{e2}
\end{split}
\end{equation}

\noindent Take one more derivative and perform the weighted energy
estimates with $W_3=A\gamma\rho^{\gamma-2}\partial^3\rho$. It is
routine to have the following high energy inequality:
\begin{equation}
\frac{1}{2}\frac{d}{dt}\{A\gamma\int\zeta\rho^{\gamma-2}(\partial^3\rho)^2
d\mathbf{x} +\int\zeta\rho|\partial^3 \mathbf{u}
|^2d\mathbf{x}\}+\frac{\mu}{2}\int\zeta|\nabla\partial^3
\mathbf{u}|^2d\mathbf{x}\leq
C_K\mathcal{E}_E+C(\mathcal{E}_E)^2+\mathcal{OL} \label{e3}
\end{equation}
We remark that $C(\mathcal{E}_E)^2$ comes from the
Gagliardo-Nirenberg inequality: $$||f||_{L^4}\leq \frac{1}{2}
||\nabla f||^\frac{3}{4}_{L^2}|| f||^\frac{1}{4}_{L^2}$$ to treat
the nonlinear terms such as
 $\int\partial^2\rho \partial_t\partial
\mathbf{u}\partial^3\mathbf{u}d\mathbf{x}$. The application of the
Gagliardo-Nirenberg inequality appears later in the article. The
spirit is same, so we omit details. Note that a overlapping term
$\mu\int\nabla\zeta\nabla\partial^3\mathbf{u}\partial^3\mathbf{u}d\mathbf{x}$,
which has overfull derivative, can be also bounded by
$C\mathcal{E}_L$ because it can be reduced to
$-\frac{\mu}{2}\int\triangle\zeta|\partial^3\mathbf{u}|^2d\mathbf{x}$
after integration by parts; a $\rho$-weighted norm of $D_x(\rho
r^2D_x(\rho r^2D_xu))$ can be controlled by $C\mathcal{E}_L$; $\rho$
has a uniform upper bound as well as a uniform  lower bound
depending only on the initial profile in the overlapping interval.
This finishes the proof of the lemma.
\end{proof}\

\begin{proof}\textit{of Theorem \ref{prop}:} From the previous lemmas and the
energy inequalities, (\ref{ei}) follows. It remains to show that
$\mathcal{E}$ is bounded. We estimate $A$ explicitly by solving the
differential inequality (\ref{ei}). By separation of variables,
\[
\int \frac{d\mathcal{E}}{C_1\mathcal{E}+C_2(\mathcal{E})^2}\leq \int
dt\;\Rightarrow\;
\ln|\frac{\mathcal{E}(t)}{\frac{C_1}{C_2}+\mathcal{E}(t)}|\leq
\ln|\frac{\mathcal{E}(0)}{\frac{C_1}{C_2}+\mathcal{E}(0)}|+C_1t\;
\Rightarrow\;\frac{\mathcal{E}(t)}{\frac{C_1}{C_2}+\mathcal{E}(t)}\leq
\frac{\mathcal{E}(0)}{\frac{C_1}{C_2}+\mathcal{E}(0)}e^{C_1t}
\]
If $e^{C_1t}<\frac{C_1+C_2\mathcal{E}(0)}{C_2\mathcal{E}(0)}$, we
can derive the following:
\[
\mathcal{E}(t)\leq\frac{C_1\mathcal{E}(0)e^{C_1t}}{C_1-C_2\mathcal{E}(0)
(e^{C_1t}-1)}\equiv A
\]
Note that
$T<\frac{1}{C_1}\ln|\frac{C_1+C_2\mathcal{E}(0)}{C_2\mathcal{E}(0)}|$.
This completes the proof of Theorem \ref{prop}.
\end{proof}\

In the next section, we verify the assumption $(K)$ to show that the
energy estimates can be closed for $0\leq t\leq T$ where $T$ is
sufficiently small.\\

\section{Weaving the Estimates}
\label{section5}

The a priori estimates in the previous sections hold under the
assumption $(K)$: (\ref{K}). In order to close the estimates, it
remains to show that the assumption can be closed under the same
energy space $\mathcal{E}$. It results from the application of the
Sobolev embedding theorems for one dimension as well as three
dimension by combining boundary and interior estimates. We remark
that the positivity of $\rho$ having a lower and upper bounds in the
overlapping region is critical, which
results from the formula (\ref{rho1}). \\

\begin{lemma}
Suppose that $\rho, u, r$ serve a smooth solution to the
Navier-Stokes-Poisson system. Then there exist $T>0$,
$C=C(\rho_{in})>0$ such that $K\leq C\{\mathcal{E}(t)
^{\frac{1}{2}}+\mathcal{E}(t)\}$ for $0\leq t\leq T$.\label{weaving}
\end{lemma}

\begin{proof}
We provide the detailed computation for a few terms in the
assumption. Other terms can be verified in the same way. First of
all, we pay attention to $\rho$. For clear presentation, we use $M$,
instead of $K$, for the bound of $|\rho r^2 D_xu+\frac{2u}{r}|$:
$|\rho r^2 D_xu+\frac{2u}{r}|\leq M$.  In the view of (\ref{rho1}),
$\rho$ can be controlled by $M$:
\[
\rho_{in}e^{-MT}\leq \rho(t,x)\leq \rho_{in}e^{MT}
\]
With these bounds of $\rho$, let us  estimate $|\rho r^2
D_xu+\frac{2u}{r}|$.
\[
\sup_{0<x<1}|\rho r^2 D_xu+\frac{2u}{r}|=\sup\{\sup_{0\leq r\leq
r_1+d}|\partial_ru+\frac{2u}{r}|,\sup_{x_1\leq x\leq 1}|\rho r^2
D_xu+\frac{2u}{r}|\}
\]
Apply the Sobolev embedding theorem:
\[
\begin{split}
\bullet\sup_{0\leq r\leq r_1+d}|\partial_ru+\frac{2u}{r}|&\leq
\sum_{i=0}^2(\int_{B_{r_1+d}}|\nabla
\partial_{\mathbf{x}}^i\mathbf{u}|^2 d\mathbf{x})^{\frac{1}{2}}\\
&\leq \sup_{0\leq r\leq
r_1+d}|\frac{1}{\sqrt{\rho}}|\sum_{i=0}^2(\int_{B_{r_1+d}}\rho|\nabla
\partial_{\mathbf{x}}^i\mathbf{u}|^2 d\mathbf{x})^{\frac{1}{2}}\\
&\leq\sup_{0\leq r\leq
r_1+d}|\frac{1}{\sqrt{\rho_{in}}}|e^{\frac{MT}{2}}
(\mathcal{E}_E(t))^{\frac{1}{2}}
\end{split}
\]
In addition to the Sobolev embedding theorem, apply the
H$\ddot{o}$lder inequality:
\[
\begin{split}
\bullet\sup_{x_1\leq x\leq 1}|\rho r^2
D_xu+\frac{2u}{r}|&\leq\int_{x_1}^1|\rho r^2
D_xu+\frac{2u}{r}|dx+\int_{x_1}^1|D_x(\rho r^2
D_xu+\frac{2u}{r})|dx\\
&\leq (\int_{x_1}^1\rho dx)^{\frac{1}{2}}(\int_{x_1}^1\{\rho
r^4|D_xu|^2+\frac{2u^2}{\rho r^2}\} dx)^{\frac{1}{2}}\\
&\;\; +(\int_{x_1}^1\frac{1}{\rho r^4}
dx)^{\frac{1}{2}}(\int_{x_1}^1\rho|r^2D_x(\rho r^2
D_xu+\frac{2u}{r})|^2 dx)^{\frac{1}{2}}\\
&\leq \sup_{x_1\leq x\leq
1}|\rho_{in}^{\frac{2-\gamma}{2}}|e^{\frac{2-\gamma}{2}MT}
\mathcal{E}_L(t)+(\frac{1}{r_1-d}-\frac{1}{R+d})^{\frac{1}{2}}\cdot
\\
&\;\;\;\;\{\sup_{x_1\leq x\leq
1}|\sqrt{\rho_{in}}|e^{\frac{MT}{2}}+\frac{1}{(r_1-d)^2}\sup_{x_1\leq
x\leq 1}|\rho_{in}^{\frac{2-\gamma}{2}}|e^{\frac{2-\gamma}{2}MT}\}
(\mathcal{E}_L(t))^{\frac{1}{2}}
\end{split}
\]
where we have used the fact $$\int_{x_1}^1\frac{1}{\rho r^4}
dx=\int_{r_1(t)}^{R(t)}\frac{1}{r^2}dr\leq
\frac{1}{r_1-d}-\frac{1}{R+d}.$$ Note that $\frac{1}{r_1-d}\leq
\frac{1}{r_0-d}\leq 1$. Combining the above two inequalities, first
we get
\[
M\leq
C_{in}e^{\frac{MT}{2}}\{(\mathcal{E}(t))^{\frac{1}{2}}+\mathcal{E}(t)\}
\]
where $C_{in}=\sup\{\sup_{0\leq r\leq
r_1+d}|\frac{1}{\sqrt{\rho_{in}}}|,\sup_{x_1\leq x\leq
1}|\sqrt{\rho_{in}}|\}$. By Taylor expansion,
\[
M-C_{in}\{(\mathcal{E}(t))^{\frac{1}{2}}+\mathcal{E}(t)\}
\sum_{k=1}^{\infty}\frac{(MT)^k}{2^kk!}\leq
C_{in}\{(\mathcal{E}(t))^{\frac{1}{2}}+\mathcal{E}(t)\}
\]
Note that for sufficiently small $T$
$$C_{in}\{(\mathcal{E}(t))^{\frac{1}{2}}+\mathcal{E}(t)\}
\sum_{k=1}^{\infty}\frac{(MT)^k}{2^kk!}\leq \frac{M}{2}$$ and
therefore for such $T$, we get the following:
\begin{equation}
M\leq C_{in}\{(\mathcal{E}(t))^{\frac{1}{2}}+\mathcal{E}(t)\} \text{
for }0\leq t\leq T.\label{M}
\end{equation}
We have the exact same bound for $|\frac{D_t\rho}{\rho}|$. Next,
following the same path, we estimate $|\rho r^2D_xu|$.
\[
\sup_{0< x< 1}|\rho r^2D_xu|\leq \sup\{\sup_{0\leq r\leq
r_1+d}|\partial_r u|,\sup_{x_1\leq x\leq 1}|\rho r^2D_xu| \}
\]
\[
\begin{split}
\bullet\sup_{0\leq r\leq r_1+d}|\partial_r u|&\leq
\sum_{i=0}^2(\int_{B_{r_1+d}}|\nabla
\partial_{\mathbf{x}}^i\mathbf{u}|^2 d\mathbf{x})^{\frac{1}{2}}\\
&\leq \sup_{0\leq r\leq
r_1+d}|\frac{1}{\sqrt{\rho}}|\sum_{i=0}^2(\int_{B_{r_1+d}}\rho|\nabla
\partial_{\mathbf{x}}^i\mathbf{u}|^2 d\mathbf{x})^{\frac{1}{2}}\\
&\leq\sup_{0\leq r\leq
r_1+d}|\frac{1}{\sqrt{\rho_{in}}}|e^{\frac{MT}{2}}
(\mathcal{E}_E(t))^{\frac{1}{2}}
\end{split}
\]
\[
\begin{split}
\bullet\sup_{x_1\leq x\leq 1}|\rho r^2D_xu|&\leq \int_{x_1}^1|\rho
r^2
D_xu|dx+\int_{x_1}^1|D_x(\rho r^2D_xu)|dx\\
&\leq (\int_{x_1}^1\rho dx)^{\frac{1}{2}}(\int_{x_1}^1\rho
r^4|D_xu|^2 dx)^{\frac{1}{2}}+(\int_{x_1}^1\frac{1}{\rho
r^4}dx)^{\frac{1}{2}} (\int_{x_1}^1\rho|r^2D_x(\rho
r^2D_xu)|^2dx)^{\frac{1}{2}}\\
&\leq\sup_{x_1\leq x\leq
1}|\rho_{in}^{\frac{2-\gamma}{2}}|e^{\frac{2-\gamma}{2}MT}
\mathcal{E}_L(t)+(\frac{1}{r_1-d}-\frac{1}{R+d})^{\frac{1}{2}}\cdot
\\
&\;\;\;\;\{\sup_{x_1\leq x\leq
1}|\sqrt{\rho_{in}}|e^{\frac{MT}{2}}+\frac{1}{(r_1-d)^2}\sup_{x_1\leq
x\leq 1}|\rho_{in}^{\frac{2-\gamma}{2}}|e^{\frac{2-\gamma}{2}MT}
+\frac{1}{(r_1-d)^2}\} (\mathcal{E}_L(t))^{\frac{1}{2}}
\end{split}
\]
Because of (\ref{M}) and Taylor expansion, we can derive the
following: for sufficiently small $T$ and a constant $C_{in}$,
\[
\sup_{0< x< 1}|\rho r^2D_xu|\leq
C_{in}\{(\mathcal{E}(t))^{\frac{1}{2}}+\mathcal{E}(t)\} \text{ for
}0\leq t\leq T.
\]
Note that it also follows that $|\frac{u}{r}|$ is bounded in
$0<x<1$. The argument works for $|\rho r^2 D_tD_xu|$ and
$|\frac{D_tu}{r}|$ in the same way.\\

Next, we estimate $|\rho^{2\gamma-1}r^2D_x\rho|$ in $x_0\leq x\leq
1$. Because the cutoff function $\chi$ values $1$ only for  $x_1\leq
x\leq 1$, $|\rho^{2\gamma-1}r^2D_x\rho|$ for $x_0\leq x\leq x_1$
should be estimated in Eulerian coordinates. Note that $r_0-d\leq
r\leq r_1+d$ covers $x_0\leq x\leq 1$.
\[
\sup_{x_0\leq x\leq 1}|\rho^{2\gamma-1}r^2D_x\rho|\leq
\sup\{\sup_{r_0-d\leq r\leq
r_1+d}|\rho^{2\gamma-2}\partial_r\rho|,\sup_{x_1\leq x\leq
1}|\rho^{2\gamma-1}r^2D_x\rho| \}
\]
Apply the $L^1$ Sobolev embedding theorem in one dimension and in
turn the H$\ddot{o}$lder inequality.
\[
\begin{split}
\bullet\sup_{r_0-d\leq r\leq r_1+d}|\rho^{2\gamma-2}\partial_r\rho|
&\leq\int_{r_0-d}^{r_1+d}|\rho^{2\gamma-2}\partial_r\rho|dr
+\int_{r_0-d}^{r_1+d}|\partial_r(\rho^{2\gamma-2}\partial_r\rho)|dr\\
&\leq(\int_{r_0-d}^{r_1+d}\frac{\rho^{3\gamma-2}}{
r^2}dr)^{\frac{1}{2}}\sum_{i=1}^2(\int_{r_0-d}^{r_1+d}\rho^{\gamma-2}
|\partial_r^i\rho|^2
r^2dr)^{\frac{1}{2}}+\int_{r_0-d}^{r_1+d}\rho^{2\gamma-3}|\partial_r\rho|^2dr\\
&\leq \{\frac{1}{r_0-d}+1\}\sup_{r_0-d\leq r\leq
r_1+d}|\rho_{in}^{\gamma-1}|e^{(\gamma-1)MT}\mathcal{E}_E(t)
\end{split}
\]
\[
\begin{split}
\bullet \sup_{x_1\leq x\leq 1}|\rho^{2\gamma-1}r^2D_x\rho|&\leq
\int_{x_1}^1|\rho^{2\gamma-1}r^2D_x\rho|dx+\int_{x_1}^1
|D_x(\rho^{2\gamma-1}r^2D_x\rho)|dx\\
&\leq(\int_{x_1}^1\rho^{2\gamma}dx)^{\frac{1}{2}}
(\int_{x_1}^1\rho^{2\gamma-2}r^4|D_x\rho|^2dx)^{\frac{1}{2}}+
(\int_{x_1}^1\frac{1}{r^4}dx)^{\frac{1}{2}}
(\int_{x_1}^1\rho^{4\gamma-2}r^8|D_x^2\rho|^2dx)^{\frac{1}{2}}\\
&\;\;\;+(2\gamma-1)\int_{x_1}^1\rho^{2\gamma-2}r^2|D_x\rho|^2dx +
(\int_{x_1}^1\frac{\rho^{2\gamma-2}}{r^5}dx)^{\frac{1}{2}}
(\int_{x_1}^1\rho^{2\gamma-2}r^4|D_x\rho|^2dx)^{\frac{1}{2}}\\
&\leq \sup_{x_1\leq x\leq
1}|\rho_{in}^{\frac{\gamma+1}{2}}|e^{\frac{\gamma+1}{2}MT}
\mathcal{E}_L(t)+\frac{1}{(r_1-d)^2}\{(\mathcal{E}_L(t))^{\frac{1}{2}}
+\mathcal{E}_L(t)\}\\
&\;\;\;+\frac{1}{(r_1-d)^{\frac{5}{2}}}\sup_{x_1\leq x\leq
1}|\rho_{in}^{\frac{\gamma-1}{2}}|e^{\frac{\gamma-1}{2}MT}\mathcal{E}_L(t)
\end{split}
\]
Thus, as in the previous cases, we conclude that
\[
\sup_{x_0\leq x\leq 1}|\rho^{2\gamma-1}r^2D_x\rho|\leq
C_{in}\{(\mathcal{E}(t))^{\frac{1}{2}}+\mathcal{E}(t)\} \text{ for
}0\leq t\leq T
\]
for small enough $T$. And hence we have completed the verification
of all the Lagrangian terms in (\ref{K}). For Eulerian terms, we
give the detail for $\sup_{0\leq r\leq
r_2-d}|\frac{\partial_r\rho}{\rho}|$. Other terms such as
$\frac{\partial_t\rho}{\rho}$ and $\partial_tu$ can be estimated in
the same way by using the change of variables:
$\partial_t=D_t-r^2\rho uD_x$  in the overlapping region to estimate
them in Lagrangian interval $x_1\leq x\leq x_2$. First we observe
that it is enough to compute $\sup_{0\leq r\leq
r_2-d}|{\partial_r\rho}|$, since $$\sup_{0\leq r\leq
r_2-d}|\frac{1}{\rho}|\leq \sup_{0\leq r\leq
r_2-d}|\frac{1}{\rho_{in}}|e^{MT}.$$ Here is the estimate of
$\partial_r\rho$.
\[
\sup_{0\leq r\leq r_2+d}|{\partial_r\rho}|= \sup\{\sup_{0\leq r\leq
r_1+d}|{\partial_r\rho}|, \sup_{r_1+d\leq r\leq
r_2-d}|{\partial_r\rho}|\}
\]
\[
\begin{split}
\bullet\sup_{0\leq r\leq
r_1+d}|{\partial_r\rho}|&\leq\sum_{|\alpha|\leq 2}(\int_{B_{r_1+d}}
|\partial_{\mathbf{x}}^{\alpha}\partial_{\mathbf{x}}\rho|^2
d\mathbf{x})^{\frac{1}{2}}\\
&\leq\sup_{0\leq r\leq
r_1+d}|\rho^{\frac{2-\gamma}{2}}|\{\sum_{|\alpha|\leq
2}(\int_{B_{r_1+d}}
\rho^{\gamma-2}|\partial_{\mathbf{x}}^{\alpha}\partial_{\mathbf{x}}
\rho|^2d\mathbf{x})^{\frac{1}{2}}\\
&\leq \sup_{0\leq r\leq r_1+d}|\rho_{in}^{\frac{2-\gamma}{2}}|
e^{\frac{2-\gamma}{2}MT}(\mathcal{E}_E(t))^{\frac{1}{2}}
\end{split}
\]
Note that $D_x(\rho
r^2D_x\rho)=\frac{2D_x\rho}{r}+r^2|D_x\rho|^2+\rho r^2D_x^2\rho$.
\[
\begin{split}
 \bullet\sup_{r_1+d\leq r\leq
r_2-d}|{\partial_r\rho}|&\leq\sup_{x_1\leq x\leq x_2} |\rho r^2
D_x\rho|\\
&\leq
\int_{x_1}^{x_2}|\rho r^2D_x\rho|dx+\int_{x_1}^{x_2}|D_x(\rho r^2D_x\rho)|dx\\
&\leq
(\int_{x_1}^{x_2}\{\rho^{4-2\gamma}+\frac{1}{\rho^{2\gamma-2}r^6}\}dx)
^{\frac{1}{2}}(\int_{x_1}^{x_2}\rho^{2\gamma-2}r^4
|D_x\rho|^2dx)^{\frac{1}{2}}\\
&\;\;\;+\sup_{x_1\leq x\leq
x_2}|\frac{1}{\rho^{2\gamma-2}}|\int_{x_1}^{x_2}\rho^{2\gamma-2}r^4
|D_x\rho|^2dx \\
&\;\;\;+(\int_{x_1}^{x_2}\frac{1}{\rho^{4\gamma-4}r^4}dx)
^{\frac{1}{2}}(\int_{x_1}^{x_2}\rho^{4\gamma-2}r^8
|D_x^2\rho|^2dx)^{\frac{1}{2}}
\end{split}
\]
This concludes the proof of the lemma.
\end{proof}\

Thus the a priori estimates can be closed at this point. In the
following sections, being inspired by the a priori estimates, we
construct local in time strong solutions to the
Navier-Stokes-Poisson system. \\

\section{Approximate Scheme} \label{section6}

Strong solutions to the free boundary problem to (\ref{nspE}) with
given initial data $\rho(0,r)=\rho_{in}(r)$ and $u(0,r)=u_{in}(r)$,
and boundary conditions (\ref{b1}), (\ref{b2}), (\ref{b3}) are
constructed by an approximate scheme, where an approximate velocity
is obtained by solving a parabolic linear PDE in Lagrangian
coordinates and the following approximate density profile is defined
by the flow generated by the approximate velocity. Due to the
coordinate singularity at the origin in Lagrangian formulation, the
corresponding Eulerian formulation is invoked and both Lagrangian
and Eulerian estimates are obtained.

Let  initial data $\rho_{in},\;u_{in}$ be given: for $0\leq r\leq
R$, $\rho_{in}(r),\;u_{in}(r)$ satisfy
\begin{equation}
\mathcal{E}|_{(\rho_{in},u_{in})}\leq A\text{ for some }A>0;\;
\rho_{in}(R)=0;\;\rho_{in}(r)>0\text{ for }0\leq
r<R;\;\int_0^R\rho_{in}r^2dr=1.\label{id}
\end{equation}
Starting with these initial data, in the following, we describe an
approximate scheme that leads to approximations $(\rho^n,u^n,r^n)$
for all $n\geq 0$. Subsequently, we study the existence, uniqueness,
and regularity of them. Approximations  are shown to be uniformly
bounded under the energy space characterized by $\mathcal{E}$, and
finally, a strong solution is obtained by taking the limit.

For given initial data $\rho_{in}\text{ and } u_{in}$ satisfying
(\ref{id}) in Eulerian coordinates, first introduce a Lagrangian
variable $x$ as follows:
\[
x\equiv \int_0^r \rho_{in}s^2ds,\;\; 0\leq x\leq 1.
\]
Then $\rho_{in}\text{ and } u_{in}$, denoted by $\rho^0\text{ and
}u^0$ respectively, can be regarded as functions of $x$. Define
$r^0$ by
\[
r^0\equiv \{3\int_0^{x}\frac{1}{\rho^0} dy\}^{\frac{1}{3}}.
\]
We would like to define the sequence $\{\rho^{n},r^n,u^n\}$
inductively for all $n\geq 0$. Suppose that $\rho^n$, $r^n$, and
$u^n$ are known functions. Firstly, consider the following linear
partial differential equations for $u^{n+1}$:
\begin{equation}
\begin{split}
D_tu^{n+1}-\mu D_{x} (\rho^{n}(r^{n})^4D_{x}u^{n+1})
+\mu\frac{2u^{n+1}}{\rho^{n}(r^{n})^2}=
-(r^{n})^2D_{x}p^{n}-\frac{4\pi x}{(r^{n})^2} \label{al}
\end{split}
\end{equation}
with the initial data $u^{n+1}(0,x)=u_{in}$ and boundary conditions
$$u^{n+1}(t,0)=0\text{ and
}(\mu\rho^{n}(r^{n})^2D_{x}u^{n+1}-p^{n})(t,1)=0.$$ (\ref{al}) is a
parabolic-type equation with degenerate coefficients at $x=0,1$, but
the singularity is either coordinate singularity or point
singularity and hence the existence, uniqueness, and regularity
follow from the classical theory. Note that  by the change of
variables
\begin{equation}
D_x=\frac{1}{\rho^n(r^n)^2}\partial_{r^n}\;\text{ and
}\;D_t=\partial_t+ (D_tr^n)\partial_{r^n} ,\label{cv}
\end{equation}
(\ref{al}) can be written in Eulerian coordinates $(t,r^n)$ as
follows:
\begin{equation}
\partial_tu^{n+1}+(D_tr^n)\partial_{r^n}u^{n+1}-\mu\frac{1}{\rho^n(r^n)^2}
\partial_{r^n}((r^n)^2\partial_{r^n}u^{n+1})+\mu\frac{2u^{n+1}}{\rho^n (r^n)^2}
=-\frac{1}{\rho^n}\partial_{r^n}p^n-\frac{4\pi x}{(r^n)^2}\label{ae}
\end{equation}
where
$D_tr^n=-\frac{1}{(r^n)^2}\int_0^x\frac{D_t\rho^n}{(\rho^n)^2}dy$
and $\partial_tx=-\rho^n(r^n)^2D_tr^n;$
$\partial_{r^n}x=\rho^n(r^n)^2.$
 Next, imitating the formula (\ref{rho1}), define $\rho^{n+1}$ by
\begin{equation}
\rho^{n+1}(t,x)\equiv \rho^0\exp\{-\int_0^t
\rho^n(r^n)^2D_xu^{n+1}+\frac{2u^{n+1}}{r^n} d\tau\}\label{rhon+1}.
\end{equation}
It is straightforward to check that $\rho^{n+1}(t,x)$ satisfies the
following equation:
\begin{equation}
D_t\rho^{n+1}+\{\rho^n(r^n)^2D_{x}u^{n+1}+\frac{2u^{n+1}}{r^n}
\}\rho^{n+1}=0.\label{rho}
\end{equation}
It reads in Eulerian coordinates as follows:
\begin{equation}
\partial_t\rho^{n+1}+(D_tr^n)
\partial_{r^n}\rho^{n+1}+\rho^{n+1}\{\partial_{r^n}u^{n+1}+\frac{2u^{n+1}}{r^n}\}=0
\label{rhoe}
\end{equation}
Lastly, we define $r^{n+1}$ by
\begin{equation}
r^{n+1}\equiv \{3\int_0^{x}\frac{1}{\rho^{n+1}}
dy\}^{\frac{1}{3}}.\label{ar}
\end{equation}

We need to make sure (\ref{al}) is solvable and (\ref{rhon+1}) and
(\ref{ar}) make sense in an appropriate sense. First, we study
(\ref{al}) in a weak formulation in Lagrangian coordinates, and
establish the regularity of weak solution. Interior regularity is
standard because (\ref{al}) is parabolic bounded away from the
boundary, while boundary regularity is obtained with weights in the
form of integrals. Once one has regularity of $u^{n+1}$ with respect
to $\rho^n,r^n$ coefficients, one can check $\rho^{n+1}, r^{n+1}$
are well-defined.  Eulerian regularity easily follows since
(\ref{al}) and (\ref{ae}) are equivalent in the interiors. We remark
 that for Eulerian regularity, $D_tr^n\partial_{r^n}\rho^{n+1}$
related terms in (\ref{ae}) can be taken care of by integrating by
parts, since $D_tr^n$ is more regular. The most important task is to
derive the uniform bound on $n$ so that we may
conclude that  the limit functions are a desired solution.\\

The rest of this section is devoted to studying weak solutions of
above approximate equations in the line of existence, uniqueness,
and regularity. We use Galerkin's method well-illustrated in
\cite{e}. The difference is that the space we will work is not a
typical Sobolev space, but inherited from the special structure of
target equations. The first purpose is to show the existence of weak
solution $u^{n+1}$ to (\ref{al})  for given $\rho^{n}$, $r^{n}$ and
$u^n$. Without confusion, we will drop the index $n$ from now on. We
assume that we have as much regularity of $\rho$ and $r$ as needed.
In particular, we keep in mind the behavior of stationary solutions:
$$\rho_0\sim (1-x)^{\frac{1}{\gamma}} \text{ if } x\sim 1 \text{ and }
r\sim x^{\frac{1}{3}}\text{ if }r\sim 0.$$ We start with Largrangian
equation (\ref{al}):
\[
D_tu-\mu D_{x} (\rho r^4D_{x}u) +\mu\frac{2u}{\rho r^2}=
-r^2D_{x}p-\frac{4\pi x}{r^2}
\]
Firstly, we define the notion of weak solution. To do so, a Hilbert
space $H$ is introduced:
\[
H=Cl\{ u\in C^{\infty}(0,1): \int_0^1\rho r^4 |D_x u|^2 +\frac{2
u^2}{\rho r^2} dx < \infty,\;u(0)=0\}
\]
It is straightforward to check $H\subset L^2(0,1)$.\\

\begin{definition}
We say $u\in L^2(0,T;H)$ with $u'\in L^2(0,T;H^{\ast})$ is a
weak solution of (\ref{al}) provided
\[
\int_0^1 u' v  dx+\mu\int_0^1\rho r^4 D_xu D_x vdx+\mu
\int_0^1\frac{2u v}{\rho r^2}dx
=\int_0^1r^2pD_xvdx+\int_0^1(\frac{2p}{\rho r}-\frac{4\pi x}{r^2}) v
dx
\]
for each $v\in H$ and a.e. time $0\leq t\leq T$, and $u(0)=u_{in}$.
$H^{\ast}$ is the dual space of $H$ and $'=D_t$.
\end{definition}\

\begin{lemma} Assume $u_{in}\in L^2$, $\rho^{-\frac{1}{2}}p\in L^2(0,T;L^2)$,
and $\rho^{\frac{1}{2}}r^{-1}x\in L^2(0,T;L^2)$.
 There exist a unique weak solution $u\in L^2(0,T;H)$ with
$u'\in L^2(0,T;H^{\ast})$ to (\ref{al}). Furthermore, there exists a
constant  $C_{\mu}>0$ such that
\begin{equation}
\max _{0\leq t\leq T}||u(t)||_{L^2}+||u||_{L^2(0,T;H)}+
||u'||_{L^2(0,T;H^{\ast})}\leq
C_{\mu}\{||u_{in}||_{L^2}+||\frac{p}{\sqrt{\rho}}||_{L^2(0,T;L^2)}+
||\sqrt{\rho}\frac{4\pi x}{r}||_{L^2(0,T;L^2)}\}.
\end{equation}
\end{lemma}

\begin{proof} Let $w_k=w_k(x)\;(k=1,2,...)$ be an orthogonal basis of $H$ and
orthonormal in $L^2$ when $t=0$, i.e., $\rho(0)=\rho_{in}$ and
$r(0)=r_{in}$. Then $\{w_k\}$ forms a basis of $H$ for $0\leq t\leq
T$, where $T$ is sufficiently small, due to smoothness of $\rho,r$.
Fix a positive integer $m$. We seek a function $u_m:[0,T]\rightarrow
H$ of the form
\begin{equation}
u_{m}(t)=\sum_{k=1}^{m}d_m^k(t)w_k,\label{form}
\end{equation}
where
\begin{equation}
d_m^k(0)=\int_0^1 u_{in}w_kdx\;\; (k=1,...,m)\label{in}
\end{equation}
and for each $k=1,...,m,\;0\leq t\leq T,$
\begin{equation}
\int_0^1 u_m' w_k dx+\mu\int_0^1\rho r^4 D_xu_m D_x w_kdx+\mu
\int_0^1\frac{2u_m w_k}{\rho
r^2}dx=\int_0^1r^2pD_xw_kdx+\int_0^1(\frac{2p}{\rho r}-\frac{4\pi
x}{r^2}) w_k dx.\label{appw}
\end{equation}\

\noindent Claim 1. For each $m$, there exists a unique function
$u_m$ of the
form (\ref{form}) satisfying (\ref{in}) and (\ref{appw}).\\

\noindent Proof of Claim 1: Note that
\[
\int_0^1 u'_m(t)w_k dx= {d_m^{k}}'(t),\;\;\mu\int_0^1\rho r^4 D_xu_m
D_x w_kdx+\mu \int_0^1\frac{2u_m w_k}{\rho
r^2}dx=\sum_{l=1}^{m}e^{kl}(t)d^l_m(t)
\]
where $e^{kl}(t)=\mu\int_0^1\rho r^4 D_xw_l D_x w_kdx+\mu
\int_0^1\frac{2w_l w_k}{\rho r^2}dx$. Let us write
$$f^k(t)=\int_0^1r^2pD_xw_kdx+\int_0^1(\frac{2p}{\rho r}-\frac{4\pi
x}{r^2}) w_k dx.$$ Then (\ref{appw}) becomes the linear system of
ODEs
\begin{equation}
{d^k_m}'(t)+\sum_{l=1}^{m}e^{kl}(t)d^l_m(t)=f^k(t),\label{ode1}
\end{equation}
subject to the initial condition (\ref{in}). According to the
standard existence theory for ordinary differential equations, there
exists a unique absolutely continuous functions $d_m^k(t)$
satisfying (\ref{in}) and (\ref{ode1}) for a.e. $0\leq t\leq T$.
And then $u_m$ defined by (\ref{form}) solves (\ref{appw}).\\

\noindent Claim 2. The following energy estimates hold:
\begin{equation}
\begin{split}
\max _{0\leq t\leq T}||u_m(t)||_{L^2}+||u_{m}||_{L^2(0,T;H)}+
||u_m'||_{L^2(0,T;H^{\ast})}&\\
\leq
C_{\mu}\{||u_{in}||_{L^2}+||\frac{p}{\sqrt{\rho}}&||_{L^2(0,T;L^2)}+
||\sqrt{\rho}\frac{4\pi x}{r}||_{L^2(0,T;L^2)}\},
\end{split}
\end{equation}
where $C_{\mu}$ is a constant independent of $m$.\\

\noindent Proof of Claim 2: Multiplying (\ref{appw}) by $d_m^k$ and
summing over $k$, we get
\[
\begin{split}
&\int_0^1 u_m' u_m dx+\mu\int_0^1\rho r^4|D_xu_m|^2+
\frac{2u_m^2}{\rho r^2}
dx=\int_0^1r^2pD_xu_mdx+\int_0^1(\frac{2p}{\rho r}-\frac{4\pi
x}{r^2}) u_m dx\\
 \Rightarrow&\; \frac{1}{2}\frac{d}{dt}||u_m||_{L^2}^2
+\mu\int_0^1\rho r^4|D_xu_m|^2+ \frac{2u_m^2}{\rho r^2} dx \leq
\frac{\mu}{2}\int_0^1\rho r^4|D_xu_m|^2+ \frac{2u_m^2}{\rho r^2}
dx\\
&\;\;\;\;\;\;\;\;\;\;\;\;\;\;\;\;\;\;\;\;\;\;\;\;\;\;
\;\;\;\;\;\;\;\;\;\;\;\;\;\;\;\;\;\;\;\;\;\;\;\;\;\;\;\;\;\;\;\;\;\;\;\;
\;\;\;\;\;\;\;\;+\frac{2}{\mu}\int_0^1A^2\rho^{2\gamma-1}+\rho\frac{16\pi^2
x^2}{r^2} dx
\end{split}
\]
By taking $t$-integral,
\[
\begin{split}
&||u_m||_{L^2}^2+\mu\int_0^t \int_0^1\rho r^4|D_xu_m|^2+
\frac{2u_m^2}{\rho r^2} dxds\leq
||u_{m}(0)||_{L^2}^2+\frac{4}{\mu}\int_0^t\int_0^1A^2\rho^{2\gamma-1}
+\rho\frac{16\pi^2 x^2}{r^2} dxds\\
\Rightarrow\;&\max _{0\leq t\leq T}||u_m(t)||^2_{L^2}+\mu\int_0^T
\int_0^1\rho r^4|D_xu_m|^2+ \frac{2u_m^2}{\rho r^2} dxds\\
&\;\;\;\;\;\;\;\;\;\;\;\;\;\;\;\;\;\;\;\;\;\;\;\;\;\;\;\;\;\;\;\;\;\;
\;\;\;\;\;\;\;\;\;\;\;\;\;\;\;\;\;\;\;\;\;\;\;\;\;\;\;\;\;\;\;\;\leq
||u_{in}||_{L^2}^2+\frac{4}{\mu}(||\frac{p}{\sqrt{\rho}}||^2_{L^2(0,T;L^2)}+
||\sqrt{\rho}\frac{4\pi x}{r}||^2_{L^2(0,T;L^2)})
\end{split}
\]
Fix $v\in H$ with $||v||_H\leq 1$. Write $v=v^1+v^2$, $v^1\in span
\{w_k\}_{k=1}^m$, and $\int_0^1 v^2 w_k dx=0,$ $k=1,...,m.$ Since
$\{w_k\}$ is orthogonal, $||v^1||_H\leq||v||_H\leq 1$. From
(\ref{appw}),
\[
\begin{split}
\int_0^1 u_m' v dx =\int_0^1u_m' v^1 dx&=\int_0^1r^2pD_x
v^1dx+\int_0^1(\frac{2p}{\rho r}-\frac{4\pi
x}{r^2}) vdx-\mu\int_0^1\rho r^4 D_x u_m D_xv^1 +\frac{2u_m v^1}{\rho r^2}dx\\
&\leq C\{||\frac{p}{\sqrt{\rho}}||_{L^2}+ ||\sqrt{\rho}\frac{4\pi
x}{r}||_{L^2}+\mu||u_m||_H\}
\end{split}
\]
Hence, $$||u_m'||_{H^{\ast}}\leq C\{||\frac{p}{\sqrt{\rho}}||_{L^2}+
||\sqrt{\rho}\frac{4\pi x}{r}||_{L^2}+\mu||u_m||_H\}\;\text{ and }$$
$$\int_0^T||u_m'||_{H^{\ast}}dt\leq C_{\mu}\{||u_{in}||_{L^2}
 +\\ ||\frac{p}{\sqrt{\rho}}||_{L^2(0,T;L^2)}+
||\sqrt{\rho}\frac{4\pi x}{r}||_{L^2(0,T;L^2)}\}.$$

Now we pass to limits as $m\longrightarrow\infty$. According to
energy estimates, $\{u_m\}_{m=1}^{\infty}$ is bounded in
$L^2(0,T;H)$, $\{u_m'\}$ bounded in $L^2(0,T;H^*)$ and therefore,
there exist a subsequence $\{u_{m_l}\}\subset \{u_m\}$ and functions
$u\in L^2(0,T;H),\;u'\in L^2(0,T; H^*)$ such that
\[
u_{m_l}\rightharpoonup u\;\;\text{weakly in }L^2(0,T;H),\;\;\;
u_{m_l}'\rightharpoonup u'\;\;\text{weakly in }L^2(0,T;H^*).
\]
Fix $N$, and consider $v\in C^1(0,T;H)$ having the form of
$v(t)=\sum_{k=1}^N d^k(t)w_k$. Choose $m\geq N$,
\[
\int_0^T\int_0^1 u_m' v dxdt+\mu\int_0^T\int_0^1\rho r^4 D_xu_m D_x
v+\frac{2 u_m v}{\rho r^2}
dxdt=\int_0^T\int_0^1r^2pD_xvdx+\int_0^1(\frac{2p}{\rho
r}-\frac{4\pi x}{r^2}) v dxdt.
\]
Setting $m=m_l$, pass to weak limits to get for all $ v\in
L^2(0,T;H),$
\[
\int_0^T\int_0^1 u' v dxdt+\mu\int_0^T\int_0^1\rho r^4 D_xu
D_xv+\frac{2 uv}{\rho r^2}
dxdt=\int_0^T\int_0^1r^2pD_xvdx+\int_0^1(\frac{2p}{\rho
r}-\frac{4\pi x}{r^2}) v dxdt.
\]
In particular,
\[
\int_0^1 u' v dx+\mu\int_0^1\rho r^4 D_xu D_xv+ \frac{2 uv}{\rho
r^2}dx=\int_0^1r^2pD_xvdx+\int_0^1(\frac{2p}{\rho r}-\frac{4\pi
x}{r^2}) v dx,\;\forall v\in H \;\text{a.e. } t.
\]
This assures the existence of weak solution to (\ref{al}).
Uniqueness of weak solutions easily follows from energy estimates:
let $u_1,\;u_2$ be two weak solutions with the same initial data,
then $u\equiv u_1-u_2$ satisfies the following:
\[
\int_0^T\int_0^1 u' v dxdt+\mu\int_0^T\int_0^1\rho r^4 D_xu
D_xv+\frac{2 uv}{\rho r^2} dxdt=0,\;\forall v\in L^2(0,T;H).
\]
Choose $v=u\in L^2(0,T;H)$, and we get
\[
\frac{1}{2}||u||_{L^2}+\mu\int_0^T\int_0^1\rho r^4 |D_x u|^2
+\frac{2u^2}{\rho r^2}dxdt=\frac{1}{2}||u(0)||_{L^2}=0,
\]
and hence $u=0$ a.e.
\end{proof}\

Next, we try to get the regularity of the weak solution $u$ obtained
in the above. First, we establish Lagrangian regularity. Eulerian
regularity is obtained by the change of variable (\ref{cv}) in the
integral form. One can study the weak solution by the same Galerkin
method in Eulerian formulation, but we skip the details here.  The
next lemma regards time regularity.\\

\begin{lemma} Assume $\sup_{0<x<1}|\frac{\rho'}{\rho}|\leq C_1$ and
$\sup_{0<x<1}|\frac{r'}{r}|\leq C_2$ for $0\leq t\leq T$. In
addition, assume $u_{in}\in H$ and $r^2D_xp+\frac{4\pi x}{r^2}\in
L^2(0,T;L^2)$. Then $u\in L^{\infty}(0,T;H),\;u'\in L^2(0,T;L^2)$
with the estimate
\begin{equation}
\begin{split}
\int_0^T||u'||^2_{L^2}dt+\sup_{0\leq t\leq T} &\mu||u(t)||_{H}^2\leq
C\{||u_{in}||_H^2+||u_{in}||_{L^2}^2 \\&+||r^2D_xp+\frac{4\pi
x}{r^2}||^2_{L^2(0,T;L^2)}+||\frac{p}{\sqrt{\rho}}||^2_{L^2(0,T;L^2)}+
||\sqrt{\rho}\frac{4\pi x}{r}||^2_{L^2(0,T;L^2)}\}.\label{tr}
\end{split}
\end{equation}
\end{lemma}

\begin{proof} Again, Galerkin method is used. We start with (\ref{al}).
Let $u_m(t)=\sum_{k=1}^m d_m^k(t)w_k$. Multiplying (\ref{appw}) by
${d_m^k}'$ and summing over $k$ to get
\[
\begin{split}
\int_0^1 u_m' u_m' dx+\mu\int_0^1\rho r^4 D_x u_m D_x u_m'+ \frac{2
u_m u_m'}{\rho r^2}& dx
=\int_0^1r^2pD_xu_m'dx+\int_0^1(\frac{2p}{\rho r}-\frac{4\pi
x}{r^2}) u_m'dx\\
&\;\;\;\;=-\int_0^1(r^2D_xp+\frac{4\pi x}{r^2})u_m'dx\\
\Rightarrow\;||u_m'||_{L^2}^2+\frac{d}{dt}\frac{\mu}{2} \int_0^1\rho
r^4 |D_x u_m|^2+\frac{2 u_m^2}{\rho r^2}dx& \leq \frac{1}{2}
||u_m'||^2_{L^2}+\frac{1}{2}
||r^2D_xp+\frac{4\pi x}{r^2}||^2_{L^2}\\
+&\underline{\frac{\mu}{2}\int_0^1(\rho'r^4+4\rho r^3r')(D_x u_m)^2
-(\frac{2\rho'}{\rho^2 r^2}+\frac{4r'}{\rho r^3})u_m^2dx}_{\star\star}\\
\end{split}
\]
Since $|\frac{\rho'}{\rho}|\leq C_1$ and $|\frac{r'}{r}|\leq C_2$,
$\star\star\leq 2\mu(C_1+C_2)||u_m||_{H}^2$, we get
\[
\begin{split}
&\int_0^T||u_m'||^2_{L^2}dt+\sup_{0\leq t\leq T}
\mu\int_0^1\rho r^4 |D_x u_m|^2+\frac{2 u_m^2}{\rho r^2}dx\\
&\leq \mu\int_0^1\rho_{in} r_{in}^4 |D_x {u_m}_{in}|^2 +\frac{2
{u_m}_{in}^2}{\rho r_{in}^2} dx+\int_0^T||r^2D_xp+\frac{4\pi
x}{r^2}||_{L^2}^2 dt
+4\mu(C_1+C_2)\int_0^T||u_m||_{H}^2dt\\
&\leq C\{||u_{in}||_H^2+||u_{in}||_{L^2}^2+||r^2D_xp+\frac{4\pi
x}{r^2}||^2_{L^2(0,T;L^2)}+||\frac{p}{\sqrt{\rho}}||^2_{L^2(0,T;L^2)}+
||\sqrt{\rho}\frac{4\pi x}{r}||^2_{L^2(0,T;L^2)}\}
\end{split}
\]
Pass to limit $m\rightarrow\infty$. (\ref{tr}) holds and the lemma
follows.
\end{proof}\

Now we would like to establish regularity in $x$ variable. Note that
(\ref{al}) is one-dimensional parabolic equation as long as $x$ is
bounded away from the boundary and hence interior regularity can be
easily shown by using standard different quotients method (see
Section 7.1 in \cite{e}), i.e.  $u\in H_{\text{loc}}^2(0,1)$. Here
$H^2$ represents the usual Sobolev space. Recall that
\begin{equation}
\mu\int_0^1\rho r^4 D_x u D_x v dx+\mu\int_0^1 \frac{2 uv}{\rho
r^2}dx =\int_0^1r^2pD_xvdx+\int_0^1(\frac{2p}{\rho r}-\frac{4\pi
x}{r^2}-u') v dx, \;\forall v\in H.\label{weak}
\end{equation}
We can now integrate by parts in (\ref{weak}) by approximating $H$
with $v\in C^{\infty}_c(0,1)\subset H$:
\[
-\mu\int_0^1D_x(\rho r^4 D_x u)  v dx+\mu\int_0^1 \frac{2 uv}{\rho
r^2}dx =-\int_0^1 (r^2D_xp+\frac{4\pi x}{r^2}+u')vdx, \;\forall v\in
C^{\infty}_c(0,1),
\]
Therefore, $u$ actually solves the PDE for a.e., and the following
estimate 
can be obtained from the
equation:
\begin{equation}
\begin{split}
\frac{\mu^2}{2}\int_0^1 \rho r^2|r^2 D_x(\rho r^2 D_xu)|^2dx\leq
\int_0^1 \rho r^2|u'|^2dx+\mu^2\int_0^1\frac{4u^2}{\rho r^2}
+4\rho r^4|D_xu|^2dx\\
+\int_0^1\rho r^2|r^2D_xp+\frac{4\pi x}{r^2}|^2 dx\label{xr}
\end{split}
\end{equation}
Note that $L^1$ estimate for $D_x(\rho r^2D_xu)$ can be obtained
from the equation. Next we claim that $u$ satisfies the boundary
condition
 $(\rho r^2 D_xu)(t,1)=0$ for $0\leq t\leq T$ in the trace sense.
Note that $D_x(\rho r^2 D_x u)$ is in $L^1(r^2dx)$ and therefore by
the Trace theorem, $\rho r^2 D_x u$ at $x=1$ is well-defined. Thus,
in (\ref{weak}), we can integrate by parts up to boundary for the
first term to get: recall $v(0)=0$,
\[
\mu(\rho r^4 D_x u  v)(t,1) -\mu\int_0^1D_x(\rho r^4 D_x u)v
dx+\mu\int_0^1 \frac{2 uv}{\rho r^2}dx =-\int_0^1
(r^2D_xp+\frac{4\pi x}{r^2}+u')vdx,\;\forall v\in H
\]
We already know that $u$ solves the PDE a.e. and so, all other terms
vanish, and therefore we obtain
 $(\rho r^2 D_xu)(t,1)=0$, the desired  boundary condition. We have
 proven spatial regularity of $u^{n+1}$:\\

\begin{lemma}
The weak solution $u$ solves (\ref{al}) and it is regular: $u\in
H^2_{\text{loc}}(0,1)$ and weak boundary regularity is given by the
estimate (\ref{xr}). Moreover, $u$ satisfies the boundary condition.
\end{lemma}\

\noindent{\bf Higher regularity:} Now let us build the higher
regularity based on the previous regularity assertion. Set
$\widetilde{u}=u'$. Differentiating (\ref{al}) with respect to $t$,
one can check that $\widetilde{u}$ is the unique weak solution of
\[
D_t\widetilde{u}-\mu D_x(\rho r^4 D_x\widetilde{u})
+\mu\frac{2\widetilde{u}}{\rho r^2}=-D_t(r^2D_xp+\frac{4\pi
x}{r^2})+\mu D_x(D_t(\rho r^4) D_x u)+\mu D_t(\frac{2}{\rho r^2})u
\]
with given compatible, well-prepared initial data. Note that each
term in the RHS is either known function or lower order derivative
terms. In the weak formulation, the spatial derivative of $p$ and
$D_t(\rho r^4)D_xu$ should be moved to test functions.\\

\begin{lemma}
The weak solution $u$ attains higher regularity as long as initial
data $u_{in}$ as well as coefficients $\rho^n,r^n$ are regular:
$u\in H^4_{\text{loc}}(0,1)$ and weak boundary regularity is
available in the integral form.
\end{lemma}

\begin{proof}
By the same argument as before, we get the regularity of
$\widetilde{u}$ with the following energy estimates:
\[
\begin{split}
\sup_{0\leq t\leq T}||\widetilde{u}||_{L^2}+||\widetilde{u}||
_{L^2(0,T;H)}+||\widetilde{u}'|| _{L^2(0,T;H*)}\leq
C(||\widetilde{f}_1||_{L^2}
+||\widetilde{u}_{in}||_{L^2}),\\
||\widetilde{u}'||_{L^2(0,T;L^2)}+\sup_{0\leq t\leq T}\sqrt{\mu}
||\widetilde{u}||_{H}\leq C(||\widetilde{u}_{in}||_H
+||\widetilde{u}_{in}||_{L^2}+||\widetilde{f}_2||_{L^2(0,T;L^2)}),
\end{split}
\]
where $\widetilde{f}_1$ is a function of $\rho,r,u, D_xu D_t\rho,D_t
r$ and $\widetilde{f}_2$ is a function of $\widetilde{f}_1$ and
$D_x\rho$. It is routine to check $\widetilde{u}\in
H^2_{\text{loc}}(0,1)$ and hence, $\widetilde{u}$ solves the above
PDE for a.e., and the following estimate (weak boundary regularity)
is obtained:
\[
\begin{split}
\frac{\mu^2}{2}\int_0^1 \rho r^2(r^2D_x(\rho
r^2D_x\widetilde{u}))^2dx \leq\int_0^1 \rho r^2(\widetilde{u}')^2dx
+\mu^2\int_0^1\frac{4\widetilde{u}^2}{\rho r^2}+4\rho
r^4|D_x\widetilde{u}|^2dx+\int_0^1\rho r^2(\widetilde{f}_2)^2 dx\\
\end{split}
\]
In the same line, letting $\widetilde{\widetilde{u}}
=\widetilde{u}'$, we get the desired regularity together with the
estimates:
\[
\begin{split}
\sup_{0\leq t\leq T}||\widetilde{\widetilde{u}}||_{L^2}
+||\widetilde{\widetilde{u}}||
_{L^2(0,T;H)}+||\widetilde{\widetilde{u}}'|| _{L^2(0,T;H*)}\leq
C(||\widetilde{\widetilde{f}}_1||_{L^2}
+||\widetilde{\widetilde{u}}_{in}||_{L^2})\\
||\widetilde{\widetilde{u}}'||_{L^2(0,T;L^2)}+ \sup_{0\leq t\leq
T}\sqrt{\mu} ||\widetilde{\widetilde{u}}||_{H}\leq
C(||\widetilde{\widetilde{u}}_{in}||_H
+||\widetilde{\widetilde{u}}_{in}||_{L^2}
+||\widetilde{\widetilde{f}}_2||_{L^2(0,T;L^2)})\\
\mu^2\int_0^1 \rho r^2(r^2D_x(\rho r^2D_x
\widetilde{\widetilde{u}}))^2 dx \leq \int_0^1 \rho
r^2(\widetilde{\widetilde{u}}')^2dx
+\mu^2\int_0^1\frac{4\widetilde{\widetilde{u}}^2}{\rho r^2}+4\rho
r^4|D_x\widetilde{\widetilde{u}}|^2dx +\int_0^1\rho
r^2(\widetilde{\widetilde{f}}_2)^2 dx
\end{split}
\]
Now it remains to investigate the spatial regularity regarding
$D_{x}^3u$ and $D_{x}^4u$. Again, it is straightforward to see
 $u\in H^4_{\text{loc}}$. To derive boundary estimates, let us go
back to (\ref{al}). Since $\rho$ vanishes at the boundary with
certain rate, we need some weight depending on $\rho$ in order to
control the spatial derivatives of $\rho$ and lower order derivative
terms. Set the RHS of (\ref{al}). Firstly, multiply by $\rho r$,
differentiate with respect to $x$, we get
\begin{equation}
\mu D_x(\rho r^3D_x(\rho r^2D_xu))=D_x(\rho rD_tu) -D_x(\rho
rf)-2\mu\underline{ D_x(\rho r^3D_x(\frac{u}{r}))}_{\ast}
\label{uxxx}
\end{equation}
\[
\ast=D_x(\rho r^2D_xu)-D_x(\frac{u}{r})=D_x(\rho r^2D_xu)
-\frac{D_xu}{r}+\frac{u}{\rho r^4}
\]
In the view of the previous weak boundary regularity, we obtain the
following estimate:
\[
\begin{split}
\frac{\mu^2}{2}\int_0^1\rho r^6|D_x(\rho r^3D_x(\rho
r^2D_xu))|^2dx\leq& \int_0^1\rho r^6(D_x(\rho
rD_tu))^2dx+\int_0^1\rho r^6 (D_x(\rho
rf))^2dx\\&+4\mu^2\int_0^1\rho r^6 (D_x(\rho
r^3D_x(\frac{u}{r})))^2dx
\end{split}
\]
Note that two other terms in the RHS are also bounded. Multiply
(\ref{uxxx}) by $\rho r^3$, differentiate in $x$, and square each
term to get
\[
\frac{\mu^2}{2}\int_0^1\rho r^6|D_x( \rho r^3D_x (\rho r^3D_x(\rho
r^2D_xu)))|^2dx\leq C
\]
This completes the Lagrangian regularity of $u$.
\end{proof}\

We remark that the boundary regularity is weak in a sense that
$r,\rho$ as weight functions vanish at $x=0, 1$ respectively.
Due to the interior Lagrangian regularity of $u$ and equivalence of
(\ref{al}) and (\ref{ae}) away from the boundary, $u$ also solves
the PDE (\ref{ae}) for a.e. 
The vacuum boundary in Eulerian coordinates $R^n(t)$ is defined by
$(3\int_0^1\frac{1}{\rho^n}dy)^{\frac{1}{3}}$ and $x=0$ corresponds
to $r^n=0$. Corresponding Eulerian regularity can be obtained by the
change of variable (\ref{cv}). We omit the details here. In the next
section, we provide the detailed energy estimates in both Lagrangian
and Eulerian coordinates with cutoff functions.

Thus $u^{n+1}$ and therefore $\rho^{n+1}$, $r^{n+1}$ are all
well-defined. They are regular in the classical sense away from the
boundary. Motivated by the a priori estimates, in the following
section, we study the behavior $u^n,\rho^n,  r^n$ under the energy
$\mathcal{E}$. It provides not only the regularity of
$\rho,u$ up to the boundary but also the unform bounds on $n$.\\

\section{Local in time Strong Solutions}
\label{section8}

Now we would like to obtain a uniform estimate of $u^n, \rho^n,r^n$
on $n$, which assures the existence of limit functions $u,\rho, r$.
From (\ref{rhon+1}), in order to get a uniform bound of approximate
densities, one needs to control $L^{\infty}$-norm of
$\rho^n(r^n)^2D_xu^{n+1}+\frac{2u^{n+1}}{r^n}$. For given $T$
sufficiently small, define $M^{n+1}$ by
\[
M^{n+1}\equiv \sup_{0\leq t\leq T;0\leq x\leq
1}|\rho^n(r^n)^2D_xu^{n+1}+\frac{2u^{n+1}}{r^n}|.
\]
Our first goal is to get a uniform bound of $M^{n+1}$ in the view of
\[
|\rho^n(r^n)^2D_xu^{n+1}+\frac{2u^{n+1}}{r^n}|\leq
\int_0^1|\rho^n(r^n)^2D_xu^{n+1}+\frac{2u^{n+1}}{r^n}|dx
+\int_0^1|D_x(\rho^n(r^n)^2D_xu^{n+1}+\frac{2u^{n+1}}{r^n})|dx.
\]
Note that the second term of the RHS may not be defined due to the
singularity of the origin. It is desirable to introduce cutoff
functions and to make use of both Lagrangian and Eulerian estimates.
To do so, we first look at the energy estimates of $u^{n+1}$
($t$-derivatives first). We use $\chi$ and $\zeta$ as cutoff
functions as in the a priori estimates. Before going any further,
define the separated $(n+1)$-th energies and the $(n+1)$-th
dissipation, resembling $\mathcal{E}$ and $\mathcal{D}$ in Section
\ref{section1}, as follows:

\[
\begin{split}
\mathcal{F}^{n+1}(t)\equiv&
\mathcal{F}^{n+1}_L(t)+\mathcal{F}^{n+1}_E(t)
\\
\equiv&\frac{1}{2}\sum_{i=0}^3\int_{x_0}^1\chi|D_t^iu^{n+1}|^2dx+
\frac{\mu}{2}\sum_{i=0}^2\int_{x_0}^1\chi\{\rho^n(r^n)^4|D_xD_t^iu^{n+1}|^2
+\frac{2|D_t^iu^{n+1}|^2}{\rho^n(r^n)^2}\}dx\\
+&\frac{1}{2}\sum_{i=0}^3\int_0^{r_2-d}\zeta\rho^n|\partial_t^iu^{n+1}|^2(r^n)^2dr^n+
\frac{\mu}{2}\sum_{i=0}^2\int_0^{r_2-d}\zeta\{|\partial_{r^n}\partial_t^iu^{n+1}|^2
+\frac{2|\partial_t^iu^{n+1}|^2}{(r^n)^2}\}(r^n)^2dr^n\\
\
\\
\mathcal{H}^{n+1}(t)\equiv&\mathcal{H}_L^{n+1}(t)+\mathcal{H}_E^{n+1}(t)
\\
\equiv&\frac{1}{\gamma-1}\li\chi(\rho^{n+1})^{\gamma-1}dx
+\frac{1}{2}\sum_{i=0}^2\int_{x_0}^1\chi
(\rho^{n+1})^{2\gamma-2}r^4|D_{x}D_t^i\rho^{n+1}|^2dx\\+&\frac{1}{2}
\sum_{i=0}^1\int_{x_0}^1\chi
(\rho^{n+1})^{4\gamma-2}r^8|D_{x}^2D_t^i\rho^{n+1}|^2dx+
\frac{1}{2}\int_{x_0}^1\chi
(\rho^{n+1})^{8\gamma-2}r^{12}|D_{x}^3\rho^{n+1}|^2dx
\\+&\frac{1}{2}\sum_{0\leq i+j\leq 3}\int_0^{r_2-d}\zeta
(\rho^{n+1})^{\gamma-2}|\partial_{r^{n}}^i\partial_t^j\rho^{n+1}|(r^{n})^2dr^{n}\\
\
\\
\mathcal{D}^{n+1}(t)\equiv&
\mu\sum_{i=0}^3\int_{x_0}^1\chi\{\rho^n(r^n)^4|D_xD_t^iu^{n+1}|^2
+\frac{2|D_t^iu^{n+1}|^2}{\rho^n(r^n)^2}\}dx+\sum_{i=1}^3\int_{x_0}^1\chi
|D_t^iu^{n+1}|^2dx\\
+&\mu\sum_{i=0}^3\int_0^{r_2-d}\zeta\{|\partial_{r^n}\partial_t^iu^{n+1}|^2
+\frac{2|\partial_t^iu^{n+1}|^2}{(r^n)^2}\}(r^n)^2dr^n+
\sum_{i=1}^3\int_0^{r_2-d}\zeta\rho^n|\partial_t^iu^{n+1}|^2(r^n)^2dr^n\\
\
\end{split}
\]
In the same spirit as in the a priori estimates, it is easy to show
that $M^n$'s are bounded by $\mathcal{F}^n$ and $\mathcal{H}^n$; our
main interest is to obtain a uniform bound. Note that $\rho^n$ and
$r^n$ can be controlled by $M^n$ for $0\leq t\leq T$:
\begin{equation}
\begin{split}
&\bullet\;\rho_{in} e^{-M^{n}T}\leq  \rho^{n}\leq \rho_{in} e^{M^{n}T}\\
&\bullet\;r_{in}e^{-\frac{M^{n}T}{3}}\leq  r^{n}\leq
r_{in}e^{\frac{M^{n}T}{3}}
\\ &\bullet \;\sup |\frac{D_t\rho^{n}}{\rho^{n}}|=|\rho^{n-1}(r^{n})^2D_xu^{n}
+\frac{2u^{n}}{r^{n-1}}|\leq M^{n}\\
&\bullet\;
D_tr^{n}=-\frac{1}{(r^{n})^2}\int_0^x\frac{D_t\rho^{n}}{\rho^{n}}
\frac{1}{\rho^{n}}dy\;\Rightarrow\; \sup
|\frac{D_tr^{n}}{r^{n}}|\leq \frac{M^{n}}{3}
\end{split}\label{nbound}
\end{equation}
In particular, the particle paths $r_i^n(t)$ emanating from $r_i$
for $i=0,1,2$ which correspond to $x_i$ can be bounded as follows:
\[
\begin{split}
&r_i^n(t)\equiv
(3\int_0^{x_i}\frac{1}{\rho^n}dy)^{\frac{1}{3}}\;\Rightarrow\;
r_{i}e^{-\frac{M^{n}T}{3}}\leq  r_i^{n}(t)\leq
r_{i}e^{\frac{M^{n}T}{3}}
\end{split}
\]
We observe that for each given $n$ and $d$, there exists a
sufficiently small $T>0$ such that $|e^{\frac{M^nT}{3}}-1|<d$, and
thus cutoff functions defined at the beginning of this article play
the same role as in the a priori estimates. We will show that
$M^n$'s are uniformly bounded and $T$ does not shrink to zero.
Before we state the result, let us speculate about the energy
defined in the above: note that $\mathcal{F}^{n+1}(t)$ does not
include some mixed, spatial derivatives of $u^{n+1}$ and
$\mathcal{H}_L^{n+1}$ does not include pure $t$-derivative terms,
which we need to close the estimates in the end. We show missing
derivative terms can be estimated in terms of $\mathcal{F}^{n+1}(t)$
and $\mathcal{H}^n(t)$ via the equations (\ref{al}) and (\ref{ae}).
Consider the following equations:
\[
\begin{split}
\mu D_x(\rho^n(r^n)^2D_xu^{n+1})=-\mu\frac{2D_xu^{n+1}}{r^n}+
\mu&\frac{2u^{n+1}}{\rho^n(r^n)^4}+\frac{D_tu^{n+1}}{(r^n)^2}+D_xp^n
+\frac{4\pi x}{(r^n)^4}\\
\Rightarrow\; 
\frac{\mu^2}{2}\int_{x_0}^1\chi
\rho^n|(r^n)^2D_x(\rho^n(r^n)^2D_xu^{n+1})|^2dx&\leq
C\frac{e^{\frac{2M^nT}{3}}}{r_0^2}\mathcal{F}_L^{n+1}(t)+C
(\sup_{x_0\leq x\leq 1}\rho_{in})e^{M^nT}\mathcal{H}^n(t)\\
&\leq C_{in,n}(\mathcal{F}^{n+1}(t)+\mathcal{H}^n(t)),
\end{split}
\]
where $C_{in,n}$ depends on initial data,
$\mathcal{F}^n,\mathcal{H}^n$. Similarly, take $D_t$ or $D_x$
derivatives of the above equation to get the estimates:
\[
\begin{split}
\frac{\mu}{2}\int_{x_0}^1\chi
\rho^n|(r^n)^2D_x(\rho^n(r^n)^2D_tD_{x}u^{n+1})|^2dx
\leq C_{in,n}(\mathcal{F}^{n+1}(t)+\mathcal{H}^n(t))\text{ and }\\
\;\mu^2\int_{x_0}^1\chi
|\rho^n(r^n)^2D_x(\rho^n(r^n)^3D_x(\rho^n(r^n)^2D_xu^{n+1}))|dx\leq
C_{in,n}(\mathcal{F}^{n+1}(t)+\mathcal{H}^n(t))
\end{split}
\]
For details, we refer to the a priori estimates. Similarly, $L^2$
norm of $\partial_{r^n}^2u^{n+1},\;\partial_t\partial_{r^n}^2u^{n+1}
,\;\partial_{r^n}^3u^{n+1}$ can be estimated in Eulerian coordinates
by using equations. For pure Lagrangian $t$-derivative terms, we
directly use the approximate continuity equation (\ref{rho}): for
$i=1,2,3,$
\[
\li\chi\frac{|D_t^i\rho^{n+1}|^2}{(\rho^{n+1})^3}dx\leq C
e^{(M^{n+1}+M^n)T}\mathcal{F}^{n+1}(t)
\]
where we have used $|\frac{\rho^{n+1}}{\rho^n}|\leq
e^{(M^{n+1}+M^n)T}$ for $0\leq t\leq T.$\\

\begin{prop} There exist $A_0,\;A_1,\;T^{\ast}>0$ such that if 
$\mathcal{F}(0)\leq A_0$, $\mathcal{H}(0)\leq A_1$ and\\
$\sup_{0\leq t\leq T^{\ast}}\mathcal{F}^n(t)\leq 2A_0$, $\sup_{0\leq
t\leq T^{\ast}}\mathcal{H}^n(t)\leq 2A_1$ then
\[
\sup_{0\leq t\leq T^{\ast}}\mathcal{F}^{n+1}(t)\leq 2A_0\text{ and
}\sup_{0\leq t\leq T^{\ast}}\mathcal{H}^{n+1}(t)\leq 2A_1.
\]
In particular, we can choose $T^{\ast}$ small enough so that there
exists $M>0$ such that   $M^{n}\leq M$ for all $n$, and for $0\leq t
\leq T^{\ast}$.\label{prop2}
\end{prop}\

To prove Proposition \ref{prop2}, we establish the following
chain-type energy inequalities. We start with $u$-part
energy $\mathcal{F}^{n+1}$.\\

\begin{lemma} Let $\rho^{n+1},u^{n+1},\rho^{n},r^n,u^n$ be regular
approximate solutions obtained in Section \ref{section6}.  Then
there exist constants $C_1,...,C_5>0$ such that for some positive
numbers $\epsilon_1,\epsilon_2>0$,
\begin{equation}
\frac{d}{dt}\mathcal{F}^{n+1}(t)+\frac{1}{8}\mathcal{D}^{n+1}(t)
\leq
C_1\{(M^n+M^{n+1})\mathcal{F}^{n+1}(t)+(\mathcal{F}^{n+1}(t))^2\}+
\mathcal{G}^{n+1}(t)+\mathcal{OL}_1^{n+1}(t),\label{e1}
\end{equation}
where
\begin{equation}
\mathcal{G}^{n+1}(t)\leq
C_2\{(1+M^{n+1}+|M^n|^2+\mathcal{F}^n(t))\mathcal{H}^n(t)+(\mathcal{H}^n(t))^{2}\}
+C_3e^{C_4(M^n+M^{n-1})T} \mathcal{F}^n(t),\label{e2}
\end{equation}
\begin{equation}
\mathcal{OL}_1^{n+1}(t)\leq
C_5\{\mathcal{F}^{n+1}(t)+\mathcal{H}^n(t)+
(\mathcal{F}^{n+1}(t))^{1+\epsilon_1}+(\mathcal{H}^n(t))^{1+\epsilon_2}\}.\label{e3}
\end{equation}\label{energyi}
\end{lemma}

 The estimate of $\rho$-part energy $\mathcal{H}^{n+1}$ is given in
 the next lemma.\\

\begin{lemma}Let $\rho^{n+1},u^{n+1},\rho^{n},r^n,u^n$ be regular
approximate solutions obtained in Section \ref{section6}.  Then
there exist constants $C_6,...,C_9>0$ such that for some positive
numbers $\epsilon_3,\epsilon_4>0$,
\begin{equation}
\begin{split}
\frac{d}{dt}\mathcal{H}^{n+1}(t)\leq
C_6(M^n+M^{n+1})\mathcal{H}^{n+1}(t)+C_7(\mathcal{F}^{n+1}(t)
+\mathcal{H}^n(t))(C_{in}e^{C_8M^{n+1}T}+\mathcal{H}^{n+1}(t))
\\+\mathcal{OL}_2^{n+1}(t), \label{e4}
\end{split}
\end{equation}
where
\begin{equation*}
\begin{split}
\mathcal{OL}_2^{n+1}(t)\leq
C_9\{\mathcal{H}^{n+1}(t)+\mathcal{H}^n(t)+(\mathcal{H}^{n+1}(t))^{1+\epsilon_3}+
(\mathcal{H}^{n}(t))^{1+\epsilon_4} \}.\label{e3}
\end{split}
\end{equation*}\label{energyii}
\end{lemma}

 We observe that it is enough to show the above energy
inequalities to prove Proposition \ref{prop2}:  Note that by
(\ref{Mn+1}) and the assumption, (\ref{e1}) can be reduced to the
following:
\[
\begin{split}
\frac{d}{dt}\mathcal{F}^{n+1}(t)\leq C_{M,t}\{\mathcal{F}^{n+1}(t)+
(\mathcal{F}^{n+1}(t))^2+A_0+A_1+(A_0+A_1)^2\}
\end{split}
\]
By solving the differential inequality, one can conclude that for
some $A_0,A_1$ and for small enough $T_1>0,$ $\sup_{0\leq t\leq
T_1}\mathcal{F}^{n+1}(t)\leq 2A_0$. Similarly, from (\ref{e4}), one
can deduce $\sup_{0\leq t\leq T_2}\mathcal{H}^{n+1}(t)\leq 2A_1$ for
small enough $T_2>0.$ Choose $T^{\ast}\equiv\min\{T_1,T_2\}$. The
moral of the proof of the above lemmas is same as in the a priori
estimates, but here, we separate the energy $\mathcal{E}$ into
$u$-part $\mathcal{F}$ and $\rho$-part $\mathcal{H}$; estimate
$\mathcal{F}$ first, and in turn $\mathcal{H}$
in accordance with the approximate scheme.\\

\begin{proof}\textit{of Lemma \ref{energyi}:} Now we are ready to
prove (\ref{e1}). The spirit is same as in the a priori estimates.
Eulerian estimates and Lagrangian estimates are to be performed
concurrently. We provide rather detailed estimates to distinguish
$\mathcal{F}^n,\mathcal{F}^{n+1},\mathcal{H}^n,\mathcal{H}^{n+1},\mathcal{G}^{n+1},
\mathcal{OL}$ etc.\\

\noindent {\bf Notation:} In the below, we use the double underline
 $\underline{\underline{\;\;\;\;}}$  to denote the terms that
 belong to $\mathcal{OL}_1^{n+1}$, and the under-brace
${\underbrace{\;\;\;\;}}$ to denote the terms that we take the sup
of. The rest of terms in the RHS will either contribute to
$\mathcal{G}^{n+1}(t)$ or be absorbed into the dissipation in the
LHS at
the last step.\\

Here is the zeroth order estimates. Multiply (\ref{al}) by $\chi
u^{n+1}$ and integrate it over $x$ to get:
\[
\begin{split}
\frac{1}{2}\frac{d}{dt}\int_{x_0}^1\chi|u^{n+1}|^2dx+\frac{3\mu}{4}\int_{x_0}^1
\chi\{\rho^n(r^n)^4|D_xu^{n+1}|^2+\frac{2|u^{n+1}|^2}{\rho^n(r^n)^2}\}dx
\leq\frac{5A^2}{\mu}\int_{x_0}^1\chi(\rho^n)^{2\gamma-1}dx\\
+\frac{16\pi^2}{\mu}\int_{x_0}^1\chi\rho^n\frac{x^2}{(r^n)^2}dx
+\underline{\underline{\int_{x_0}^{x_1}\chi'(r^n)^2p^nu^{n+1}dx-
\mu\int_{x_0}^{x_1}\chi'u^{n+1}\rho^n(r^n)^4D_xu^{n+1}dx}}
\end{split}
\]\\

\noindent Multiply (\ref{ae}) by $\zeta u^{n+1}(r^n)^2$ and
integrate it over $r^n$ to get:
\[
\begin{split}
\frac{1}{2}\frac{d}{dt}\int_0^{r_2-d}\zeta\rho^n|u^{n+1}|^2(r^n)^2dr^n
+\frac{3\mu}{4}\int_0^{r_2-d}\zeta\{|\partial_{r^n}u^{n+1}|^2+\frac{2|u^{n+1}|^2}
{(r^n)^2}\}(r^n)^2dr^n\\\leq\frac{3A^2}{\mu}\int_0^{r_2-d}\zeta(\rho^n)^{2\gamma}
(r^n)^2dr^n+\frac{16\pi^2}{\mu}\int_0^{r_2-d}\zeta\frac{x^2}{(r^n)^2}(\rho^n)^2
(r^n)^2dr^n\\+\frac{2}{\mu}\int_0^{r_2-d}\zeta\underbrace{\rho^n|D_tr^n|^2}\rho^n
|u^{n+1}|^2
(r^n)^2dr^n+\frac{1}{2}\int_0^{r_2-d}\zeta\underbrace{\frac{\partial_t\rho^n}
{\rho^n}}\rho^n|u^{n+1}|^2(r^n)^2dr^n\\
\underline{\underline{-\mu\int_{r_1+d}^{r_2-d}(\partial_{r^n}\zeta)
(\partial_{r^n}u^{n+1})
u^{n+1}(r^n)^2dr^n+\int_{r_1+d}^{r_2-d}(\partial_{r^n}\zeta)p^n
u^{n+1}(r^n)^2dr^n}}
\end{split}
\]\\

\noindent Next, we estimate one spatial derivative of $u^{n+1}$.
Multiply (\ref{al}) by $\chi D_tu^{n+1}$ and integrate it over $x$
to get:
\[
\begin{split}
\frac{\mu}{2}\frac{d}{dt}\int_{x_0}^1\chi\{\rho^n
(r^n)^4|D_xu^{n+1}|^2+\frac{2|u^{n+1}|^2}{\rho^n (r^n)^2}\}dx
+\frac{3}{4}\int_{x_0}^1\chi|D_tu^{n+1}|^2 dx\\
\leq\frac{\mu}{2}\int_{x_0}^1\chi
\underbrace{\frac{D_t(\rho^n(r^n)^4)}{\rho^n(r^n)^4}}\rho^n(r^n)^4|D_xu^{n+1}|^2dx
+\mu\int_{x_0}^1\chi
\underbrace{\rho^n(r^n)^2D_t(\frac{1}{\rho^n(r^n)^2})}\frac{|u^{n+1}|^2}
{\rho^n(r^n)^2}dx\\
+ \int_{x_0}^1\chi
(r^n)^4|D_xp^n|^2dx+\frac{\mu}{4}\li\chi\frac{|D_tu^{n+1}|^2}{\rho^n(r^n)^2}dx
+\frac{1}{\mu}\int_{x_0}^1\chi\frac{16\pi^2 x^2\rho^n}{(r^n)^2}dx\\
\underline{\underline{-\mu\int_{x_0}^{x_1}\chi' \rho^n(r^n)^4D_x
u^{n+1}D_tu^{n+1}dx}}
\end{split}
\]\\

\noindent Multiply (\ref{ae}) by $\zeta
\rho^n\partial_tu^{n+1}(r^n)^2$ and integrate it over $r^n$ to get:
\[
\begin{split}
\frac{\mu}{2}\frac{d}{dt}\int_0^{r_2-d}\zeta\{|\partial_{r^n}u^{n+1}|^2
+\frac{2|u^{n+1}|^2}{(r^n)^2}\}(r^n)^2dr^n+\frac{3}{4}\int_0^{r_2-d}\zeta\rho^n
|\partial_tu^{n+1}|^2(r^n)^2dr^n\\
\leq
6\int_0^{r_2-d}\zeta\underbrace{\rho^n|D_tr^n|^2}|\partial_{r^n}u^{n+1}|^2
(r^n)^2dr^n
+6\int_0^{r_2-d}\frac{1}{\rho^n}|\partial_{r^n}p^n|^2(r^n)^2dr^n\\
+6\int_0^{r_2-d}\zeta\frac{16\pi^2 x^2}{(r^n)^2}\rho^ndr^n
\underline{\underline{-\mu\int_{r_1+d}^{r_2-d}(\partial_{r^n}\zeta)
(\partial_t u^{n+1})(\partial_{r^n}u^{n+1})(r^n)^2dr^n}}
\end{split}
\]\\

\noindent Next, in order to estimate one temporal derivative of
$u^{n+1}$, take $D_t$ of (\ref{al}), multiply it by $\chi
D_tu^{n+1}$, and integrate to get:
\[
\begin{split}
\frac{1}{2}\frac{d}{dt}\int_{x_0}^1\chi|D_tu^{n+1}|^2dx+\frac{3\mu}{4}\int_{x_0}^1\chi
\{\rho^n(r^n)^4|D_tD_{x}u^{n+1}|^2+\frac{2|D_tu^{n+1}|^2}{\rho^n(r^n)^2}\}dx\\
\leq
2\mu\int_{x_0}^1\chi\underbrace{|\frac{D_t(\rho^n(r^n)^4)}{\rho^n
(r^n)^4}|^2}\rho^n(r^n)^4|D_xu^{n+1}|^2dx+8\mu\int_{x_0}^1\chi
\underbrace{|\rho^n(r^n)^2D_t(\frac{1}{\rho^n(r^n)^2})|^2}
\frac{|u^{n+1}|^2}{\rho^n(r^n)^2}dx\\+\frac{2}{\mu}\int_{x_0}^1\chi
{\frac{|D_t((r^n)^2p^n)|^2}
{\rho^n(r^n)^4}}dx+\frac{2}{\mu}\int_{x_0}^1\chi{\rho^n(r^n)^2|D_t(\frac{2p^n}
{\rho^nr^n})|^2}dx+\frac{32\pi^2}{\mu}\int_{x_0}^1\chi\rho^n(r^n)^2|D_t(\frac{x}
{(r^n)^2})|^2dx\\
\underline{\underline{-\mu\int_{x_0}^{x_1}\chi'\rho^n(r^n)^2D_tu^{n+1}D_{x}D_tu^{n+1}dx-\mu\int_{x_0}^{x_1}
\chi'D_t(\rho^n(r^n)^4)D_tu^{n+1}D_xu^{n+1}dx}}\\
\underline{\underline{+\int_{x_0}^{x_1}\chi'D_t((r^n)^2p^n)
D_tu^{n+1}dx}}
\end{split}
\]\\

\noindent Likewise, take $\partial_t$ of (\ref{ae}), multiply it by
$\zeta\rho^n \partial_tu^{n+1}(r^n)^2$, and integrate to get:
\[
\begin{split}
\frac{1}{2}\frac{d}{dt}\int_0^{r_2-d}\zeta\rho^n|\partial_tu^{n+1}|^2(r^n)^2dr^n
+\frac{3\mu}{4}\int_0^{r_2-d}\zeta\{|\partial_t\partial_{r^n}u^{n+1}|^2
+\frac{2|\partial_tu^{n+1}|^2}{(r^n)^2}\}(r^n)^2dr^n\\
\leq\frac{2}{\mu}\int_0^{r_2-d}\zeta\underbrace{\rho^n|D_tr^n|^2}\rho^n
|\partial_tu^{n+1}|^2(r^n)^2dr^n
+\frac{3}{\mu}\int_0^{r_2-d}\zeta|\partial_tp^n|^2(r^n)^2dr^n
\\
-\frac{1}{2}\int_0^{r_2-d}\zeta\underbrace{\frac{\partial_t\rho^n}{\rho^n}}
\rho^n|\partial_tu^{n+1}|^2(r^n)^2dr^n
-\frac{1}{\mu}\int_0^{r_2-d}\zeta\underbrace{\frac{|\partial_t(\rho^nD_tr^n)|^2}
{\rho^n}} \rho^n|\partial_tu^{n+1}|^2
(r^n)^2dr^n\\
+\frac{\mu}{4}\int_0^{r_2-d}\zeta|\partial_{r^n}u^{n+1}|^2
(r^n)^2dr^n+\frac{32\pi^2}{\mu}\int_0^{r_2-d}\zeta
x^2|\partial_t\rho^n|^2dr^n+
\frac{32\pi^2}{\mu}\int_0^{r_2-d}\zeta(\rho^n)^4(r^n)^4|D_tr^n|^2dr^n\\
\underline{\underline{-\mu\int_{r_1+d}^{r_2-d}(\partial_{r^n}\zeta)\partial_t\partial_{r^n}u^{n+1}
\partial_tu^{n+1}(r^n)^2dr^n-\int_{r_1+d}^{r_2-d}(\partial_{r^n}\zeta)\partial_tp^n
\partial_{t}u^{n+1}(r^n)^2dr^n}}
\end{split}
\]\\

\noindent Next we estimate $D_tD_x u^{n+1}$. Multiply
$D_t$(\ref{al}) by $\chi D_t^2u^{n+1}$ and integrate to get:
\[
\begin{split}
\frac{\mu}{2}\frac{d}{dt}\int_{x_0}^1\chi\{\rho^n(r^n)^4|D_{t}D_xu^{n+1}|^2
+\frac{2|D_tu^{n+1}|^2}{\rho^n(r^n)^2}\}dx+\frac{3}{4}\int_{x_0}^1\chi
|D_{t}^2u^{n+1}|^2dx\\ \leq \frac{\mu}{4}\li\chi\{\rho^n(r^n)^4
|D_t^2D_xu^{n+1}|^2+\frac{2|D_t^2u^{n+1}|^2}{\rho^n(r^n)^2} \}dx+\mu
\li\chi
\underbrace{|\frac{D_t(\rho^n(r^n)^4)}{\rho^n(r^n)^4}|^2}\rho^n(r^n)^4|D_xu^{n+1}|^2dx \\
+\mu\li\chi\underbrace{|\rho^n(r^n)^2D_t(\frac{1}{\rho^n(r^n)^2})|^2}
\frac{2|u^{n+1}|^2}{\rho^n(r^n)^2}dx+\frac{\mu}{2}\li\chi
\underbrace{\frac{D_t(\rho^n(r^n)^4)}{\rho^n(r^n)^4}}
\rho^n(r^n)^4|D_tD_xu^{n+1}|^2dx \\
+\frac{\mu}{2}\li\chi
\underbrace{\rho^n(r^n)^2D_t(\frac{1}{\rho^n(r^n)^2})}
\frac{2|D_tu^{n+1}|^2}{\rho^n(r^n)^2}dx
+\li\chi|D_t((r^n)^2D_xp^n)|^2dx+\frac{1}{\mu}\li\chi
16\pi^2\rho^n(r^n)^2 x^2\cdot\\|D_t(\frac{1}{(r^n)^2})|^2dx
\underline{\underline{-\mu\int_{x_0}^{x_1}\chi'\rho^n(r^n)^4D_tD_{x}u^{n+1}
D_{t}^2u^{n+1}dx-\mu\int_{x_0}
^{x_1}\chi'D_t(\rho^n(r^n)^4)D_xu^{n+1}D_{t}^2u^{n+1}dx}}
\end{split}
\]\\

\noindent Corresponding Eulerian estimates for
$\partial_t\partial_{r^n}u^{n+1}$ are given as follows:
\[
\begin{split}
\frac{\mu}{2}\frac{d}{dt}\int_0^{r_2-d}\zeta\{|\partial_t\partial_{r^n}u^{n+1}|^2
+\frac{2|\partial_tu^{n+1}|^2}{(r^n)^2}\}(r^n)^2dr^n+\frac{3}{4}\int_0^{r_2-d}\zeta\rho^n
|\partial_t^2u^{n+1}|^2(r^n)^2dr^n\\
\leq
4\int_0^{r_2-d}\zeta\underbrace{|\frac{\partial_t\rho^n}{\rho^n}|^2}\rho^n
|\partial_tu^{n+1}|^2(r^n)^2dr^n
+8\int_0^{r_2-d}\zeta\underbrace{|\rho^nD_tr^n|^2}|\partial_t\partial_{r^n}u^{n+1}|^2(r^n)^2dr^n\\+
8\int_0^{r_2-d}\zeta\underbrace{\frac{|\partial_t(\rho^nD_tr^n)|^2}{\rho^n}}|
\partial_{r^n}u^{n+1}|^2(r^n)^2dr^n+4\int_0^{r_2-d}\zeta\frac{1}{\rho^n}
|\partial_t\partial_{r^n}p^n|^2(r^n)^2dr^n\\
+64\pi^2\int_0^{r_2-d}\zeta\frac{|\partial_t(\rho^n
x)|^2}{\rho^n(r^n)^2}dr^n
\underline{\underline{-\mu\int_{r_1+d}^{r_2-d}(\partial_{r^n}\zeta)(\partial_t\partial_{r^n}u^{n+1})
\partial_t^2u^{n+1}(r^n)^2dr^n}}
\end{split}
\]\\

\noindent Similarly, we can perform the higher order energy
estimates. Here is the estimate of $D_t^2u^{n+1}$. Multiply
$D_t^2$(\ref{al}) by $\chi D_t^2u^{n+1}$ and integrate:
\[
\begin{split}
\frac{1}{2}\frac{d}{dt}\int_{x_0}^1\chi|D_{t}^2u^{n+1}|^2dx+\frac{3\mu}{4}\int_{x_0}^1\chi
\{\rho^n(r^n)^4|D_t^2D_{x}u^{n+1}|^2+\frac{2|D_{t}^2u^{n+1}|^2}{\rho^n(r^n)^2}\}dx\\
\leq
12\mu\int_{x_0}^1\chi\underbrace{|\frac{D_t(\rho^n(r^n)^4)}{\rho^n(r^n)^4}|^2}
\rho^n(r^n)^4|D_tD_{x}u^{n+1}|^2dx
+3\mu\int_{x_0}^1\chi\frac{|D_{t}^2(\rho^n(r^n)^4)|^2}{(\rho^n)^3
(r^n)^8}\underbrace{|\rho^n(r^n)^2D_xu^{n+1}|^2}dx\\
+32\mu
\int_{x_0}^1\chi\underbrace{|\rho^n(r^n)^2D_t(\frac{1}{\rho^n(r^n)^2})|^2}
\frac{|D_tu^{n+1}|^2}{\rho^n(r^n)^2}dx+8\mu\int_{x_0}^{1}\chi
\rho^n(r^n)^2|D_{t}^2(\frac{1}{\rho^n(r^n)^2})|^2
\underbrace{|u^{n+1}|^2}dx\\+\frac{3}{\mu}\int_{x_0}^1
\chi\frac{|D_{t}^2((r^n)^2p^n)|^2}{\rho^n(r^n)^4}dx
+\frac{2}{\mu}\int_{x_0}^1\chi\rho^n(r^n)^2|D_{t}^2(\frac{2p^n}{\rho^nr^n})|^2dx+
\frac{32\pi^2}{\mu}\int_{x_0}^1\chi\rho^n(r^n)^2|D_{t}^2(\frac{x}{(r^n)^2})|^2dx
\\ \underline{\underline{-\mu\int_{x_0}^{x_1}\chi'\rho^n(r^n)^4D_t^2D_{x}u^{n+1}
D_{t}^2u^{n+1}dx
-2\mu\int_{x_0}^{x_1}\chi'D_t(\rho^n(r^n)^4)D_tD_{x}u^{n+1}D_{t}^2u^{n+1}dx}}\\
\underline{\underline{-\mu\int_{x_0}^{x_1}\chi'D_{t}^2(\rho^n(r^n)^4)D_xu^{n+1}
D_{t}^2u^{n+1}dx+\int_{x_0}^{x_1}\chi'D_{t}^2((r^n)^2p^n)D_{t}^2u^{n+1}dx}}
\end{split}
\]\\

\noindent The Eulerian estimate of $\partial_t^2u^{n+1}$: multiply
$\partial_t^2$(\ref{ae}) by $\zeta\rho^n\partial_t^2u^{n+1}(r^n)^2$
and integrate:
\[
\begin{split}
\frac{1}{2}\frac{d}{dt}\int_0^{r_2-d}\zeta\rho^n|\partial_t^2u^{n+1}|^2(r^n)^2dr^n
+\frac{3\mu}{4}\int_0^{r_2-d}\chi\{|\partial_t^2\partial_{r^n}u^{n+1}|^2
+\frac{2|\partial_t^2u^{n+1}|^2}{(r^n)^2}\}(r^n)^2dr^n\\
\leq-\frac{3}{2}\int_0^{r_2-d}\zeta
\underbrace{\frac{\partial_t\rho^n}{\rho^n}}\rho^n|\partial_t^2u^{n+1}|^2(r^n)^2dr^n
+\frac{2}{\mu}\int_0^{r_2-d}\zeta\underbrace{|\rho^nD_tr^n|^2}
|\partial_t^2u^{n+1}|^2(r^n)^2dr^n\\
+\int_0^{r_2-d}\zeta\underbrace{|\partial_tu^{n+1}|^2}
\rho^n|\partial_t^2u^{n+1}|^2(r^n)^2dr^n+\int_0^{r_2-d}
\underbrace{|\partial_{r^n}u^{n+1}|^2}\rho^n
|\partial_t^2u^{n+1}|^2(r^n)^2dr^n
\\+\int_0^{r_2-d}\zeta \frac{(\partial_t^2\rho^n)^2}{\rho^n}(r^n)^2dr^n
+\frac{10}{\mu}\int_0^{r_2-d}\zeta|\partial_t^2p^n|^2(r^n)^2dr^n+\frac{16\pi^2}{\mu}
\int_0^{r_2-d}\zeta |\partial_t^2(\rho^nx)|^2dr^n\\
+\int_0^{r_2-d}\zeta\frac{|\partial_t^2(\rho^nD_tr^n)|^2}{\rho^n}
(r^n)^2dr^n+\frac{4}{\mu}\int_0^{r_2-d}\zeta
\underbrace{\frac{|\partial_t(\rho^nD_tr^n)|^2}{\rho^n}}
\rho^n|\partial_t^2u^{n+1}|^2(r^n)^2dr^n\\
\underline{\underline{-\mu\int_{r_1+d}^{r_2-d}(\partial_{r^n}\zeta)\partial_t^2u^{n+1}\partial_t^2
\partial_{r^n}u^{n+1}(r^n)^2dr^n+\int_{r_1+d}^{r_2-d}(\partial_{r^n}\zeta)\partial_t^2
p^n\partial_t^2u^{n+1}(r^n)^2dr^n}}
\end{split}
\]\\

\noindent Now we turn into $D_t^2D_xu^{n+1}$. Multiply
$D_t^2$(\ref{al}) by $\chi D_t^3u^{n+1}$ and integrate:
\[
\begin{split}
\frac{\mu}{2}\frac{d}{dt}\int_{x_0}^1\chi\{\rho^n(r^n)^4|D_t^2D_{x}u^{n+1}|^2
+\frac{2|D_{t}^2u^{n+1}|^2}{\rho^n(r^n)^2}\}dx+\frac{3}{4}
\int_{x_0}^1\chi|D_{t}^3u^{n+1}|^2dx\\
\leq\frac{\mu}{4}\li\chi\{\rho^n(r^n)^4
|D_t^3D_xu^{n+1}|^2+\frac{2|D_t^3u^{n+1}|^2}{\rho^n(r^n)^2} \}dx
+\frac{\mu}{2}\int_{x_0}^1\chi
\underbrace{\frac{D_t(\rho^n(r^n)^4)}{\rho^n(r^n)^4}}
\rho^n(r^n)^4|D_{t}^2D_xu^{n+1}|^2dx\\
+\frac{\mu}{2}\li\chi\underbrace{\rho^n(r^n)^2D_t(\frac{1}{\rho^n(r^n)^2})}
\frac{2|D_t^2u^{n+1}|^2}{\rho^n(r^n)^2}dx+8\mu\li\chi
\underbrace{|\frac{D_t(\rho^n(r^n)^4)}{\rho^n(r^n)^4}|^2}\rho^n(r^n)^4
|D_tD_xu^{n+1}|^2dx
\\
+2\mu\int_{x_0}^1\chi\frac{|D_{t}^2(\rho^n(r^n)^4)|^2}{(\rho^n)^3
(r^n)^8}\underbrace{|\rho^n(r^n)^2D_xu^{n+1}|^2}dx+4\mu\li\chi
\underbrace{|\rho^n(r^n)^2D_t(\frac{1}{\rho^n(r^n)^2})|^2}
\frac{2|D_tu^{n+1}|^2}{\rho^n(r^n)^2}dx\\
+8\mu\li\chi\rho^n(r^n)^2|D_t^2(\frac{1}{\rho^n(r^n)^2})|^2
\underbrace{|u^{n+1}|^2}dx
+2\int_{x_0}^1\chi|D_{t}^2((r^n)^2D_xp^n)|^2dx\\+32\pi^2\int_{x_0}^1
\chi x^2|D_{t}^2(\frac{1}{(r^n)^2})|^2dx
\underline{\underline{-\mu\int_{x_0}^{x_1}\chi'\rho^n(r^n)^4D_t^2D_{x}
u^{n+1}D_{t}^3u^{n+1}dx}}\\
\underline{\underline{-2\mu
\int_{x_0}^{x_1}\chi'D_t(\rho^n(r^n)^4)D_tD_{x}u^{n+1}D_{t}^3u^{n+1}dx
-\mu\int_{x_0}^{x_1}\chi'D_{t}^2(\rho^n(r^n)^4)D_xu^{n+1}
D_{t}^3u^{n+1}dx}}
\end{split}
\]\\

\noindent For $\partial_t^2\partial_{r^n}u^{n+1}$, multiply
$\partial_t^2$(\ref{ae}) by $\zeta\rho^n\partial_t^3u^{n+1}(r^n)^2$
and integrate:
\[
\begin{split}
\frac{\mu}{2}\frac{d}{dt}\int_0^{r_2-d}\zeta\{|\partial_t^2\partial_{r^n}u^{n+1}|^2
+\frac{2|\partial_t^2u^{n+1}|^2}{(r^n)^2}\}(r^n)^2dr^n+\frac{1}{2}\int_0^{r_2-d}\zeta
\rho^n|\partial_t^3u^{n+1}|^2(r^n)^2dr^n\\
\leq
16\int_0^{r_2-d}\zeta\underbrace{|\frac{\partial_t\rho^n}{\rho^n}|^2}
\rho^n |\partial_t^2u^{n+1}|(r^n)^2dr^n+4\int_0^{r_2-d}\zeta
\frac{|\partial_t^2\rho^n|^2}{\rho^n}
\underbrace{|\partial_tu^{n+1}|^2}(r^n)^2dr^n\\
+8\int_0^{r_2-d}\zeta\underbrace{\rho^n|D_tr^n|^2}
|\partial_t^2\partial_{r^n}u^{n+1}|^2(r^n)^2dr^n
+32\int_0^{r_2-d}\zeta\underbrace{\frac{|\partial_t(\rho^nD_tr^n)|^2}{\rho^n}}
|\partial_t\partial_{r^n}u^{n+1}|^2(r^n)^2dr^n\\+8\int_0^{r_2-d}\zeta
\frac{|\partial_t^2(\rho^nD_tr^n)|^2}
{\rho^n}\underbrace{|\partial_{r^n}u^{n+1}|^2}(r^n)^2dr^n+2\int_0^{r_2-d}\zeta
\frac{|\partial_t^2\partial_{r^n}p^n|^2}{\rho^n}(r^n)^2dr^n\\+32\pi^2
\int_0^{r_2-d}\zeta\frac{|\partial_t^2(\rho^n
x)|^2}{\rho^n(r^n)^2}dr^n
\underline{\underline{-\mu\int_{r_1+d}^{r_2-d}(\partial_{r^n}\zeta)(\partial_t^2\partial_{r^n}u^{n+1})
(\partial_{t}^3u^{n+1})(r^n)^2dr^n}}
\end{split}
\]\\

\noindent Here is the estimate of the full time derivative terms.
Multiply $D_t^3$(\ref{al}) by $\chi D_t^3u^{n+1}$ and integrate to
get:
\[
\begin{split}
\frac{1}{2}\frac{d}{dt}\int_{x_0}^1\chi|D_{t}^3u^{n+1}|^2dx+\frac{3\mu}{4}
\int_{x_0}^1\chi
\{\rho^n(r^n)^4|D_t^3D_{x}u^{n+1}|^2+\frac{2|D_{t}^3u^{n+1}|^2}{\rho^n(r^n)^2}\}dx
\leq
36\mu\int_{x_0}^1\chi\cdot\\
\underbrace{|\frac{D_t(\rho^n(r^n)^4)}{\rho^n(r^n)^4}|^2}
\rho^n(r^n)^4 |D_t^2D_{x}u^{n+1}|^2dx
+36\mu\int_{x_0}^1\chi\frac{|D_{t}^2(\rho^n(r^n)^4)|^2}
{(\rho^n)^3(r^n)^8} \underbrace{|\rho^n(r^n)^2D_tD_{x}u^{n+1}|^2}dx
\\+4\mu\int_{x_0}^1\chi\frac{|D_{t}^3(\rho^n(r^n)^4)|^2}{(\rho^n)^3(r^n)^8}
\underbrace{|\rho^n(r^n)^2D_{x}u^{n+1}|^2}dx
+\frac{4}{\mu}\int_{x_0}^1\chi\frac{|D_{t}^3((r^n)^2p^n)|^2}{\rho^n(r^n)^4}dx\\
+90\mu\int_{x_0}^1\chi\underbrace{|\rho^{n}(r^n)^2D_t(\frac{1}{\rho^n(r^n)^2})|^2}
\frac{|D_{t}^2u^{n+1}|^2}{\rho^n(r^n)^2}dx
+90\mu\int_{x_0}^1\chi\rho^{n}(r^n)^2|D_{t}^2(\frac{1}{\rho^n(r^n)^2})|^2
\underbrace{|D_{t}u^{n+1}|^2}dx\\+10\mu\int_{x_0}^1\chi\rho^{n}(r^n)^2
|D_{t}^3(\frac{1}{\rho^n(r^n)^2})|^2\underbrace{|u^{n+1}|^2}dx+\frac{5}{2\mu}
\int_{x_0}^1\chi\rho^n(r^n)^2|D_{t}^3(\frac{2p^n}{\rho^nr^n})|^2dx\\
+\frac{40\pi^2}{\mu}\int_{x_0}^1\chi\rho^n(r^n)^2x^2|D_{t}^3(\frac{1}{(r^n)^2})|^2dx
\underline{\underline{-\mu\int_{x_0}^{x_1}\chi'
\rho^n(r^n)^4D_t^3D_{x}u^{n+1}D_{t}^3u^{n+1}dx}}\\
\underline{\underline{-3\mu\int_{x_0}^{x_1}\chi'
D_t(\rho^n(r^n)^4)D_t^2D_{x}u^{n+1}D_{t}^3u^{n+1}dx
-3\mu\int_{x_0}^{x_1}\chi'D_{t}^2(\rho^n(r^n)^4)D_tD_{x}u^{n+1}D_{t}^3u^{n+1}dx}}\\
\underline{\underline{-\mu\int_{x_0}^{x_1}\chi'D_{t}^3(\rho^n(r^n)^4)D_xu^{n+1}
D_{t}^3u^{n+1}dx+\int_{x_0}^{x_1}\chi'D_{t}^3((r^n)^2p^n)D_{t}^3u^{n+1}dx}}
\end{split}
\]\\

\noindent Finally, we are ready to the Eulerian estimate of
$\partial_t^3u^{n+1}$. In order to control the nonlinear terms
involving higher derivatives, displaying three-dimensional feature,
we employ the Gagliardo-Nirenberg inequality. Multiply
$\partial_t^3$(\ref{ae}) by $\zeta\rho^n\partial_t^3u^{n+1}(r^n)^2$
and integrate to get
\[
\begin{split}
\frac{1}{2}\frac{d}{dt}\int_0^{r_2-d}\zeta\rho^n|\partial_t^3u^{n+1}|^2(r^n)^2dr^n
+\frac{\mu}{2}\int_{0}^{r_2-d}\zeta\{|\partial_t^3\partial_{r^n}u^{n+1}|^2
+\frac{2|\partial_t^3u^{n+1}|^2}{(r^n)^2}\}(r^n)^2dr^n\\
\leq-\frac{5}{2}\int_0^{r_2-d}\zeta\underbrace{\frac{\partial_t\rho^n}
{\rho^n}}\rho^n|\partial_t^3u^{n+1}|^2(r^n)^2dr^n+\frac{2}{\mu}\int_0^{r_2-d}\zeta
\underbrace{\rho^n|D_tr^n|^2}\rho^n|\partial_t^3u^{n+1}|^2(r^n)^2dr^n\\
+36\int_0^{r_2-d}\zeta\underbrace{\frac{|\partial_t(\rho^nD_tr^n)|^2}{\rho^n}}
|\partial_t^2\partial_{r^n}u^{n+1}|^2(r^n)^2dr^n+4\int_0^{r_2-d}\zeta
\frac{|\partial_t^3(\rho^nD_tr^n)|^2}{\rho^n}\underbrace{|\partial_{r^n}u^{n+1}|^2}
(r^n)^2dr^n
\\
+\frac{1}{8}\int_0^{r_2-d}\zeta\rho^n|\partial_t^3u^{n+1}|^2(r^n)^2dr^n
+4\int_0^{r_2-d}\zeta\frac{|\partial_t^3\rho^n|^2}{\rho^n}
\underbrace{|\partial_tu^{n+1}|^2}(r^n)^2dr^n\\
+\frac{6}{\mu}\int_0^{r_2-d}\zeta|\partial_t^3p^n|^2(r^n)^2dr^n+\frac{64\pi^2}{\mu}
\int_0^{r_2-d}\zeta|\partial_t^3(\rho^n x)|^2dr^n
\\
\underbrace{-3\int_0^{r_2-d}\zeta\partial_t^2\rho^n\partial_t^2u^{n+1}
\partial_t^3u^{n+1} (r^n)^2dr^n+3\int_0^{r_2-d}\zeta
\partial_t^2(\rho^nD_tr^n)\partial_t\partial_{r^n}u^{n+1}\partial_t^3u^{n+1}
(r^n)^2dr^n}_{\star\star}\\
 \underline{\underline{-\mu\int_{r_1+d}^{r_2-d}(\partial_{r^n}\zeta)
(\partial_t^3u^{n+1})
(\partial_t^3\partial_{r^n}u^{n+1})(r^n)^2dr^n+\int_{r_1+d}^{r_2-d}
(\partial_{r^n}\zeta)(\partial_t^3p^n)(\partial_t^3u^{n+1})(r^n)^2dr^n}}
\end{split}
\]\\

\noindent For $\star\star$, first by the Cauchy-Schwarz inequality,
we write it as
\[
\begin{split}
\star\star\leq\frac{1}{8}\int_0^{r_2-d}\zeta\rho^n
|\partial_t^3u^{n+1}|^2(r^n)^2dr^n+6\int_0^{r_2-d}\zeta
|\frac{\partial_t^2\rho^n}{\sqrt{\rho^n}}|^4(r^n)^2dr^n+6\int_0^{r_2-d}\zeta
|\partial_t^2u^{n+1}|^4(r^n)^2dr^n\\
+6\int_0^{r_2-d}\zeta|\frac{\partial_t^2(\rho^nD_tr^n)}{\sqrt{\rho^n}}|^4
(r^n)^2dr^n+6\int_0^{r_2-d}\zeta
|\partial_t\partial_{r^n}u^{n+1}|^4(r^n)^2dr^n
\end{split}
\]\\

\noindent The first term will be absorbed into the LHS  and we apply
the Gagliardo-Nirenberg inequality $$||f||_{L^4}\leq
\frac{1}{2}||\nabla
f||_{L^2}^\frac{3}{4}||f||_{L^2}^{\frac{1}{4}}\;\text{ in
}\mathbb{R}^3$$ to the rest of terms. For instance, we have the
following:
\[
\int_0^{r_2-d}\zeta
|\frac{\partial_t^2\rho^n}{\sqrt{\rho^n}}|^4(r^n)^2dr^n\leq
(\int_0^{r_2-d}\zeta|\partial_{r^n}(\frac{\partial_t^2\rho^n}
{\sqrt{\rho^n}})|^2(r^n)^2dr^n)^{\frac{3}{2}} (\int_0^{r_2-d}\zeta
|\frac{\partial_t^2\rho^n} {\sqrt{\rho^n}}|^2(r^n)^2dr^n
)^{\frac{1}{2}}
\]
\[
\int_0^{r_2-d}\zeta
|\partial_t\partial_{r^n}u^{n+1}|^4(r^n)^2dr^n\leq
(\int_0^{r_2-d}\zeta
|\partial_t\partial_{r^n}^2u^{n+1}|^2(r^n)^2dr^n)^{\frac{3}{2}}
(\int_0^{r_2-d}\zeta
|\partial_t\partial_{r^n}u^{n+1}|^2(r^n)^2dr^n)^{\frac{1}{2}}
\]\\

\noindent Now we combine all the estimates that we have obtained so
far to get:
\[
\begin{split}
\frac{d}{dt}\mathcal{F}^{n+1}(t)+\mathcal{D}^{n+1}(t)\leq &
\frac{7}{8}\mathcal{D}^{n+1}(t)+C\sup|\underbrace{\{\text{under-braced
terms}\}}|\mathcal{F}^{n+1}(t)\\&+ C(\mathcal{F}^{n+1}(t))^2
+\mathcal{G}^{n+1}(t)+\underline{\underline{\{\text{double
underlined terms}\}}}
\end{split}
\]\\

\noindent Next, we list the terms of energy type of the  previous
n-th step: $\mathcal{G}^{n+1}(t)$.
\[
\begin{split}
\sum_{i=0}^3\li\chi\frac{|D_t^i((r^n)^2p^n)|^2}{\rho^n(r^n)^4}dx,\;
\sum_{i=0}^3\li\chi\rho^n(r^n)^2|D_t^i(\frac{2p^n}{\rho^nr^n})|^2dx,\;
\sum_{i=0}^3\li\chi\rho^n(r^n)^2 |D_t^i(\frac{x}{(r^n)^2})|^2dx,\\
\sum_{j=0}^2\li\chi|D_t^j((r^n)^2D_xp^n)|^2dx;\;\;
\sum_{i=0}^3\int\zeta|\partial_t^ip^n|^2(r^n)^2dr^n,\;
\sum_{j=0}^2\int\zeta\frac{|\partial_t^j\partial_{r^n}p^n|^2}{\rho^n}(r^n)^2dr^n,\\
\sum_{i=0}^3\int\zeta|\partial_t^i(x\rho^n)|^2dr^n,\;
\sum_{j=0}^2\int\zeta\frac{|\partial_t^j(x\rho^n)|^2}{\rho^n(r^n)^2}dr^n,\;
\sum_{i=1}^3\int\zeta\frac{|\partial_t^i(\rho^nD_tr^n)|^2}{\rho^n}(r^n)^2dr^n,\\
\sum_{j=2}^3\int\zeta\frac{|\partial_t^j\rho^n|^2}{\rho^n}(r^n)^2dr^n,\;
(\int_0^{r_2-d}\zeta|\partial_{r^n}(\frac{\partial_t^2\rho^n}
{\sqrt{\rho^n}})|^2(r^n)^2dr^n)^{\frac{3}{2}} (\int_0^{r_2-d}\zeta
|\frac{\partial_t^2\rho^n} {\sqrt{\rho^n}}|^2(r^n)^2dr^n
)^{\frac{1}{2}},\\
\end{split}
\]
\[
\begin{split}
(\int_0^{r_2-d}\zeta|\partial_{r^n}(\frac{\partial_t^2(\rho^nD_tr^n)}
{\sqrt{\rho^n}})|^2(r^n)^2dr^n)^{\frac{3}{2}} (\int_0^{r_2-d}\zeta
|\frac{\partial_t^2(\rho^nD_tr^n)} {\sqrt{\rho^n}}|^2(r^n)^2dr^n
)^{\frac{1}{2}}.
\end{split}
\]\\

\noindent Claim 1. There exists a constant $C_4$ so that
$$\mathcal{G}^{n+1}(t)\leq C_4\{(1+M^{n+1}+|M^n|^2+\mathcal{F}^n(t))
\mathcal{H}^n(t)+(\mathcal{H}^n(t))^{2}\}
+C_5e^{C_6(M^n+M^{n-1})T} \mathcal{F}^n(t).$$\\

\noindent Proof of Claim 1: Recall (\ref{nbound}). And note that the
dynamics of $r^n$ follows from $\rho^n$:
\[
\begin{split}
&D_tr^n=-\frac{1}{(r^n)^2}\int_0^x\frac{D_t\rho^n}{(\rho^n)^2}dy;\\
&D_t^2r^n=-\frac{2|D_tr^n|^2}{r^n}
-\frac{1}{(r^n)^2}\int_0^x\frac{D_t^2\rho^n}{(\rho^n)^2}dy+
\frac{2}{(r^n)^2}\int_0^x\frac{(D_t\rho^n)^2}{(\rho^n)^3}dy;\\
&D_t^3r^n=-\frac{4D_t^2r^nD_tr^n}{r^n}+2r^n(\frac{D_tr^n}{r^n})^3+
\frac{2D_tr^n}{(r^n)^3}\int_0^x\frac{D_t\rho^n}{(\rho^n)^2}dy-
\frac{1}{(r^n)^2}\int_0^x\frac{D_t^3\rho^n}{(\rho^n)^2}dy\\
&\;\;\;\;\;\;\;\;\;\;\;\;
+\frac{6}{(r^n)^2}\int_0^x\frac{D_t^2\rho^nD_t\rho^n}{(\rho^n)^3}dy
-\frac{4D_tr^n}{(r^n)^3}\int_0^x\frac{(D_t\rho^n)^2}{(\rho^n)^3}dy
-\frac{6}{(r^n)^2}\int_0^x\frac{(D_t\rho^n)^3}{(\rho^n)^4}dy.
\end{split}
\]
Because of (\ref{nbound}) and
\[
\begin{split}
&|\int_0^x\frac{D_t^2\rho^n}{(\rho^n)^2}
dy|\leq(\int_0^x\frac{|D_t^2\rho^n|^2}{|\rho^n|^3}dy)^{\frac{1}{2}}
(\int_0^x \frac{1}{\rho^n}dy)^{\frac{1}{2}}\leq C
e^{(M^n+M^{n-1})\frac{T}{2}}
(\mathcal{F}^n)^{\frac{1}{2}}(r^n)^{\frac{3}{2}},\\
&|\int_0^x\frac{D_t^3\rho^n}{(\rho^n)^2} dy|\leq C
e^{(M^n+M^{n-1})\frac{T}{2}}(\mathcal{F}^n)^{\frac{1}{2}}(r^n)^{\frac{3}{2}},
\end{split}
\]
one gets the following:
\[
|D_t^2r^n|\leq C|M^n|^2r^n+Ce^{(M^n+M^{n-1})\frac{T}{2}}
\frac{(\mathcal{F}^n)^{\frac{1}{2}}}{(r^n)^{\frac{1}{2}}},\;\;|D_t^3r^n|\leq
C|D_t^2r^n|.
\]
Now let us look at Lagrangian terms. Here we provide the details for
the first and the last terms in the list. The first one deals with
pure $t$-derivative terms.
\[
\begin{split}
\bullet&\li\chi(\rho^n)^{2\gamma-1}dx\leq C_{in}e^{\gamma
M^nT}\li\chi(\rho^n)^{\gamma-1}dx\leq C_{in}e^{\gamma
M^nT}\mathcal{H}^n(t)\\
\bullet&\li\chi\frac{|D_t((r^n)^2(\rho^n)^{\gamma})|^2}{\rho^n(r^n)^4}dx\leq
2\li\chi
|\frac{D_tr^n}{r^n}|^2(\rho^n)^{2\gamma-1}dx+2\gamma^2\li\chi
|\frac{D_t\rho^n}{\rho^n}|^2 (\rho^n)^{2\gamma-1}dx\\
&\;\;\;\;\;\;\;\;\;\;\;\;\;\;\;\;\;\;\;\;\;\;
\;\;\;\;\;\;\;\;\;\;\;\;\;\;\;\;\;\;\;\leq C_{in}|M^n|^2e^{\gamma
M^nT}\mathcal{H}^n(t)\\
\bullet&\li\chi\frac{|D_t^2((r^n)^2(\rho^n)^{\gamma})|^2}{\rho^n(r^n)^4}dx\leq
C\li\chi\{|\frac{D_t^2r^n}{r^n}|^2+|\frac{D_tr^n}{r^n}|^2
|\frac{D_t\rho^n}{\rho^n}|^2+|\frac{D_t\rho^n}{\rho^n}|^4\}(\rho^n)^{2\gamma-1}dx\\
&\;\;\;\;\;\;\;\;\;\;\;\;\;\;\;\;\;\;\;\;\;\;
\;\;\;\;\;\;\;\;\;\;\;\;\;\;\;\;\;\;\;\;\;+C\li\chi\rho^{2\gamma}
\frac{|D_t^2\rho^n|^2}{(\rho^n)^3}dx\\
&\;\;\;\;\;\;\;\;\;\;\;\;\;\;\;\;\;\;\;\;\;\;
\;\;\;\;\;\;\;\;\;\;\;\;\;\;\;\;\;\;\;\leq
C_{in}(|M^n|^2+|M^n|^4+e^{(2M^n+M^{n-1})T}\frac{\mathcal{F}^n(t)}{(r_0)^3})
e^{\gamma M^nT}\mathcal{H}^n(t)\\
&\;\;\;\;\;\;\;\;\;\;\;\;\;\;\;\;\;\;\;\;\;\;
\;\;\;\;\;\;\;\;\;\;\;\;\;\;\;\;\;\;\;\;\;+C_{in}e^{(M^n+M^{n-1})T}e^{2\gamma
M^nT}\mathcal{F}^n(t)\\
\bullet&\li\chi\frac{|D_t^3((r^n)^2(\rho^n)^{\gamma})|^2}{\rho^n(r^n)^4}dx\leq
C_{in}(|M^n|^2+|M^n|^4+e^{(2M^n+M^{n-1})T}\frac{\mathcal{F}^n(t)}{(r_0)^3})
e^{\gamma M^nT}\mathcal{H}^n(t)\\
&\;\;\;\;\;\;\;\;\;\;\;\;\;\;\;\;\;\;\;\;\;\;
\;\;\;\;\;\;\;\;\;\;\;\;\;\;\;\;\;\;\;\;\;+C_{in}e^{(M^n+M^{n-1})T}e^{2\gamma
M^nT}\mathcal{F}^n(t)
\end{split}
\]
The last one treats $x$-derivative terms. Due to the dynamics of
$r^n$, we obtain the following:
\[
\begin{split}
\bullet&\li\chi
|(r^n)^2D_xp^n|^2dx=A^2\gamma^2\li\chi(\rho^n)^{2\gamma-2}
(r^n)^4|D_x\rho^n|^2 \leq C\mathcal{H}^n(t)\\
\bullet&\li\chi |D_t((r^n)^2D_xp^n)|^2dx\leq
C(1+|M^n|^2)\mathcal{H}^n(t)\\
\bullet&\li\chi |D_t^2((r^n)^2D_xp^n)|^2dx\leq
C(1+|M^n|^2+|M^n|^4+e^{(2M^n+M^{n-1})T}\frac{\mathcal{F}^n(t)}{(r_0)^3})
\mathcal{H}^n(t)
\end{split}
\]
For Eulerian terms, first note that
\[
\partial_{r^n}=\frac{\rho^n(r^n)^2}{\rho^{n-1}(r^{n-1})^2}\partial_{r^{n-1}},\;\;
\rho^n(r^n)^2dr^n=\rho^{n-1}(r^{n-1})^2dr^{n-1}.
\]
Hence, by change of variables, we obtain, for instance
\[
\begin{split}
\bullet&\int \zeta|p^n|^2(r^n)^2dr^n=A^2\gamma^2\int\zeta
(\rho^n)^{2\gamma-1}\rho^{n-1} (r^{n-1})^2dr^{n-1}\leq C_{in}
e^{((\gamma-1)M^n+M^{n-1})T} \mathcal{H}^n(t) \\
\bullet&\int\zeta\frac{|\partial_{r^n}p^n|^2}{\rho^n}(r^n)^2dr^n=A^2\gamma^2\int\zeta
\frac{(\rho^n)^{\gamma}}{\rho^{n-1}}|\frac{r^n}{r^{n-1}}|^4(\rho^n)
^{\gamma-2}|\partial_{r^{n-1}}\rho^n|^2(r^{n-1})^2dr^{n-1}\\
&\;\;\;\;\;\;\;\;\;\;\;\;\;\;\;\;\;\;\;\;\;\;\;\;\;\;\;\;\;\;\;\;\;\;\leq
C_{in}e^{(\frac{3\gamma+4}{3}M^n+\frac{7}{3}M^{n-1})T}\mathcal{H}^n(t)\\
\bullet&(\int_0^{r_2-d}\zeta|\partial_{r^n}(\frac{\partial_t^2\rho^n}
{\sqrt{\rho^n}})|^2(r^n)^2dr^n)^{\frac{3}{2}} (\int_0^{r_2-d}\zeta
|\frac{\partial_t^2\rho^n} {\sqrt{\rho^n}}|^2(r^n)^2dr^n
)^{\frac{1}{2}}\leq C_{in}e^{C(M^{n}+M^{n-1})T}(\mathcal{H}^n(t))^2
\end{split}
\]
Other $t$-derivative terms can be estimated in the same way. This finishes the
proof of Claim 1. \\

 Here are the terms, under-braced in the above estimates,
needed to be point-wisely estimated ($L^{\infty}$) to close the
above energy estimates:
\[
\begin{split}
&\text{Lagrangian terms}(x_0\leq x\leq 1)\Rightarrow
|\frac{D_t\rho^n}{\rho^n}|,\;|\frac{D_tr^n}{r^n}|,\;
|u^{n+1}|,\;|\rho^n(r^n)^2D_xu^{n+1}|,\;|\rho^n(r^n)^2D_{t}D_xu^{n+1}|\\
&\text{Eulerian terms}(0\leq r\leq r_2-d)\Rightarrow
|\frac{\partial_t\rho^n}{\rho^n}|,\;\rho^n|D_tr^n|^2,\;
\frac{|\partial_t(\rho^nD_tr^n)|^2}{\rho^n},
\;|u^{n+1}|,\;|\partial_{r^n}u^{n+1}|,\;|\partial_tu^{n+1}|
\end{split}
\]
As we noted in (\ref{nbound}), any $\rho^n$, $r^n$ related terms can
be bounded by $\mathcal{F}^n$ and $\mathcal{H}^n$. We need to show
that other $u^{n+1}$ related terms can be bounded by
$\mathcal{F}^{n+1}$ as well as $\mathcal{F}^n$, $\mathcal{H}^n$. Now
let us estimate $M^{n+1}$:
\[
\begin{split}
&\sup |\rho^n(r^n)^2D_xu^{n+1}+\frac{2u^{n+1}}{r^n}|\leq
\underline{\int_{x_0}^1\chi|\rho^n(r^n)^2D_xu^{n+1}+\frac{2u^{n+1}}{r^n}|dx}_{(i)}
\\
&\;\;\;\;
+\underline{\int_{x_0}^1\chi|D_x(\rho^n(r^n)^2D_xu^{n+1}+\frac{2u^{n+1}}{r^n})|dx}_{(ii)}
+(\int_0^{r_2-d}\zeta|\partial_{r^n}u^{n+1}+\frac{2u^{n+1}}{r^n}|^2(r^n)^2dr^n)
^{\frac{1}{2}}\\
&\;\;\;\;+(\int_0^{r_2-d}\zeta|\nabla_n(\partial_{r^n}u^{n+1}
+\frac{2u^{n+1}}{r^n})|^2(r^n)^2dr^n)
^{\frac{1}{2}}+(\int_0^{r_2-d}\zeta|\nabla_n^2(\partial_{r^n}u^{n+1}
+\frac{2u^{n+1}}{r^n})|^2(r^n)^2dr^n)^{\frac{1}{2}}
\end{split}
\]
\[
\begin{split}
(i)&\leq C(\int_{x_0}^1\chi\rho^ndx)^{\frac{1}{2}}
(\int_{x_0}^1\chi\{\rho^n(r^n)^4|D_xu^{n+1}|^2+\frac{2|u^{n+1}|^2}
{\rho^n(r^n)^2}\}dx)^{\frac{1}{2}}\\
(ii)&=\int_{x_0}^1\chi|\frac{D_tu^{n+1}}{(r^n)^2}+D_xp^n+\frac{4\pi
x}{(r^n)^4}|dx\\
&\leq
(\int_{x_0}^1\chi\frac{1}{(r^n)^4}dx)^{\frac{1}{2}}\{(\int_{x_0}^1\chi
|D_tu^{n+1}|^2dx)^{\frac{1}{2}}
+(\int_{x_0}^1\chi(r^n)^4|D_xp^n|^2dx)^{\frac{1}{2}}\}\\&\;\;\;\;+(\int_{x_0}^1
\chi\frac{1}{\rho^n(r^n)^4}dx)^{\frac{1}{2}}(\li\chi
\frac{16\pi^2x^2\rho^n}{(r^n)^4}dx)^{\frac{1}{2}}
\end{split}
\]
Hence,
\begin{equation}
\begin{split}
M^{n+1}\leq Ce^{\frac{M^n}{2}T}
(\li\chi\rho_{in}dx)^{\frac{1}{2}}(\mathcal{F}^{n+1})^{\frac{1}{2}}+
e^{\frac{2M^n}{3}T}(\li\chi\frac{1}{r_{in}^4}dx)^{\frac{1}{2}}
\{(\mathcal{F}^{n+1})^{\frac{1}{2}}
+(\mathcal{H}^n)^{\frac{1}{2}}\}\\
+e^{\frac{7M^n}{6}T}(\li\chi\frac{1}{\rho_{in}r_{in}^4}dx)^{\frac{1}{2}}
(\mathcal{H}^n)^{\frac{1}{2}} +Ce^{\frac{M^n}{2}T}\sup_{0\leq r\leq
r_2-d}|\frac{1}{\sqrt{\rho_{in}}}|(\mathcal{F}^{n+1})^{\frac{1}{2}}\label{Mn+1}
\end{split}
\end{equation}\\

\noindent Next we estimate $L^{\infty}$ bound of
$|\rho^n(r^n)^2D_xu^{n+1}|(=|\partial_{r^n}u^{n+1}|)$:
\[
\begin{split}
&\sup|\rho^n(r^n)^2D_xu^{n+1}|\leq\int_{x_0}^1\chi|\rho^n(r^n)^2D_xu^{n+1}|dx+
\int_{x_0}^1\chi|D_x(\rho^n(r^n)^2D_xu^{n+1})|dx\\
&\;\;\;\;\;\;\;+(\int_0^{r_2-d}\zeta|\partial_{r^n}u^{n+1}|^2(r^n)^2dr^n)^{\frac{1}{2}}
+(\int_0^{r_2-d}\zeta|\nabla_n\partial_{r^n}u^{n+1}|^2(r^n)^2dr^n)^{\frac{1}{2}}\\
&\;\;\;\;\;\;\;+(\int_0^{r_2-d}\zeta|\nabla_n^2\partial_{r^n}u^{n+1}|^2(r^n)^2dr^n)^{\frac{1}{2}}\\
 &\;\;\leq
(\int_{x_0}^1\chi\rho^ndx)^{\frac{1}{2}}(\int_{x_0}^1\chi\rho^n(r^n)^4
|D_xu^{n+1}|^2dx)^{\frac{1}{2}}+\int_{x_0}^1\chi|\frac{D_tu^{n+1}}{(r^n)^2}
+D_xp^n+\frac{4\pi
x}{(r^n)^4}-\underline{D_x(\frac{2u^{n+1}}{r^n})}_{\star}|dx\\
&\;\;\;\;\;\;\;+(\int_0^{r_2-d}\zeta|\partial_{r^n}u^{n+1}|^2(r^n)^2dr^n)^{\frac{1}{2}}
+(\int_0^{r_2-d}\zeta|\nabla_n\partial_{r^n}u^{n+1}|^2(r^n)^2dr^n)^{\frac{1}{2}}\\
&\;\;\;\;\;\;\;+(\int_0^{r_2-d}\zeta|\nabla_n^2\partial_{r^n}u^{n+1}|^2(r^n)^2dr^n)
^{\frac{1}{2}}
\end{split}
\]
For $\star$, first separate the integral into two by the quotient
rule, and use H$\ddot{o}$lder's inequality.
\[
\begin{split}
 &\star\leq\int_{x_0}^1\chi|\frac{D_xu^{n+1}}{r^n}|dx+
\int_{x_0}^1\chi|\frac{u^{n+1}}{\rho^n(r^n)^4}|dx\\
&\;\;\leq(\int_{x_0}^1\chi \rho^n(r^n)^4|D_xu^{n+1}|^2dx)
^{\frac{1}{2}}(\int_{x_0}^1\chi
\frac{1}{\rho^n(r^n)^6}dx)^{\frac{1}{2}}+(\int_{x_0}^1\chi\frac{|u^{n+1}|^2}
{\rho^n(r^n)^2} dx)
^{\frac{1}{2}}(\int_{x_0}^1\chi\frac{1}{\rho^n(r^n)^6}
dx)^{\frac{1}{2}}
\end{split}
\]
Hence, we get
\[
\begin{split}
|\rho^n(r^n)^2D_xu^{n+1}|\leq Ce^{\frac{M^n}{2}T}
(\li\chi\rho_{in}dx)^{\frac{1}{2}}(\mathcal{F}^{n+1})^{\frac{1}{2}}+
e^{\frac{2M^n}{3}T}(\li\chi\frac{1}{r_{in}^4}dx)^{\frac{1}{2}}
\{(\mathcal{F}^{n+1})^{\frac{1}{2}}
+(\mathcal{H}^n)^{\frac{1}{2}}\}\\
+e^{\frac{7M^n}{6}T}
(\li\chi\frac{1}{\rho_{in}r_{in}^4}dx)^{\frac{1}{2}}
(\mathcal{H}^n)^{\frac{1}{2}}+e^{\frac{3M^n}{2}T}(\li\chi\frac{1}{\rho_{in}r_{in}^6}dx)^{\frac{1}{2}}
(\mathcal{F}^{n+1})^{\frac{1}{2}}  \\
+Ce^{\frac{M^n}{2}T}\sup_{0\leq r\leq
r_2-d}|\frac{1}{\sqrt{\rho_{in}}}|(\mathcal{F}^{n+1})^{\frac{1}{2}}
\end{split}
\]

Similarly, $L^{\infty}$ bound of other terms can be obtained. The
estimate of $\mathcal{F}^{n+1}$ indeed can be closed by
$\mathcal{F}^{n+1},\mathcal{F}^{n}, \mathcal{H}^n$. It completes
(\ref{e1}).
\end{proof}\

Next we move onto $\mathcal{H}^{n+1}$.\\

\begin{proof}\textit{of Lemma \ref{energyii}:} Note that the
point-wise estimates of $\rho^{n+1}$ and $r^{n+1}$, as in
(\ref{nbound}), can be obtained by using $M^{n+1}$:
\[
\begin{split}
\rho_{in} e^{-M^{n+1}T}\leq \rho^{n+1}\leq \rho_{in} e^{M^{n+1}T};\;
r_{in}e^{-\frac{M^{n+1}T}{3}}\leq r^{n+1}\leq
r_{in}e^{\frac{M^{n+1}T}{3}};
\\ |\frac{D_t\rho^{n+1}}{\rho^{n+1}}|=|\rho^n(r^n)^2D_xu^{n+1}
+\frac{2u^{n+1}}{r^n}|\leq M^{n+1};\\
D_tr^{n+1}=-\frac{1}{(r^{n+1})^2}\int_0^x\frac{D_t\rho^{n+1}}{\rho^{n+1}}
\frac{1}{\rho^{n+1}}dy\;\Rightarrow\;
|\frac{D_tr^{n+1}}{r^{n+1}}|\leq \frac{M^{n+1}}{3}.
\end{split}
\]
Now we are ready to construct (\ref{e4}). Observe that we need all
the mixed derivatives estimates both in Lagrangian and in Eulerian
formulations because of overlapping terms, and note that the
overlapping region holds one-dimensional feature. Before we derive
Eulerian energy estimates, we complete $L^2$-type of Lagrangian
estimates, in particular of spatial derivatives of $\rho^{n+1}$. In
order to do so, we turn to the continuity equation (\ref{rho}).
First, here is the zeroth order estimate. Multiply (\ref{rho}) by
$(\rho^{n+1})^{\gamma-2}$ and integrate:
\[
\begin{split}
&\frac{1}{\gamma-1}\frac{d}{dt}\li\chi
(\rho^{n+1})^{\gamma-1}dx+\li\chi
(\rho^n(r^n)^2D_xu^{n+1}+\frac{2u^{n+1}}{r^n})(\rho^{n+1})^{\gamma-1}dx=0\\
\Rightarrow&\frac{1}{\gamma-1}\frac{d}{dt}\li\chi
(\rho^{n+1})^{\gamma-1}dx\leq M^{n+1}\li\chi
(\rho^{n+1})^{\gamma-1}dx
\end{split}
\]
To estimate $D_x\rho^{n+1}$, take $D_x$ of (\ref{rho}), multiply by
$(\rho^{n+1})^{\alpha}(r^n)^4D_x\rho^{n+1}$, and integrate:
\[
\begin{split}
\frac{1}{2}\frac{d}{dt}\int_{x_0}^1\chi(\rho^{n+1})^{\alpha}(r^n)^4|D_x\rho^{n+1}|^2dx
-\frac{\alpha+2}{2}\int_{x_0}^1\chi(\rho^{n+1})^{\alpha-1}D_t\rho^{n+1}
(r^n)^4|D_x\rho^{n+1}|^2dx\\
-2\li\chi(\rho^{n+1})^{\alpha}D_tr^n(r^n)^3|D_x\rho^{n+1}|^2dx\\
+\int_{x_0}^1\chi(\rho^{n+1})^{\alpha+1}(r^n)^4D_x(\rho^n(r^n)^2D_xu^{n+1}
+\frac{2u^{n+1}}{r^n})D_x\rho^{n+1}dx=0\\
\Rightarrow\frac{1}{2}\frac{d}{dt}\int_{x_0}^1\chi
(\rho^{n+1})^{\alpha}(r^n)^4|D_x\rho^{n+1}|^2dx \leq
\frac{(\alpha+2)M^{n+1}+2M^n+\epsilon}{2}\int_{x_0}^1\chi(\rho^{n+1})^{\alpha}
(r^n)^4|D_x\rho^{n+1}|^2dx\\+\frac{1}{2\epsilon}\int_{x_0}^1\chi
(\rho^{n+1})^{\alpha+2}|\underline{(r^n)^2D_x(\rho^n(r^n)^2D_xu^{n+1}
+\frac{2u^{n+1}}{r^n})}|^2dx\\
=\frac{1}{\mu}\{D_tu^{n+1}+(r^n)^2D_xp^n+\frac{4\pi x}{(r^n)^2}\}\\
\leq
C(M^{n+1}+M^n)\mathcal{H}^{n+1}(t)+C_{in}e^{CM^{n+1}T}(\mathcal{F}^{n+1}(t)
+\mathcal{H}^n(t))
\end{split}
\]
We may choose $\alpha=2\gamma-2$. $t$-derivatives do not destroy the
structure of the equation and thus estimates of
$D_{t}D_x\rho^{n+1},D_t^2D_{x}\rho^{n+1}$ follow in the similar
fashion with the same weight. Now we take one more spatial
derivative:
\[
\begin{split}
D_{t}D_x^2\rho^{n+1}=-\frac{1}{\mu}\{\frac{D_{t}D_xu^{n+1}}{(r^n)^2}-\frac{2D_tu^{n+1}}
{\rho^n(r^n)^5}+D_{x}^2p^n+\frac{4\pi}{(r^n)^4}-\frac{16\pi
x}{\rho^n(r^n)^7}\}\rho^{n+1}\\
-\frac{2}{\mu}\{\frac{D_tu^{n+1}}{(r^n)^2}+D_xp^n+\frac{4\pi
x}{(r^n)^4}\}D_x\rho^{n+1}-\{\rho^n(r^n)^2D_xu^{n+1}+\frac{2u^{n+1}}{r^n}\}
D_{x}^2\rho^{n+1}
\end{split}
\]
Multiply by $(\rho^{n+1})^{\alpha_1}(r^n)^8D_{x}^2\rho^{n+1}$ and
integrate:
\[
\begin{split}
\frac{1}{2}\frac{d}{dt}\int_{x_0}^1\chi(\rho^{n+1})^{\alpha_1}(r^n)^8
|D_{x}^2\rho^{n+1}|^2dx
=\frac{\alpha_1+2}{2}\int_{x_0}^1\chi(\rho^{n+1})^{\alpha_1-1}D_t\rho^{n+1}(r^n)^8
|D_{x}^2\rho^{n+1}|^2dx\\
+4\li\chi(\rho^{n+1})^{\alpha_1}D_tr^n(r^n)^7|D_{x}^2\rho^{n+1}|^2dx\\
-\frac{1}{\mu}\int_{x_0}^1\chi\{(r^n)^2D_{t}D_xu^{n+1}-\frac{2D_tu^{n+1}}
{\rho^n r^n}+(r^n)^4D_{x}^2p^n+4\pi-\frac{16\pi
x}{\rho^n(r^n)^3}\}(\rho^{n+1})^{\alpha_1+1}(r^n)^4D_{x}^2\rho^{n+1}dx\\
-\frac{2}{\mu}\int_{x_0}^1\chi\{D_tu^{n+1}+(r^n)^2D_xp^n+\frac{4\pi
x}{(r^n)^2}\}(r^n)^2D_x\rho^{n+1}(\rho^{n+1})^{\alpha_1}(r^n)^4D_{x}^2\rho^{n+1}
dx
\end{split}
\]
Let $$K^{n+1}\equiv\sup_{x_0\leq x\leq 1}
|(\rho^{n+1})^{\frac{\alpha_1}{2}}(r^n)^2D_x\rho^{n+1}|.$$ Apply the
Sobolev embedding theorem and the Cauchy-Schwarz inequality to get
\[
\begin{split}
K^{n+1}&\leq
\int_{x_0}^1|(\rho^{n+1})^{\frac{\alpha_1}{2}}(r^n)^2D_x\rho^{n+1}|dx+
\int_{x_0}^1|D_x((\rho^{n+1})^{\frac{\alpha_1}{2}}(r^n)^2D_x\rho^{n+1})|dx\\
&\leq
(\int_{x_0}^1(\rho^{n+1})^{\alpha_1}(r^n)^4|D_x\rho^{n+1}|^2dx)^{\frac{1}{2}}
+\int_{x_0}^1
(\rho^{n+1})^{\frac{\alpha_1}{2}-1}(r^n)^2|D_x\rho^{n+1}|^2dx\\
&\;\;\;+2\li\frac{(\rho^{n+1})^{\frac{\alpha_1}{2}}}{\rho^n
r^n}|D_x\rho^{n+1}| dx+(\int_{x_0}^1
(\rho^{n+1})^{\alpha_1}(r^n)^4|D_{x}^2\rho^{n+1}|^2dx)^{\frac{1}{2}}.
\end{split}
\]
Note that $|\frac{\rho^{n}}{\rho^{n+1}}|\text{ or
}|\frac{\rho^{n+1}}{\rho^{n}}|\leq e^{(M^{n+1}+M^n)T}$. Hence, we
may choose $\frac{\alpha_1}{2}-1=2\gamma-2$, and with this
$\alpha_1$, we get $$K^{n+1}\leq (1+C_{in}e^{\gamma
M^{n+1}T})(\mathcal{H}^{n+1}(t))^{\frac{1}{2}}$$ and also the
following:
\[
\begin{split}
\frac{1}{2}\frac{d}{dt}\int_{x_0}^1\chi(\rho^{n+1})^{\alpha_1}(r^n)^8
|D_{x}^2\rho^{n+1}|^2dx
\leq(\frac{\alpha_1+2}{2}M^{n+1}+2M^n+\epsilon)\int_{x_0}^1\chi(\rho^{n+1})^{\alpha_1}
(r^n)^8|D_{x}^2\rho^{n+1}|^2dx\\
+\frac{1}{2\mu\epsilon}\int_{x_0}^1\chi(\rho^{n+1})^{\alpha_1+2}
|(r^n)^2D_{t}D_xu^{n+1}-\frac{2D_tu^{n+1}}
{\rho^nr^n}+(r^n)^4D_{x}^2p^n+4\pi-\frac{16\pi
x}{\rho^n(r^n)^3}|^2dx\\
+\frac{4(K^{n+1})^2}{\mu\epsilon}\int_{x_0}^1\chi
|D_tu^{n+1}+(r^n)^2D_xp^n+\frac{4\pi x}{(r^n)^2}|^2dx\\
\leq
C(M^{n+1}+M^n)\mathcal{H}^{n+1}(t)+(C_{in}e^{CM^{n+1}T}+\mathcal{H}^{n+1}(t))
(\mathcal{F}^{n+1}(t) +\mathcal{H}^n(t))
\end{split}
\]
$D_{t}D_x^2\rho^{n+1}$ can be estimated in the same way with the
same weight. Take one more $D_x$ to get
\[
\begin{split}
D_{t}D_x^3\rho^{n+1}=-\frac{1}{\mu}\{\frac{D_{t}D_x^2u^{n+1}}{(r^n)^2}
-\frac{4D_{t}D_xu^{n+1}} {\rho^n(r^n)^5}-2D_tu^{n+1}D_x(\frac{1}
{\rho^n(r^n)^5})+D_{x}^3p^n\\+D_x(\frac{4\pi}{(r^n)^4}-\frac{16\pi
x}{\rho^n(r^n)^7})\}\rho^{n+1}-\frac{3}{\mu}\{\frac{D_{t}D_xu^{n+1}}{(r^n)^2}
-\frac{2D_tu^{n+1}}
{\rho^n(r^n)^5}+D_{x}^2p^n+\frac{4\pi}{(r^n)^4}-\frac{16\pi
x}{\rho^n(r^n)^7}\}D_x\rho^{n+1}\\
-\frac{3}{\mu}\{\frac{D_tu^{n+1}}{(r^n)^2}+D_xp^n+\frac{4\pi
x}{(r^n)^4}\}D_{x}^2\rho^{n+1}-\{\rho^n(r^n)^2D_xu^{n+1}+\frac{2u^{n+1}}{r^n}\}
D_{x}^3\rho^{n+1}
\end{split}
\]
Now consider weights $(\rho^{n+1})^{\alpha_2}(r^n)^{12}$. The
problematic term is \\$\int_{x_0}^1\chi
(\rho^{n+1})^{\alpha_2}(r^n)^{12}D_{x}^2p^nD_x\rho^{n+1}D_{x}^3\rho^{n+1}dx$:
\[
\begin{split}
&A\gamma\underbrace{\int_{x_0}^1\chi
(\rho^{n+1})^{\alpha_2}(\rho^n)^{\gamma-1}(r^n)^{12}D_{x}^2\rho^{n}D_x\rho^{n+1}
D_{x}^3\rho^{n+1}dx}_{(1)}\\
&\;\;\;+A\gamma(\gamma-1)\underbrace{\int_{x_0}^1\chi
(\rho^{n+1})^{\alpha_2}(\rho^n)^{\gamma-2}(r^n)^{12}(D_x\rho^n)^2D_x\rho^{n+1}
D_{x}^3\rho^{n+1}dx}_{(2)}
\end{split}
\]
We estimate $(1)$ and $(2)$ respectively as follows:
\[
\begin{split}
&(1)\leq
K^{n+1}\int_{x_0}^1\chi(\rho^{n+1})^{\frac{\alpha_2-\alpha_1}{2}}
(\rho^n)^{\gamma-1}(r^n)^4|D_{x}^2\rho^{n}|
(\rho^{n+1})^{\frac{\alpha_2}{2}}(r^n)^6|D_{x}^3\rho^{n+1}|dx\\
&\;\;\;\;\;\;\;\;\;\;\Rightarrow
\frac{\alpha_2-\alpha_1}{2}+\gamma-1\geq \frac{\alpha_1}{2}\\
&(2)\leq
K^{n+1}K^{n}\int_{x_0}^1\chi(\rho^{n+1})^{\frac{\alpha_2-\alpha_1}{2}}
(\rho^n)^{\gamma-2-\frac{\alpha_1}{2}}(r^n)^2|D_{x}\rho^{n}|
(\rho^{n+1})^{\frac{\alpha_2}{2}}(r^n)^6|D_{x}^3\rho^{n+1}|dx\\
&\;\;\;\;\;\;\;\;\;\;\Rightarrow
\frac{\alpha_2-\alpha_1}{2}+\gamma-2-\frac{\alpha_1}{2}\geq
\frac{\alpha}{2}
\end{split}
\]
Hence, we choose $\alpha_2$ as follows: $\alpha_2\geq
2\alpha_1+\alpha+4-2\gamma\geq
2(4\gamma-2)+2\gamma-2+4-2\gamma=8\gamma-2.$ On the other hand, in
Eulerian coordinates, we use the approximate equation (\ref{rhoe})
to get $L^2$-type of estimates. First, here is the zeroth order
estimate:
\[
\begin{split}
\frac{1}{2}\frac{d}{dt}\int\zeta(\rho^{n+1})^{\beta+2}(r^n)^2dr^n
-\frac{\beta}{2}\int\zeta\partial_t\rho^{n+1}(\rho^{n+1})^{\beta+1}(r^n)^2dr^n\\
+\underline{\int\zeta
D_tr^n(\rho^{n+1})^{\beta+1}\partial_{r^n}\rho^{n+1}(r^n)^2dr^n}_{\star}
+\int\zeta
(\rho^{n+1})^{\beta+2}\{\partial_{r^n}u^{n+1}+\frac{2u^{n+1}}{r^n}\}(r^n)^2dr^n=0
\end{split}
\]
For $\star$ :
\[
\begin{split}
(\beta+2)\cdot\star=-\int\partial_{r^n}\zeta
D_tr^n(\rho^{n+1})^{\beta+2}
(r^n)^2dr^n-\int\zeta\partial_{r^n}(D_tr^n)(\rho^{n+1})^{\beta+2}
(r^n)^2dr^n\\-\int\zeta D_tr^n(\rho^{n+1})^{\beta+2} 2r^ndr^n
\end{split}
\]
Note that $D_tr^n$ is more or less a zeroth order term in the
following sense:
\[
\begin{split}
\partial_{r^n}(D_tr^n)=\frac{2}{(r^n)^3}\int_0^x\frac{D_t\rho^n}{(\rho^n)^2}dy-
\frac{1}{(r^n)^2}\frac{D_t\rho^n}{(\rho^{n})^2}\rho^n(r^n)^2=-\frac{2D_tr^n}{r^n}-
\frac{D_t\rho^n}{\rho^n}\\
\partial_t(D_tr^n)=\frac{1}{(r^n)^2}\frac{D_t\rho^n}{(\rho^{n})^2}\rho^n(r^n)^2
D_tr^n=\frac{D_t\rho^n}{\rho^n}D_tr^n
\end{split}
\]
Let $\partial$ be either temporal or spatial derivative
($\partial_t\rho^{n+1}=D_t\rho^{n+1}-D_tr^n\partial_{r^n}\rho^{n+1}$),
\[
\begin{split}
\partial_t\partial\rho^{n+1}+\partial(D_tr^n)\partial_{r^n}\rho^{n+1}+
D_tr^n\partial_{r^n}\partial\rho^{n+1}+\partial\rho^{n+1}
\{\partial_{r^n}u^{n+1}+\frac{2u^{n+1}}{r^n}\}\\
+\rho^{n+1}\partial
\{\partial_{r^n}u^{n+1}+\frac{2u^{n+1}}{r^n}\}=0
\end{split}
\]
Multiply by $(\rho^{n+1})^{\beta}\partial\rho^{n+1}(r^n)^2$ and
integrate it to get
\[
\begin{split}
\frac{1}{2}\frac{d}{dt}\int\zeta
(\rho^{n+1})^{\beta}|\partial\rho^{n+1}|^2(r^n)^2dr^n-\frac{\beta}{2}\int\zeta
(\rho^{n+1})^{\beta-1}\partial_t\rho^{n+1}|\partial\rho^{n+1}|^2(r^n)^2dr^n
+\int\zeta(\rho^{n+1})^{\beta}\cdot\\
\partial(D_tr^n)\partial_{r^n}\rho^{n+1}
\partial\rho^{n+1}(r^n)^2dr^n+\underline{\int\zeta(\rho^{n+1})^{\beta+1}
\partial\rho^{n+1}D_tr^n\partial_{r^n}\partial\rho^{n+1}(r^n)^2dr^n}_{\star}
+\int\zeta(\rho^{n+1})^{\beta}\cdot\\
\{\partial_{r^n}u^{n+1}+\frac{2u^{n+1}}{r^n}\}
|\partial\rho^{n+1}|^2(r^n)^2dr^n+\int\zeta(\rho^{n+1})^{\beta+1}\partial
\{\partial_{r^n}u^{n+1}+\frac{2u^{n+1}}{r^n}\}\partial\rho^{n+1}(r^n)^2dr^n=0
\end{split}
\]
For $\star$, we get
\[
\begin{split}
2\cdot\star=-\int\partial_{r^n}\zeta(\rho^{n+1})^{\beta}D_tr^n|\partial\rho^{n+1}|^2(r^n)^2dr^n
-\beta\int\zeta(\rho^{n+1})^{\beta-1}\partial_{r^n}\rho^{n+1}
D_tr^n|\partial\rho^{n+1}|^2(r^n)^2dr^n\\
-\int\zeta(\rho^{n+1})^{\beta}\partial_{r^n}(D_tr^n)|\partial\rho^{n+1}|^2(r^n)^2dr^n
-\int\zeta(\rho^{n+1})^{\beta}D_tr^n|\partial\rho^{n+1}|^22r^ndr^n
\end{split}
\]
We may choose $\beta=\gamma-2$. Similarly, the energy estimates of
higher order derivatives of $\rho^{n+1}$ can be performed. Note that
$\partial^3\rho^{n+1}$, namely up to the third derivatives, needs to
be estimated in order to close the energy estimates. We may choose
the same weights $\beta=\gamma-2$. We refer Lemma \ref{weaving} for
closing the estimates.

Now, based on the above estimates, let us examine $\rho^{n+1}$
integrals with respect to $dr^{n+1}$, which we used at the previous
steps. By the definition of $r^n$ and $r^{n+1}$,
$\rho^n(r^n)^2dr^n=\rho^{n+1}(r^{n+1})^2dr^{n+1}$ and thus
$(r^{n+1})^2dr^{n+1}=\frac{\rho^{n}}{\rho^{n+1}}(r^n)^2dr^n$. Now we
want to show that $dr^{n+1}$ integrals can be controlled by $dr^n$
integrals derived as above. This can be done easily because
\[
\int\zeta(\rho^{n+1})^{2\gamma}(r^{n+1})^2dr^{n+1}=
\int\zeta(\rho^{n+1})^{2\gamma}\frac{\rho^{n}}{\rho^{n+1}}(r^n)^2dr^n.
\]
Note that $|\frac{\rho^{n}}{\rho^{n+1}}|\leq e^{(M^{n+1}+M^n)T}$.
 This completes the proof of Lemma
\ref{energyii}.
\end{proof}\

\begin{proof}\textit{of Theorem \ref{thm}:} By Proposition \ref{prop2}, now
we can take $n\longrightarrow \infty$ and we get limits
$u(t,x)=u(t,r)$, $\rho(t,x)=\rho(t,r)$, $r(t,x)$ as well as the
lagrangian transformation $x=\int_0^r\rho
s^2ds;\;r(t,x)=(3\int_0^x\frac{1}{\rho}dy)^{\frac{1}{3}};
\;D_tr(t,x)=u(t,x)$. Furthermore, those limits serve the solution to
the problem. Since $\mathcal{E}(t)\sim \lim_{n\rightarrow
\infty}\{\mathcal{F}^n(t)+\mathcal{H}^n(t)\}$, the energy bound is
easily obtained. Now it remains to show the uniqueness. Let
$(\rho_1,u_1,r_1)$ and $(\rho_2,u_2,r_2)$ be two strong solutions to
(\ref{nspL}) satisfying the same initial conditions. Note that it is
enough to show that $u_1=u_2$, since it right away implies $r_1=r_2$
from the dynamics of $r$: $D_tr=u$ and thus $\rho_1=\rho_2$ as well.
Now let us consider the momentum equations for $u_1$ and $u_2$ in
Lagrangian coordinates:
\[
\begin{split}
D_tu_1-\mu
D_x(\rho_1r_1^4D_xu_1)+\mu\frac{2u_1}{\rho_1r_1^2}=-r_1^2D_x
p_1-\frac{4\pi x}{r_1^2};\\
D_tu_2-\mu
D_x(\rho_2r_2^4D_xu_2)+\mu\frac{2u_2}{\rho_2r_2^2}=-r_2^2D_x
p_2-\frac{4\pi x}{r_2^2}.
\end{split}
\]
Subtracting one from another, we get the equation for $u_1-u_2$ as
follows:
\[
\begin{split}
D_t(u_1-u_2)-\mu
D_x(\rho_1r_1^4D_x(u_1-u_2))+\mu\frac{2(u_1-u_2)}{\rho_1r_1^2}=-r_1^2D_x
p_1-\frac{4\pi x}{r_1^2}+r_2^2D_x p_2+\frac{4\pi x}{r_2^2}\\+\mu
D_x((\rho_1r_1^4-\rho_2r_2^4)D_xu_2)-2\mu u_2(\frac{1}{\rho_1r_1^2}
-\frac{1}{\rho_2r_2^2})
\end{split}
\]
Multiply by $u_1-u_2$ and integrate to get
\[
\begin{split}
&\frac{1}{2}\frac{d}{dt}\int_0^1|u_1-u_2|^2 dx+\mu\int_0^1\rho_1
r_1|D_x(u_1-u_2)|^2+\frac{2|u_1-u_2|^2}{\rho_1r_1^2}dx\\
&=\int_0^1(r_1^2p_1-r_2^2p_2)D_x(u_1-u_2)dx +\int_0^1
(\frac{2p_1}{\rho_1r_1}-\frac{2p_2}{\rho_2r_2})(u_1-u_2)dx-4\pi\int_0^1
(\frac{x}{r_1^2}-\frac{x}{r_2^2})(u_1-u_2)dx\\
&\;\;-\mu\int_0^1(\rho_1r_1^4-\rho_2r_2^4)D_xu_2D_x(u_1-u_2)dx
-2\mu\int_0^1u_2(\frac{1}{\rho_1r_1^2} -\frac{1}{\rho_2r_2^2})
(u_1-u_2)dx\\
&\leq (\underbrace{\int_0^1
\frac{1}{\rho_1r_1^4}|r_1^2p_1-r_2^2p_2|^2
dx}_{(i)})^{\frac{1}{2}}(\int_0^1\rho_1r_1^4|D_x(u_1-u_2)|^2
dx)^{\frac{1}{2}}\\
&\;\;+\{(2\underbrace{\int_0^1
\rho_1r_1^2|\frac{p_1}{\rho_1r_1}-\frac{p_2}{\rho_2r_2}|^2
dx}_{(ii)})^{\frac{1}{2}}+(8\pi^2\underbrace{\int_0^1\rho_1r_1^2x^2
|\frac{1}{r_1^2}-\frac{1}{r_2^2}|^2
dx}_{(iii)})^{\frac{1}{2}}\}(\int_0^1\frac{2|u_1-u_2|^2}{\rho_1r_1^2}
dx)^{\frac{1}{2}}\\
&\;\;+\mu(\underbrace{\int_0^1
\frac{1}{\rho_1r_1^4}|(\rho_1r_1^4-\rho_2r_2^4)D_xu_2|^2
dx}_{(iv)})^{\frac{1}{2}}(\int_0^1\rho_1r_1^4|D_x(u_1-u_2)|^2
dx)^{\frac{1}{2}}\\
&\;\;+\mu(2\underbrace{\int_0^1\rho_1r_1^2u_2^2|\frac{1}{\rho_1r_1^2}
-\frac{1}{\rho_2r_2^2}|^2
dx}_{(v)})^{\frac{1}{2}}(\int_0^1\frac{2|u_1-u_2|^2}{\rho_1r_1^2}
dx)^{\frac{1}{2}}
\end{split}
\]
Now we would like to estimate $(i)-(v)$ in terms of $u_1-u_2$. In
order to do so, we use the explicit formula for $\rho_1$, $\rho_2$:
\[
\rho_1(t,x)=\rho_{in}(x)e^{-\int_0^t\rho_1r_1^2D_xu_1+\frac{2u_1}{r_1}
d\tau};\;\rho_2(t,x)=\rho_{in}(x)e^{-\int_0^t\rho_2r_2^2D_xu_2+\frac{2u_2}{r_2}
d\tau}
\]
Here we provide the detail for $(iv)$ and other terms can be
estimated in the same way. First $(iv)$ can be bounded by
\begin{equation}
(iv)\leq M^2 \int_0^1
\frac{1}{\rho_1r_1^4}|\frac{\rho_1r_1^4-\rho_2r_2^4}{\rho_2r_2^2}|^2
dx=M^2\int_0^1
\frac{1}{\rho_1}|\frac{\rho_1r_1^2}{\rho_2r_2^2}-\frac{r_2^2}{r_1^2}|^2
dx\label{vi}
\end{equation}
where $M$ is the bound of $\rho_2r_2^2D_xu_2$. Note that
\[
\begin{split}
|\frac{\rho_1r_1^2}{\rho_2r_2^2}-\frac{r_2^2}{r_1^2}|^2
&=|e^{\int_0^t\rho_2r_2^2D_xu_2-\rho_1 r_1^2D_xu_1 d\tau}-
e^{\int_0^t \frac{2u_2}{r_2}-\frac{2u_1}{r_1}d\tau}|^2\\
&\leq 2|\int_0^t\rho_2r_2^2D_xu_2-\rho_1 r_1^2D_xu_1 d\tau -\int_0^t
\frac{2u_2}{r_2}-\frac{2u_1}{r_1}d\tau |^2\\
&\leq 4|\int_0^t \rho_1r_1^2
D_x(u_1-u_2)d\tau|^2+4|\int_0^t(\rho_1r_1^2-\rho_2r_2^2) D_xu_2
d\tau|^2\\
&\;\;+4|\int_0^t\frac{2(u_1-u_2)}{r_1}d\tau|^2
+4|\int_0^t2u_2(\frac{1}{r_1}-\frac{1}{r_2})d\tau|^2\\
&\leq 4(\int_0^t
\rho_1d\tau)(\int_0^t\rho_1r_1^4|D_x(u_1-u_2)|^2d\tau)+4M^2|\int_0^t
\frac{\rho_1r_1^2}{\rho_2r_2^2}-1d\tau|^2\\
&\;\;+4(\int_0^t\rho_1d\tau)(\int_0^t\frac{4|u_1-u_2|^2}{\rho_1r_1^2}d\tau)
+4M^2|\int_0^t\frac{r_2}{r_1}-1 d\tau|^2
\end{split}
\]
Since
\begin{equation*}
\begin{split}
4M^2|\int_0^t \frac{\rho_1r_1^2}{\rho_2r_2^2}-1d\tau|^2\leq
8M^2t\int_0^t|\frac{\rho_1r_1^2}{\rho_2r_2^2}-\frac{r_2^2}{r_1^2}|^2d\tau+
8M^2t\int_0^t|\frac{r_2^2}{r_1^2}-1|^2 d\tau,\\
|\frac{r_2}{r_1}+1|\leq 1+e^{\frac{2M}{3}T}\;\;\; \text{ and }\;\;\;
4M^2|\int_0^t\frac{r_2}{r_1}-1 d\tau|^2\leq
4M^2t\int_0^t|\frac{r_2}{r_1}-1|^2 d\tau,
\end{split}
\end{equation*}
we get, for $0\leq t\leq T$ where $T$ is sufficiently small to be
fixed,
\[
\begin{split}
|\frac{\rho_1r_1^2}{\rho_2r_2^2}-\frac{r_2^2}{r_1^2}|^2&\leq
8M^2T\int_0^t|\frac{\rho_1r_1^2}{\rho_2r_2^2}-\frac{r_2^2}{r_1^2}|^2d\tau
+4M^2T(3+4e^{\frac{4M}{3}T})\int_0^t|\frac{r_2}{r_1}-1|^2 d\tau\\
&+4(\int_0^t \rho_1d\tau)(\int_0^t\rho_1r_1^4|D_x(u_1-u_2)|^2+
\frac{4|u_1-u_2|^2}{\rho_1r_1^2}d\tau).
\end{split}
\]
On the other hand, in the same way, one can derive the similar
inequality for $|\frac{r_2}{r_1}-1|^2$:
\[
|\frac{r_2}{r_1}-1|^2\leq 2|\int_0^t\frac{u_2}{r_2}-\frac{u_1}{r_1}
d\tau|^2\leq 2 M^2T\int_0^t|\frac{r_2}{r_1}-1|^2 d\tau+2
(\int_0^t\rho_1d\tau)(\int_0^t\frac{|u_1-u_2|^2}{\rho_1r_1^2}d\tau)
\]
By the Gronwall inequality, we get
\[
\begin{split}
|\frac{\rho_1r_1^2}{\rho_2r_2^2}-\frac{r_2^2}{r_1^2}|^2+|\frac{r_2}{r_1}-1|^2
\leq  9 (\int_0^T \rho_1d\tau)(\int_0^T\rho_1r_1^4|D_x(u_1-u_2)|^2+
\frac{2|u_1-u_2|^2}{\rho_1r_1^2}d\tau)\cdot\\
(1+2M^2T(7+8e^{\frac{4M}{3}T})t e^{2M^2T(7+8e^{\frac{4M}{3}T})t})\;
\end{split}
\]
Taking this into account (\ref{vi}), $(iv)$ can be controlled as
follows:
\[
\begin{split}
(iv)&\leq C_{M,T}\int_0^1 \frac{1}{\rho_1}(\int_0^T
\rho_1d\tau)(\int_0^T\rho_1r_1^4|D_x(u_1-u_2)|^2+
\frac{2|u_1-u_2|^2}{\rho_1r_1^2}d\tau)dx\\
&\leq C_{M,T} \int_0^T \int_0^1\rho_1r_1^4|D_x(u_1-u_2)|^2+
\frac{2|u_1-u_2|^2}{\rho_1r_1^2} dxd\tau
\end{split}
\]
where we have used $|\frac{1}{\rho_1}\int_0^T \rho_1d\tau|\leq
Te^{2MT}$. Following the same path, one can derive the following:
\[
(i)+(ii)+(iii)+(iv)+(v)\leq C_{M,T} \int_0^T
\int_0^1\rho_1r_1^4|D_x(u_1-u_2)|^2+
\frac{2|u_1-u_2|^2}{\rho_1r_1^2} dxd\tau
\]
Hence, we obtain
\[
\begin{split}
\frac{1}{2}\frac{d}{dt}\int_0^1|u_1-u_2|^2
dx+&\frac{\mu}{2}\int_0^1\rho_1
r_1|D_x(u_1-u_2)|^2+\frac{2|u_1-u_2|^2}{\rho_1r_1^2}dx\\
&\leq C_{M,T}\int_0^T \int_0^1\rho_1r_1^4|D_x(u_1-u_2)|^2+
\frac{2|u_1-u_2|^2}{\rho_1r_1^2} dxd\tau
\end{split}
\]
Now we integrate over $[0,t]$ to get  for $0\leq t\leq T$ where $T$
is sufficiently small,
\[
\{\frac{1}{2}\int_0^1|u_1-u_2|^2 dx\}(t)\leq
\{\frac{1}{2}\int_0^1|u_1-u_2|^2 dx\}(0),
\]
and therefore the uniqueness follows.
\end{proof}\



\noindent\textbf{Acknowledgments:} The author would like to thank
\textsc{Yan Guo} for many stimulating discussions. She also would
like to give many thanks to \textsc{Constantine Dafermos} and
\textsc{Chau-Hsing Su} for helpful
discussions. \\

\end{document}